\documentclass{article}
\usepackage{times}
\usepackage[hyperindex=true,pageanchor=true,hyperfigures=true,backref=false]{hyperref} 
\usepackage{amsmath}
\usepackage{amssymb}

\usepackage{url}
\usepackage{dsfont} 
\DeclareSymbolFontAlphabet{\Bbb}{AMSb}
\newlength{\myleftmargin}
\usepackage{amsmath}
\usepackage{amssymb}
\usepackage{amsthm}
\usepackage{color}
\DeclareSymbolFontAlphabet{\Bbb}{AMSb}

\newtheorem{theorem}{Theorem}[section]
\newtheorem{lemma}[theorem]{Lemma}
\newtheorem{proposition}[theorem]{Proposition}
\newtheorem{corollary}[theorem]{Corollary}
\newtheorem{definition}[theorem]{Definition}
\newtheorem*{assumptionsk}{Assumption K}
\newtheorem*{assumptionsc}{Assumption CK}
\newtheorem*{assumptionsx}{Assumption X}
\newtheorem*{assumptionsks}{Assumption KS}

\newtheorem{example}[theorem]{Example}

\newcommand{\atob}[2]{\emph{#1)} $\Rightarrow$ \emph{#2)}.} 
\newcommand{\aeqb}[2]{\emph{#1)} $\Leftrightarrow$ \emph{#2)}.} 
\newcommand{\ada}[1]{\emph{#1).}}

\newenvironment{proofof}[1]{\noindent{\bf Proof of #1:}}{\qed\medskip}

\newlength{\fixboxwidth}
\definecolor{darkgreen}{rgb}{0,0.6,0}

\newcommand{\ca}[1]{{\cal #1}}

\newcommand{\mycdot}{\,\cdot\,}

\newcommand{\N}{\mathbb{N}}

\newcommand{\R}{\mathbb{R}}
\newcommand{\Rd}{\mathbb{R}^d}

\newcommand{\E}{\mathbb{E}}
\DeclareMathOperator{\var}{Var}
\newcommand{\dx}[1]{\hspace*{0.25ex}d\hspace*{-0.15ex}#1}
\newcommand{\dxy}[2]{\dx{#1}(#2)}
\newcommand{\dP}[1]{\dxy P{#1}}

\newcommand{\dn}[1]{\dxy \nu{#1}}

\renewcommand{\a}{\alpha}
\renewcommand{\b}{\beta}
\newcommand{\g}{\gamma}

\renewcommand{\d}{\delta}

\newcommand{\e}{\varepsilon}
\newcommand{\eps}{\epsilon}

\renewcommand{\k}{\kappa}
\newcommand{\lb}{\lambda}
\newcommand{\x}{\xi}

\newcommand{\s}{\sigma}

\renewcommand{\t}{\tau}

\renewcommand{\P}{\Phi}
\newcommand{\om}{\omega}
\newcommand{\Om}{\Omega}

\DeclareMathOperator{\spann}{span}

\DeclareMathOperator{\ran}{ran}
\DeclareMathOperator{\rank}{rank}
\DeclareMathOperator{\interior}{int}

\DeclareMathOperator{\id}{id}
\newcommand{\eins}{\boldsymbol{1}}

\newcommand{\snorm}[1]{\Vert #1 \Vert}
\newcommand{\mnorm}[1]{\bigl\Vert \, #1 \, \bigr\Vert}
\newcommand{\bnorm}[1]{\Bigl\Vert \, #1 \, \Bigr\Vert}

\newcommand{\inorm}[1]{\Vert #1 \Vert_\infty}

\newcommand{\sLx}[2]{{\ca L_{#1}(#2)}}
\newcommand{\sLo}[1]{\sLx 1 {#1}}
\newcommand{\sLt}[1]{\sLx 2 {#1}}
\newcommand{\sLtn}{\sLx 2 \nu}

\newcommand{\Lx}[2]{{L_{#1}(#2)}}
\newcommand{\Lt}[1]{\Lx 2 {#1}}
\newcommand{\Ltn}{\Lx 2 \nu}

\newcommand{\seqdia}[7]{
\setlength{\unitlength}{1ex}
\begin{picture}(70,10)
\put(-4,4){\makebox(9,2)[rt]{$#1$}}
\put(17,4){\makebox(9,2)[rt]{$#2$}}
\put(40,4){\makebox(9,2)[lt]{$#3$}}
\put(61,4){\makebox(9,2)[lt]{$#4$}}
\put(6,5){\vector(1,0){12}}
\put(27,5){\vector(1,0){12}}
\put(60,5){\vector(-1,0){12}}
\put(9,6){\makebox(6,2){$#5$}}
\put(30,6){\makebox(6,2){$#6$}}
\put(52,6){\makebox(6,2){$#7$}}
\end{picture}}

\newcommand{\tridia}[6]{
\setlength{\unitlength}{1ex}
\begin{picture}(60,20)
\put(20,14){\makebox(9,2)[rt]{$#1$}}
\put(51,14){\makebox(9,2)[lt]{$#2$}}
\put(35.5,2){\makebox(9,2){$#3$}}
\put(30,15){\vector(1,0){20}}
\put(30,14){\vector(1,-1){9}}
\put(41,5){\vector(1,1){9}}
\put(37,16){\makebox(6,2){$#4$}}
\put(27.5,8){\makebox(6,2)[rt]{$#5$}}
\put(47.5,8){\makebox(6,2)[lt]{$#6$}}
\end{picture}}

\newcommand{\quadiasn}[8]{
\setlength{\unitlength}{1ex}
\begin{picture}(60,20)
\put(20,15){\makebox(9,2)[rt]{$#1$}}
\put(45,15){\makebox(9,2)[lt]{$#2$}}
\put(20,2){\makebox(9,2)[rt]{$#3$}}
\put(45,2){\makebox(9,2)[lt]{$#4$}}
\put(30,16){\vector(1,0){14}}
\put(28,14){\vector(0,-1){9}}
\put(46,5){\vector(0,1){9}}
\put(30,3){\vector(1,0){14}}
\put(34,17){\makebox(6,2){$#5$}}
\put(21,8.5){\makebox(6,2)[rt]{$#6$}}
\put(47,8.5){\makebox(6,2)[lt]{$#7$}}
\put(34,0){\makebox(6,2){$#8$}}
\end{picture}}

\newcommand{\eHnx}[1] {{[H]_\sim^{#1}}}
\newcommand{\ebbHnx}[1] {{[\bar H]_\sim^{#1}}}
\newcommand{\eHnb} {{\eHnx \b}}

\newcommand{\Hxy}[2]{{H_{#1}^{#2}}}
\newcommand{\bHxy}[2]{{\hat H_{#1}^{#2}}}
\newcommand{\bbHxy}[2]{{\bar H_{#1}^{#2}}}
\newcommand{\kxy}[2]{{k_{#1}^{#2}}}
\newcommand{\bkxy}[2]{{\hat k_{#1}^{#2}}}

\newcommand{\HSb}{\Hxy S\b}
\newcommand{\bHSb}{\bHxy S\b}
\newcommand{\HTb}{\Hxy T\b}

\newcommand{\kSb}{\kxy S\b}
\newcommand{\bkSb}{\bkxy S\b}
\newcommand{\kTb}{\kxy T\b}

\title{Convergence Types and Rates  in Generic Karhunen-Lo{\`e}ve Expansions with Applications 
to Sample Path Properties}

\author{Ingo Steinwart\\
Institute for Stochastics and Applications\\
Faculty 8: Mathematics and Physics\\
University of Stuttgart\\
D-70569 Stuttgart Germany \\
\texttt{\small ingo.steinwart@mathematik.uni-stuttgart.de}
}

\begin{document}

\maketitle

\begin{abstract}
 We establish a Karhunen-Lo{\`e}ve expansion for generic centered, second order stochastic processes,
 which does not rely on topological assumptions.
 We further 
 investigate in which norms the expansion converges and derive exact average rates of convergence for these norms.
 For 
Gaussian processes we additionally prove certain sharpness results
in terms of the norm. 
Moreover, we investigate when the generic Karhunen-Lo{\`e}ve expansion 
can be used to construct reproducing kernel Hilbert spaces (RKHSs)
 containing the paths of a version of the process. 
 We further illustrate how the general theory can be applied, even in the absence of
 an explicitly known Karhunen-Lo{\`e}ve expansion,
 by comparing the
smoothness of the paths with the smoothness of the functions
contained in the RKHS of the covariance function and 
by discussing some small ball probabilities.
Key tools for our results are a recently shown generalization
of Mercer's theorem, spectral properties of the covariance
integral operator,  interpolation spaces of the real method, and for the smoothness results, entropy numbers of 
embeddings between classical function spaces.
\end{abstract}

\section{Introduction}\label{sec:intro}

Given a  real-valued, centered stochastic process $(X_t)_{t\in T}$ with finite second moments,
the covariance function $k:T\times T\to \R$ defined by 
$k(s,t) := \E X_sX_t$ is positive semi-definite.
Consequently, there exists a reproducing kernel Hilbert space (RKHS) $H$ on $T$ for which
$k$ is the (reproducing) kernel.
It is well-known that there are intimate relationships between 
$H$ and the stochastic process.

One such relation   is described by the 
classical Lo{\`e}ve isometry $\Psi:\Lx 2 X\to H$ defined by 
$\Psi(X_t)= k(t,\cdot)$, where 
$\Lx 2 X$ denotes the $\Lx 2 P$-closure of the space spanned by $(X_t)_{t\in T}$.
In particular, if    $(e_i)_{i\in I}$ is 
an arbitrary orthonormal basis (ONB) of $H$, then the process enjoys the  representation, 
\begin{equation}\label{Loeve-iso-representation}
   X_t = \sum_{i\in I} \xi_i e_i(t)\, ,
\end{equation}
where $(\xi_i)_{i\in I}$ is the family of uncorrelated random variables
given by $\xi_i:= \Psi^{-1}(e_i)$, and the 
convergence is, for each $t\in T$, unconditional in $\Lx 2 P$, that is, independent of the order of summation.

Not surprisingly, 
 the relationship between the process and its RKHS is, for Gaussian processes,
even closer, as the finite dimensional distributions of the process are completely 
determined by $k$. For example, the isometry $\Psi$
can be used  to define stochastic integrals,
see e.g.~\cite[Chapter 7]{Janson97}, and if $H$ is separable,
the representation \eqref{Loeve-iso-representation} converges also $P$-almost surely for each $t$
and $(\xi_i)_{i\in I}$ is a family of independent, standard normal random variables, see 
e.g.~\cite[Theorem 8.22]{Janson97}.
Last but not least, if $T$ is a compact metric space and $(X_t)_{t\in T}$ has continuous paths,
we even have $P$-almost surely   uniform convergence in $t$, 
see \cite[Theorem 3.8]{Adler90}. Note that unlike the convergence in \eqref{Loeve-iso-representation},
uniform convergence in $t$ makes it possible to represent the \emph{paths} of the process by a series expansion.

For compact and metric $T$ it is also possible to obtain path representations for more general processes.
Indeed,
%
%
if $\nu$ is a strictly positive and finite Borel measure on $T$
and $k$ is continuous,
the famous Karhunen-Lo{\`e}ve expansion
shows that there is
an ONB $(e_i)_{i\in I}$ of $H$ that is 
also orthogonal in $\Ltn$ and, additionally to \eqref{Loeve-iso-representation}, we have
\begin{equation}\label{KL-expansion-intro}
   X(\om) = \sum_{i\in I}  \xi_i(\om) e_i \, ,
\end{equation}
where the series converges unconditionally in $\Ltn$ 
for $P$-almost all $\om\in \Om$. In addition, we  
have $X = \sum_{i\in I} \xi_i\otimes e_i$ with unconditional convergence in $\Lx 2 {P\otimes \nu}$.
The form of convergence in \eqref{KL-expansion-intro} again allows for 
a series expansion of the paths of the process, but
unlike for the above mentioned Gaussian processes, only with $\Ltn$-convergence.
Unfortunately, however, the assumptions needed for \eqref{KL-expansion-intro} are significantly more
restrictive than those for \eqref{Loeve-iso-representation}, and thus a natural question 
is to ask for weaker assumptions ensuring a path representation \eqref{KL-expansion-intro}.
In addition, $\Ltn$-convergence is a rather weak form of convergence so that is seems to be desirable to 
replace it by stronger notions of convergence such as uniform convergence in $t$.

Another, rather different relationship between the process and its RKHS is in terms of quadratic mean smoothness.
For example, 
 if $T$ is a metric space, then the process is continuous in quadratic mean, if 
and only if its kernel $k$ is continuous. Moreover, a similar statement is true for quadratic mean 
differentiability. We refer to  \cite[p.~63]{BeTA04} and \cite[p.~65ff]{Yaglom87} for details.
Of course, smoothness in quadratic mean is not related to the smoothness of the
paths of the process. However, considering the path expansion \eqref{KL-expansion-intro}
it seems natural to ask to which extend the paths inherit smoothness properties 
from $H$, or from the ONB $(e_i)_{i\in I}$.
Probably, the first attempt in this direction is to check whether the paths 
are $P$-almost surely contained in $H$. Unfortunately, this is, in general not true. Indeed,
for Gaussian processes with infinite dimensional RKHS the paths are $P$-almost surely \emph{not}
contained in $H$, see \cite[Corollary 7.1]{LuBe01a} and also \cite{Parzen63a}. A natural
and well-studied next question is to 
look for larger Banach spaces $E$ that do contain the paths almost surely. 
For example, continuity and boundedness of the paths can be easily described by 
suitable spaces $E$. 
In view of the considered relationship between path properties and $H$ 
one may also
ask
for larger RKHSs $\bar H$ that do contain the paths almost surely. The first result in this
direction goes back to Driscoll, see \cite{Driscoll73a}. Namely, he essentially showed:

\begin{theorem}\label{driscoll}
   Let $(T,d)$ be a separable metric space and  $(X_t)_{t\in T}$ be a centered and 
   continuous Gaussian process, whose kernel $k$ is continuous. Then for all RKHS $\bar H$ on $T$ having a continuous 
   kernel, the following statements are equivalent:
   \begin{enumerate}
      \item Almost all paths of the process are contained in $\bar H$.
      \item We have $H\subset \bar H$ and the embedding $\id: H\to \bar H$ is Hilbert-Schmidt.
   \end{enumerate}
\end{theorem}

Since being Hilbert-Schmidt is  a rather strong notion of compactness, Driscoll's theorem shows that 
suitable spaces $\bar H$ need to be significantly larger than $H$, at least for Gaussian processes
satisfying the assumptions above. In particular, if we try to  describe smoothness properties 
of the paths by suitable Sobolev spaces $\bar H$, this result suggests that the paths should be 
rougher than the functions in $H$.
More recently, Luki{\'c} and   Beder have shown, see \cite[Theorem 5.1]{LuBe01a},
that for \emph{arbitrary} centered, second-order stochastic process  $(X_t)_{t\in T}$ condition \emph{ii)}
implies the existence of a version  $(Y_t)_{t\in T}$ whose paths are almost surely contained in $\bar H$,
and for generic Gaussian processes \cite[Corollary 7.1]{LuBe01a} shows \emph{i)} $\Rightarrow$ \emph{ii)}.
Furthermore, they provide examples of non-Gaussian processes, for which the implication 
\emph{i)} $\Rightarrow$ \emph{ii)} does \emph{not} hold, and  they also present modifications 
\emph{i')} and \emph{ii')} of \emph{i)} and \emph{ii)},
for which we have \emph{i')} $\Rightarrow$ \emph{ii'')} in the general case, see 
\cite[Theorem 3.1 and Corollary 3.1]{LuBe01a} for details.
Summarizing these results, it seems fair to say that we already have reasonably good means to \emph{test}
whether a given RKHS $\bar H$ contains the paths of our process almost surely.
Except for a couple of specific examples, however, 
very little is known 
how 
 to \emph{construct} such an $\bar H$, or even
whether such an $\bar H$ \emph{exists} 
cf.~\cite[p.~255ff]{Lukic04a}.

It turns out in this paper, that all the questions raised above are related to each other by a rather general 
form of Mercer's theorem and its consequences, which have been recently presented in \cite{StSc12a}.
Before we go into details in the next sections let us briefly outline our main results. To this end
let us assume in the following that we have a $\s$-finite measure $\nu$ on $T$ and a centered, second order process
   $(X_t)_{t\in T}$ with $X\in \sLx 2 {P\otimes \nu}$. It turns out that for such processes, $H$
   is ``contained'' in $\Ltn$ and 
the ``embedding'' $H \to \Ltn$ is
   Hilbert-Schmidt, which makes the results from  \cite{StSc12a} readily applicable.
   Here we use the quotation marks, since we actually need to consider equivalence classes 
   to properly define the embedding. As a matter of fact, the entire theory of \cite{StSc12a} foots on the 
	careful differentiation between functions and their equivalence classes, and thus we need
	to adopt 
   the somewhat pedantic notation of \cite{StSc12a} later in the paper. For now, however, let us ignore these   differences
   for the informal description of our main results:
\begin{itemize}
   \item The  Karhunen-Lo{\`e}ve expansion
   \eqref{KL-expansion-intro} always holds for the process $(X_t)_{t\in T}$, and in this case 
   $(\xi_i)_{i\in I}$ is an orthonormal system (ONS) in $\Lx 2 P$, $(e_i)_{i\in I}$ is an ONS in $H$,  
   and $(\mu_i^{-1/2}e_i)_{i\in I}$ is an ONS in $\Ltn$. Here $(\mu_i)_{i\in I}$ denotes the eigenvalue sequence 
   of the integral operator $T_k:\Ltn\to \Ltn$ associated to $k$, which turns out to be summable.
   Conversely, every such triple $(\xi_i)_{i\in I}$,  $(e_i)_{i\in I}$, and $(\mu_i)_{i\in I}$ gives 
   a  centered, second order process
   $(X_t)_{t\in T}$ with $X\in \sLx 2 {P\otimes \nu}$ via \eqref{KL-expansion-intro}  and this representation is unique.

   \item If the embedding $H \to \Ltn$ is, in a certain sense, more compact than Hilbert-Schmidt, then 
   almost all paths of the process are contained in  suitable interpolation spaces between   $\Ltn$ and $H$.
   Moreover, \eqref{KL-expansion-intro} converges in these interpolation spaces, too, and the average rate 
   of this convergence can be exactly described by the tail behavior of the eigenvalue sequence.
   Finally, for Gaussian processes the results are sharp in the considered scale of interpolation spaces.
   
   \item Under even stronger compactness assumptions on the   embedding $H \to \Ltn$, some of the 
   interpolation spaces are RKHSs and there exists a version of the process 
   having almost all its paths in these RKHSs. 
   
   \item Using the eigenvalue sequence of $T_k$, small ball probabilities of Gaussian processes
   with respect 
   to the above mentioned interpolation spaces can be estimated by extending known techniques.
   
   \item If $T\subset \Rd$ is a bounded and open subset with suitable boundary conditions,
   and $H$ is embedded into a (fractional) Sobolev space $W^m(T)$ with $m>d/2$, then almost all 
   paths are in the fractional Sobolev space $W^{m-d/2-\eps}(T)$, where   $\e>0$ is arbitrary. 
   In other words, the paths of $X$ are about $d/2$-less smooth than 
   the functions in $H$. 
   Again, for 
   Gaussian processes this  turns out to be sharp and small ball probabilities can be estimated.
\end{itemize}

  Describing path properties of a process in terms of spaces is not only a stochastic question
  in its own interest, but also important for other areas.
  For example, 
  certain non-parametric Bayesian methods for regression problems, called ``Gaussian processes'',
  use Gaussian processes as a prior, see \cite{RaWi06}. Understanding small spaces that 
  contain all paths of the prior process is then important for the mathematical analysis,
  as these spaces determine both the approximation properties of the non-parametric method as well
  as its statistical properties, see e.g.~\cite{PiWuLiMuWo07a,VaZa08a,VaZa11a}.
  In this regard note that one of the strengths of these methods is that they can be considered 
  on general input spaces, which translates into general index sets $T$ in our terminology.
  Similarly, certain spatial
   statistical methods require knowledge on the paths properties in terms of spaces, we refer to
  \cite{Stein99,GnKlSc10a,ScScSc13a} and the references in these articles.

  Karhunen-Lo{\`e}ve expansions, their speed of convergence, and 
   their relation to path properties
  have recently been
  considered in  the area of numerics for stochastic partial differential equations, too.
   Without going into details we refer to  
   \cite{ScTo06a,HoSc13a,LoPoSh14,LaScXXa} and the references therein.
  Moveover, \cite{LiPaWo12a,LiPaWo12b,Xu14a} consider 
  Karhunen-Lo{\`e}ve expansions to investigate how well processes can be approximated 
  by linear schemes.
  In all these papers, the eigenvalue behavior of the operator $T_k$ is crucial
  to estimate speeds of convergence. Similarly, 
    $L_2$-small ball probabilities for Gaussian processes
    can be described by the eigenvalue asymptotics,
    see 
  e.g.~the survey \cite{LiSh01a} and the references mentioned therein.
  In this respect recall that the eigenvalue behavior may be known even if the exact
   eigenvalues and eigenfunctions are unknown. For example, 
  for the fractional Brownian motion asymptotics are 
  determined in \cite{Bronski03a}, while eigenvalue estimates for certain integrated processes
   can be found in \cite{GaHaTo03a,NaNi04a}. Moreover, 
  eigenvalue estimates for covariance functions of tensor type are derived  in \cite{KaNaNi08a,KaNa14a}.
  All these papers also apply their eigenvalue estimates to small ball probabilities.

Classical examples of explicit  Karhunen-Lo{\`e}ve expansions include those 
of the Wiener process, the Brownian bridge, and the Ornstein-Uhlenbeck process.
Finding such
explicit  expansions requires to solve an eigenvalue problem associated to the integral operator 
 of the covariance function $k$, which, in general, can be viewed as a difficult problem.
 This may be the reason, why so far only a few 
 of such 
 explicit expansions are known.
 Recently, however, this question has regained attraction. For example, 
 \cite{Deheuvels06a} derive an explicit expansion for mean-centered Wiener processes in terms of Bessel functions,
 and \cite{Deheuvels07a} extends these considerations to a mean-centered Brownian bridge.
 Multivariate versions of these results are given in \cite{DePeYo06a}.
 Weighted and unweighted Karhunen-Lo{\`e}ve expansions of so-called $\a$-Wiener bridges 
  have been established in \cite{BaIg11a}, while \cite{Liu13a} considers multi-dimensional sums of 
  independent Wiener processes and Brownian bridges.
  Further examples of recently obtained explicit expansions can be found in \cite{Pycke01a,Istas06a,Pycke07a,DeMa08a,AiLiLi12a,LiHuMa14a}.

The rest of this paper is organized as follows: In Section \ref{sec:prelim} some concepts from \cite{StSc12a}
are recalled and some additional results are presented. The generic  Karhunen-Lo{\`e}ve expansion
is established in Section \ref{sec:KL-generic} and Section \ref{sec:path-in-interpol} contains
the results that are related to stronger notions of convergence in  the  Karhunen-Lo{\`e}ve expansion.
In Section \ref{sec:paths-in-rkhs} we continue these investigations with the focus on instances, where
the interpolation spaces are RKHSs. The Sobolev space related results are presented as applications of
the general theory in  Sections \ref{sec:path-in-interpol} and \ref{sec:paths-in-rkhs}, while 
the small ball probabilities can be found at the end of Section \ref{sec:path-in-interpol}.
Section \ref{sec:final-remaks} contains a few final remarks.
All proofs as well as some auxiliary results can be found in Section \ref{sec:proofs}.

\section{Preliminaries}\label{sec:prelim}

Let us begin by introducing some notations used throughout this paper. To this end, let 
 $(T,\ca B, \nu)$ be a measure space. 
 As usual, $\ca B$ 
is called $\nu$-complete, 
if, for 
every $A\subset T$ for which there exists an $N\in \ca B$ 
such that $A\subset N$ and $\nu(N)=0$, we have  $A\in \ca B$. In 
this case we say that
$(T,\ca B, \nu)$ is complete. 

For $S\subset T$ we denote the indicator function of $S$ by $\eins_S$.
Moreover, for an
$f:S\to \R$ we denote its 
zero-extension by $\hat f$, that is, $\hat f(t) := f(t)$ for all $t\in S$ and $\hat f(t) := 0$
otherwise.

As usual, $\sLx 2 \nu$ denotes the set of all measurable functions $f:T\to \R$ such that 
$
\snorm f_{\sLx 2 \nu} :=\int |f|^2 \, d\nu <\infty
$. 
For $f\in \sLx 2 \nu$, we further write 
$$
[f]_\sim:= \bigl\{ g\in \sLx 2\nu: \nu(\{f\neq g\})=0\bigr\}
$$
for the $\nu$-equivalence class of $f$.
Let $\Ltn := \sLx 2 \nu_{/\sim}$ be the corresponding quotient space
 and $\snorm \cdot_\Ltn$ be its norm.
For an arbitrary, non-empty index set $I$ and $p\in (0,\infty)$, 
we denote,  the space 
of all $p$-summable real-valued families by
$\ell_p(I)$.

Given two non-negative sequences $(a_i)_{i\geq 1}$ and $(b_i)_{i\geq 1}$
 we write 
 $a_i\preceq b_i$, if there exists a constant $c\in (0,\infty)$ such that 
 $a_i \leq c b_i$ for all $i\geq 1$. Moreover, we write $a_i \asymp b_i$, if we have both 
$a_i\preceq b_i$ and $b_i\preceq a_i$. 
Finally, we write $a_i\sim b_i$, if $\lim_{i\to \infty} a_i/b_i = 1$.

%

In the following,  we say that a
 Banach space $F$ is continuously embedded into a Banach space $E$,
if $F\subset E$ and the identity map $\id :F\to E$ is continuous. In this case, we sometimes write $F\hookrightarrow E$.

Let us now recall some properties of reproducing kernel Hilbert spaces (RKHSs),
and their 
interaction with measures from \cite{StSc12a}.
To this end, let $(T,\ca B, \nu)$ be a measure space
and $k:T\times T\to \R$ be a measurable (reproducing) kernel with RKHS $H$,
see e.g.~\cite{BeTA04,Wendland05,StCh08}
for more information about these spaces. 
Recall that in this case the RKHS $H$ 
consists of measurable functions $T\to \R$. In the following, we 
say that  $H$ is embedded into $\Lt \nu$, if all $f\in H$ are measurable with $[f]_\sim \in \Lt\nu$  and 
the linear operator  
\begin{eqnarray*}
	I_k: H & \to & \Lt \nu\\
	f&\mapsto & [f]_\sim
\end{eqnarray*}
is  continuous. We write $[H]_\sim$ for its image, that is $[H]_\sim := \{[f]_\sim: f\in H\}$.
Moreover, we say that $H$ is compactly embedded into $\Lt \nu$, if $I_k$ is compact.
For us, the most interesting class of compactly embedded RKHSs $H$
are those whose kernel $k$ satisfies
\begin{equation}\label{kernel-diag-int}
\snorm k_{\sLtn}:=  \biggl( \int_T k(t,t) d\nu(t)\biggr)^{ 1 /2} < \infty\, .
\end{equation}
For these kernels, the embedding $I_k:\to H$ is actually Hilbert-Schmidt, 
see e.g.~\cite[Lemma 2.3]{StSc12a}. Finally note that $\snorm k_{\sLtn}<\infty$ 
is always satisfied for bounded kernels as long as $\nu$ is a finite measure.

Now assume that $H$ is embedded into $\Lt \nu$.
Then one can show, see e.g.~\cite[Lemma 2.2]{StSc12a}, that 
 the adjoint $S_k:= I_k^*:\Ltn\to H$ of the embedding $I_k$ satisfies 
\begin{equation}\label{Iks}
S_k f(t)  = \int_T  k(t,t') f(t') d\nu(t')\, , \qquad \qquad f\in \Lt\nu, t\in T\, .
\end{equation}
We write $T_k := I_k \circ S_k$ for the resulting integral operator $T_k:\Ltn\to\Ltn$.
Clearly, $T_k$ is self-adjoint and positive, and if $H$ is compactly embedded, then 
$T_k$ is also compact, so that the classical spectral theorem for compact, self-adjoint
operators can be applied. In our situation, however, the spectral theorem can be refined, 
as we will see in Theorem \ref{spectral} below.
In order to formulate this theorem, we say that 
an at most countable family $(\a_i)_{i\in I}\subset (0,\infty)$ converges to $0$
if either $I=\{1,\dots,n\}$ or $I=\N:=\{1,2,\dots\}$ and $\lim_{i\to \infty} \a_i = 0$.
Analogously, when we consider  an at most countable family $(e_i)_{i\in I}$, 
we always assume without loss of generality that either $I=\{1,\dots,n\}$ or $I=\N$.

With these preparation we can now state the following spectral theorem for $T_k$, which
is an abbreviated version of \cite[Lemma 2.12]{StSc12a}:

\begin{theorem}\label{spectral}
Let $(T,\ca B, \nu)$ be a measure space and $k$ be a measurable kernel on $T$ whose  RKHS $H$ 
 is compactly embedded into $\Lt \nu$. 
Then there exists an at most countable  family 
$(\mu_i)_{i\in I}\subset (0,\infty)$
converging to 0 with $\mu_1\geq \mu_2\geq \dots> 0$ and a family $(e_i)_{i\in I}\subset H$ 
such that:
\begin{enumerate}
 \item The family $(\sqrt{\mu_i}e_i)_{i\in I}$ is an ONS in $H$ and $([e_i]_\sim)_{i\in I}$ is an ONS in $\Ltn$.
 \item The operator $T_k$ enjoys the following spectral representation, which is convergent in $\Ltn$:
	\begin{equation}\label{Tk-spectral}
	T_kf = \sum_{i \in I} \mu_i \bigl\langle f,[e_i]_\sim\bigr\rangle_{\Lx 2 \nu} [e_i]_\sim \, ,
	\qquad \qquad f\in \Ltn\, .
	\end{equation}
\item The following identities hold:
\begin{eqnarray}
	\label{Sk-EW}
	\mu_i e_i  & = & S_k [e_i]_{\sim}\, , \qquad \qquad i\in I \\ \label{Sk-ker}
	\ker S_k &= & \ker T_k\\ \label{Sk-ran}
	\overline{\ran S_k} &= &\overline{\spann\{\sqrt{\mu_i}e_i:i\in I  \}}\\ \label{Sks-ran-or}
	\overline{\ran S_k^*} &= &\overline{\spann\{[e_i]_\sim:i\in I  \}}\\ \label{Sks-ker}
	\ker S_k^*  &= & (\overline{\ran S_k})^\perp \\ \label{Sks-ran}
	\overline{\ran S_k^*} & = & (\ker S_k)^\perp\, ,
\end{eqnarray}
where the closures and orthogonal complements are taken in the spaces the objects are naturally contained in, that is,
\eqref{Sk-ran} and \eqref{Sks-ker} are considered in $H$, while 
\eqref{Sks-ran-or} and \eqref{Sks-ran} are considered in $\Ltn$. 
\end{enumerate}
\end{theorem}

The following set of assumptions, which is frequently used throughout this paper, essentially summarizes 
some notations from Theorem \ref{spectral}.

\begin{assumptionsk}
    Let $(T,\ca B, \nu)$ be a measure space and $k$ be a measurable kernel on $T$ whose  RKHS $H$ 
 is compactly embedded into $\Lt \nu$. Furthermore, let 
$(\mu_i)_{i\in I}$ and $(e_i)_{i\in I}$ be as in Theorem \ref{spectral}.
\end{assumptionsk}

For later use we note that if Assumption K is satisfied, then we automatically have 
\begin{equation}\label{diag-sum}
	\sum_{i\in I} \mu_i e_i^2(t)  \leq  k(t,t)\, , \qquad \qquad \qquad t\in T
\end{equation}
by \cite[Theorem 4.20]{StCh08}.
With the help of these families 
$(\mu_i)_{i\in I}$ and $(e_i)_{i\in I}\subset H$, some spaces and new kernels were defined in \cite{StSc12a},
which we need to recall since they are essential for this work.
To begin with,  
\cite[Equation (36)]{StSc12a} introduced, for $\b\in (0, 1]$, 
the subspace 
\begin{align}\label{def-power}
[H]^\b_\sim := \biggl\{ \sum_{i\in I}  a_i \mu_i^{\b/2} [e_i]_\sim: (a_i)\in \ell_2(I)    \biggr\} 
\end{align}
of $\Ltn$
and equipped it with the Hilbert space norm 
\begin{equation}\label{Hnk-norm}
   \bnorm {\sum_{i\in I} a_i \mu_i^{\b/2} [e_i]_\sim}_{[H]_\sim^\b } := \snorm {(a_i)}_{\ell_2(I) }\, .
\end{equation}
It is easy to verify that 
$(\mu_i^{\b/2} [e_i]_\sim)_{i\in I}$ is an ONB of $[H]_\sim^\b$ and that the set $\eHnb$ is independent 
of the particular choice of the family
$(e_i)_{i\in I}\subset H$ in Theorem \ref{spectral}.
In particular, \cite[Theorem 4.6]{StSc12a}
showed that 
\begin{equation}\label{interpol-repres}
 \ran T_k^{\b/2} = \eHnb = \bigl[\Lx 2 \nu , [H]_\sim\bigr]_{\b,2}   \,,
\end{equation}
where $T_k^{\b/2}$ denotes the $\b/2$-power of the operator $T_k$ defined, as usual, by its spectral representation,
and $[\Lx 2 \nu , [H]_\sim]_{\b,2}$ stands for the interpolation space of the standard real interpolation
method, see e.g.~\cite[Definition 1.7 on page 299]{BeSh88}. In addition, \cite[Theorem 4.6]{StSc12a}
showed that the norms of $\eHnb$ and $[\Lx 2 \nu , [H]_\sim]_{\b,2}$ are equivalent.
In other words, modulo equivalence of norms, $\eHnb$ \emph{is} the interpolation space $[\Lx 2 \nu , [H]_\sim]_{\b,2}$.

The spaces 
$\ran T_k^{\b/2}$ and $\eHnb$ almost naturally appear for the description of sample path properties,
see Section \ref{sec:path-in-interpol}, but unfortunately,
we rarely have a description of them that, unlike \eqref{def-power}, does not involve the eigenfunctions and -values.
Since the latter are notoriously difficult to find,  sample path descriptions solely in terms of 
$\eHnb$ would therefore be of little practical value.
On the other hand,  interpolation spaces of the real method are a long studied object and in many cases we have
alternative descriptions of these spaces. We refer to the case of $H$ being a Sobolve space and $\nu$ being the Lebesgue
measure for the probably most classical example in which such alternative descriptions exists.
Consequently, the key feature of \eqref{interpol-repres} in view of this discussion is that 
it will link sample path properties to the rich theory of interpolation spaces via 
the naturally appearing spaces 
$\ran T_k^{\b/2}$ and $\eHnb$.

In \cite[Section 4]{StSc12a} it was shown that
under certain circumstances
$[H]_\sim^\b$ is actually the image of an RKHS under $[\mycdot]_\sim$. To recall the construction
of this RKHS, let us assume that we have a measurable $S\subset T$ with $\nu(T\setminus S)=0$ and 
\begin{equation}\label{k-series-diag}
 \sum_{i\in I} \mu_i^\b e_i^2(t) < \infty\, , \qquad \qquad t\in S 
\end{equation}
We write 
 $\hat e_i := \eins_{S} e_i$ for all $i\in I$. 
Clearly, this gives 
$
 \sum_{i\in I} \mu_i^\b \hat e_i^2(t) < \infty 
$ for all $t\in T$.
Based on this and the fact that $([\hat e_i]_\sim)_{i\in I}$ is an ONS of $\Ltn$,
\cite[Lemma 2.6]{StSc12a} showed that 
\begin{equation}\label{rkhs-induc-form}
 \bHSb := \biggl\{ \sum_{i\in I}  a_i \mu_i^{\b/2}\, \hat e_i: (a_i)\in \ell_2(I)    \biggr\} 
\end{equation}
equipped with the norm 
\begin{equation}\label{rkhs-induc-norm}
\bnorm {\sum_{i\in I} a_i \mu_i^{\b/2}\,  \hat e_i}_{\bHSb} := \snorm {(a_i)}_{\ell_2(I) }
\end{equation}
is a separable RKHS,
which is compactly embedded into $\Ltn$.
Moreover, the family $(\mu_i^{\b/2}\hat e_i)_{i\in I}$ is an ONB of $\bHSb$ and the (measurable) kernel $\bkSb$ of $\bHSb$
is given by the pointwise convergent series representation
\begin{equation}\label{kern-sum-basic}
\bkSb(t,t') = \sum_{i\in I} \mu_i^\b \hat e_i(t)\hat e_i(t')\, , \qquad \qquad t,t'\in T\, .
\end{equation}
Recall that
 $\bkSb$ and its RKHS $\bHSb$ are actually independent  
of the particular choice of the family
$(e_i)_{i\in I}\subset H$ in Theorem \ref{spectral}, 
see \cite[Proposition 4.2]{StSc12a}.
This justifies the chosen notation.
Furthermore, note that in general, $\hat k_T^1$ does \emph{not} equal $k$, and, of course, the same is true
for the resulting RKHSs $\bHxy T1$ and $H$. In fact, \cite[Theorem 3.3]{StSc12a} 
shows that  $k=\hat k_T^1$ holds, if and only if $I_k:H\to \Ltn$ is injective,
and a sufficient condition for the latter will be presented in Lemma \ref{Ik-inject}.
Finally, 
for $\a\geq \b$, we have $\bHxy S\a \hookrightarrow \bHSb$ by the definition of 
the involved norms, cf.~also the proof of
\cite[Lemma 4.3]{StSc12a}.

In the following, we write $\kSb:S\times S\to \R$ for the restriction of $\bkSb$ onto $S\times S$
and we denote the RKHS of $\kSb$ by $\HSb$. 
Formally, the spaces $\bHSb$, $\HSb$ and $[H]^\b_\sim$ are different. 
Not surprisingly, however, they are all isometrically isomorphic to each other via natural operators.
The corresponding results are collected  in Lemma \ref{basic-isometries}.

%


%
%

Let us now recall conditions, which ensure \eqref{k-series-diag} for a set $S$ of full measure.
To begin with, note that we find such an $S$ if $\sum_{i\in I}\mu_i^\b<\infty$, since a simple 
calculation based on Beppo Levi's theorem shows
\begin{equation}\label{when-kS}
   \int_T \sum_{i\in I} \mu_i^\b e_i^2(t) \dxy\nu t
   =  \sum_{i\in I} \mu_i^\b \int_Te_i^2(t) \dxy\nu t
   = \sum_{i\in I}\mu_i^\b<\infty\, .
\end{equation}
Moreover, in this case we obviously have $\snorm\bkSb_{\sLtn}<\infty$.
Interestingly, the converse implication is also true, namely 
\cite[Proposition 4.4]{StSc12a} showed
that we have  $\sum_{i\in I}\mu_i^\b<\infty$, if and only if \eqref{k-series-diag} holds for  a set $S$ of full measure
and the resulting kernel $\bkSb$ satisfies $\snorm\bkSb_{\sLtn}<\infty$. 
Moreover, using the theory of liftings \cite[Theorem 5.3]{StSc12a} showed that
\eqref{k-series-diag} holds for a set $S$ of full measure, if 
$\nu$ is a $\s$-finite measure for which $\ca B$ is complete
and
\begin{equation}\label{linfty-inclus}
   [H]_\sim^\b \hookrightarrow \Lx \infty \nu\, .
\end{equation}
Note that this
 sufficient condition is particularly interesting when combined with \eqref{interpol-repres},
since the inclusion $[\Ltn , [H]_\sim]_{\b,2}\hookrightarrow \Lx \infty \nu$
may be known in specific situations.
Finally,  \cite[Theorem 5.3]{StSc12a} actually showed that
the inclusion \eqref{linfty-inclus} holds, if and only if 
\eqref{k-series-diag} holds for  a set $S$ of full measure \emph{and} 
the resulting kernel $\bkSb$ is bounded. 

Our next goal is to investigate under which conditions \eqref{k-series-diag} actually holds for $S:=T$.
%
%
To this end, let us now assume that we have a topology $\t$ on $T$.
The following definition introduces some notions of continuity for $k$.

\begin{definition}\label{k-conts-def}
 Let $(T,\tau)$ be a topological space and $k$ be a kernel on $T$ with RKHS $H$. Then we say that
$k$ is:
\begin{enumerate}
\item $\t$-continuous, if $k$ is continuous with respect to the product topology $\t\otimes \t$.
 \item  separately $\t$-continuous, if
$k(t,\cdot):T\to \R$ is $\t$-continuous for all $t\in T$.
\item weakly $\t$-continuous, if all $f\in H$ are $\t$-continuous.
\end{enumerate}
\end{definition}

Clearly, $\t$-continuous kernels are  separately $\t$-continuous.
Moreover, it is a well-known fact that given a  $\t$-continuous kernel $k$
its canonical feature map $\P:T\to H$ defined by $\P(t) := k(t,\cdot)$
is $\t$-continuous, see e.g.~\cite[Lemma 4.29]{StCh08}, and hence   the reproducing
property $f= \langle f, \P(\cdot)\rangle_H$, which holds for all $f\in H$, shows that $k$ is 
also weakly $\t$-continuous. Moreover, \cite[Lemma 4.28]{StCh08} shows that bounded, 
separately $\t$-continuous kernels are weakly $\t$-continuous, too. In this regard note that
even on $T= [0,1]$
not every bounded, separately $\t$-continuous kernel is continuous, see \cite{Lehto50a}.

Let us now introduce two topologies on $T$ generated by $k$ and its RKHS $H$.
The first one is the topology $\t_k$ generated by the well-known 
 pseudo-metric $d_k$ on $T$ defined by
\begin{displaymath}
 d_k(t,t') :=\snorm{\P(t) - \P(t')}_H\, , \qquad \qquad t,t'\in T.
\end{displaymath}
Obviously, this pseudo-metric is a metric if and only if the canonical feature map $\P:T\to H$
is injective, and  this is also 
the only case in which $\t_k$ is Hausdorff.
Less often used is another topology on $T$
that is related to $k$,
namely the initial topology $\t(H)$ generated by the set of functions $H$. In other words, $\t(H)$ is the smallest
topology on $T$ for which all $f\in H$ are continuous, that is, for which $k$ is weakly $\t$-continuous.
More information on these topologies can be found in Lemma \ref{simple-topol}.

In the following, we sometimes need measures $\nu$ 
that are strictly positive on all non-empty $\t(H)$-open sets. Such measures are introduced in the following definition.

\begin{definition}
   Let $(T,\ca B, \nu)$ be a measure space and $k$ be a kernel on $T$ with RKHS $H$
   such that  $\t(H)\subset \ca B$.
   Then $\nu$ is called  $k$-positive, if, for all non-empty $O\in \t(H)$, we have $\nu(O)>0$.
\end{definition}

The notion of $k$-positive measures generalizes that of strictly positive measures.
Indeed, if $(T,\t)$ is a topological space,  and $\ca B:= \s(\t)$ is the 
corresponding Borel $\s$-algebra, then a measure $\nu$ on $\ca B$ 
is strictly positive, if $\nu(O)>0$ for all non-empty $O\in \t$.
Now assume that we have a (weakly)-$\t$-continuous 
kernel $k$ on $T$. Then we find $\t(H)\subset \t \subset \ca B$, and thus 
$\nu$ is also $k$-positive.

Note that if $H$ is separable and $k$ is both bounded and $\ca B\otimes \ca B$-measurable, 
then every $f\in H$ is $\ca B$-measurable, see e.g.~\cite[Lemma 4.25]{StCh08},
and hence $\s(H)\subset \ca B$.
By part \emph{iii)} of Lemma \ref{simple-topol} we will thus find $\t(H)\subset \s(H)\subset \ca B$.
In other words, the assumption  $\t(H)\subset \ca B$, which will occur frequently, is automatically satisfied
for such $H$.

%
%

The following simple lemma gives a first glance at the importance of $k$-positive measures.

\begin{lemma}\label{Ik-inject}
 Let $(T,\ca B, \nu)$ be a measure space and $k$ be a kernel on $T$ with RKHS $H$ such that  $\t(H)\subset \ca B$.
If  $\nu$ is $k$-positive, then
	$I_k:H\to \Ltn$ is injective and $k=k_T^1$.
\end{lemma}

Let us now collect a set of assumptions frequently used when dealing with 
$k$-positive measures.

\begin{assumptionsc}
  Let $(T, \ca B, \nu)$ be a $\s$-finite and complete measure space and $k$ be a kernel on $T$ 
with RKHS $H$ such that 
 $\t(H)\subset \ca B$ and $\nu$ is $k$-positive. Furthermore, Assumption K is satisfied. 
\end{assumptionsc}

With these preparations we are now in the position to improve the result on bounded $\kSb$
from \cite[Theorem 5.3]{StSc12a}.

\begin{theorem}\label{bounded-powers-of-cont-kernels}
Let Assumption CK be satisfied.
%
 Furthermore, assume that 
  for some  $0<\b\leq  1$, we have 
\begin{equation}\label{Linfty-embed}
   \bigl[\Ltn , [H]_\sim\bigr]_{\b,2}\hookrightarrow \Lx \infty {\nu}\, .
\end{equation}
  Then,   \eqref{k-series-diag} holds for $S:=T$, the resulting kernel $\kTb$  is bounded, and 
$\t(\HTb)= \t(H)$.
\end{theorem}

Note that under the assumptions of Theorem \ref{bounded-powers-of-cont-kernels} we also 
have $\sum_{i\in I}\mu_i^\b<\infty$ provided that $\nu$ is \emph{finite}, see \cite[Theorem 5.3]{StSc12a}.
In addition, there is a partial converse,
which does not need any continuity assumption.
Indeed, if
  we have $\sup_{i\in I} \inorm{e_i}< \infty$, then a simple estimate shows that 
 $\sum_{i\in I}\mu_i^\b<\infty$  implies  
 \eqref{k-series-diag}  for $S:=T$, and  the resulting kernel $\kTb$ turns out to be bounded.

To illustrate the theorem above, let us assume that 
$(T,\t)$ is a topological space. In addition, let $\ca B$ be a $\s$-algebra on $T$ and $\nu$ be a $\s$-finite
and strictly positive 
measure on $\ca B$ such that $\ca B$ is $\nu$-complete and $\t\subset \ca B$.
If $k$ is a weakly $\t$-continuous kernel on $T$, we then obtain $\t(H)\subset \t\subset \ca B$,
where $H$ is the RKHS of $k$.
Consequently, if  $H$ is 
compactly embedded into $\Ltn$, and, for some  $0<\b\leq  1$, we have \eqref{Linfty-embed},
then the assumptions of Theorem \ref{bounded-powers-of-cont-kernels} are satisfied, and 
hence 
$\kTb$ is defined and bounded.
Moreover, we have $\t(\HTb)= \t(H)\subset \t$, that is, $\kTb$ is weakly $\t$-continuous. 
In other words, modulo the technical assumptions of Theorem \ref{bounded-powers-of-cont-kernels},
the embedding \eqref{Linfty-embed} ensures that $\kTb$ is defined and inherits the weak continuity 
from $k$.

Since in Section \ref{sec:paths-in-rkhs}
we will  investigate inclusions between powers of  RKHSs in more detail,
let us introduce  some more notations. To this end, we fix two kernels $k_1, k_2$ on $T$ with 
corresponding RKHSs $H_1$ and $H_2$.
Following \cite{LuBe01a} we say that $k_2$ \emph{dominates} $k_1$ and 
 write $k_1\leq k_2$, if 
$H_1\subset H_2$ and the natural inclusion
operator $I_{k_1,k_2}:H_1\to H_2$ is continuous. In this case, the adjoint 
operator $I_{k_1,k_2}^*:H_2\to H_1$ exists and is continuous. In analogy to our previous
notations, we write $S_{k_1,k_2} := I_{k_1,k_2}^*$.
Moreover, we speak of nuclear dominance and write 
$k_1\ll k_2$, if $k_1\leq k_2$ and $I_{k_1,k_2}\circ S_{k_1,k_2}$ is nuclear.
Lemma \ref{nuclear-dom} characterizes when $\kxy S1 \ll \kxy S\b$ holds.

Many of our results are formulated in terms of the eigenvalues $(\mu_i)_{i\in I}$,
but
determining these eigenvalues in 
a specific situation 
 is often an extremely difficult task. For many of our results, we need,
however, only the \emph{asymptotic behavior} of the eigenvalues. It is well-known,
see e.g.~\cite{CaSt90,EdTr96}, that this behavior can often be determined by 
entropy numbers.
%
Our next goal is to 
make this statement precise. 
To this end, recall that 
%
the $i$-th (dyadic) entropy number of a compact, linear operator $T:E\to F$ 
between Banach spaces $E$ and $F$
is defined by 
\begin{displaymath}
   \e_{i}(T)
:=
\inf\biggl\{ \e>0: \exists\, y_{1},\dots,y_{2^{i-1}}\in F \mbox{ such that } TB_E\subset \bigcup_{j=1}^{2^{i-1}} (y_j + \e B_F) 
\biggr\} \, .
\end{displaymath}
Note that in the literature these numbers are usually denoted by $e_i(T)$, instead. Since
this in conflict with our notation for eigenvectors, we departed from this convention.
For an introduction to these numbers we refer to the above mentioned books
\cite{CaSt90,EdTr96}. 

Now the following, somewhat folklore, result compares the eigenvalues 
$(\mu_i)_{i\in I}$ with the entropy numbers of $I_k$. Note that the latter are often asymptotically known, see 
\eqref{sob-entrop} below for an example.

\begin{lemma}\label{entropy-versus-ew}
   Let Assumption K be satisfied. Then, for all $i\in I$,  we have 
   \begin{equation}\label{inverse-car-ineq}
      \mu_i \leq 4 \e_i^2(I_k)\, .
   \end{equation}
  Moreover, for all $\b>0$, there exists a constant $c_\b>0$ such that 
   \begin{equation}\label{carl-ineq}
      \sum_{i=1}^\infty \e_i^{2\b}(I_k) \leq c_\b \sum_{i\in I} \mu_i^\b
   \end{equation}
   In particular, for all $\b>0$ we have $\sum_{i\in I} \mu_i^\b< \infty$ if and only if 
   $\sum_{i=1}^\infty \e_i^{2\b}(I_k)<\infty$.
   Similarly, for all $\b>0$ we have 
   \begin{equation}
    \mu_i\preceq i^{-\b} \qquad \qquad \Longleftrightarrow \qquad \qquad \e_i(I_k) \preceq  i^{-\b/2} 
   \end{equation}
    as well as 
    $\mu_i\asymp i^{-\b}$
   if and only if $\e_i(I_k) \asymp  i^{-\b/2}$.
\end{lemma}

Finally, to describe    some higher order smoothness properties of functions, we fix 
a non-empty open and bounded $T\subset \Rd$ that satisfies the strong local Lipschitz condition
of \cite[p.~83]{AdFo03}. 
Note that the strong local Lipschitz condition is
satisfied for e.g.~the interior of $[0,1]^d$ or open Euclidean balls. 
We write $\Lx 2 T$ for the $L_2$-space with respect to the Lebesgue measure on $T$. 
For $m\in \N_0$ and $p\in [1,\infty]$ we denote the
Sobolev space of smoothness $m$ by  $W^{m,p}(T)$, that is 
\begin{displaymath}
 W^{m,p}(T)
 := \bigl\{ f\in \Lx p T:  D^{(\a)} \mbox{ exists and } D^{(\a)}f \in \Lx p T 
\mbox{ for all } \a\in \N_0^d  \mbox{ with } |\a|\leq m\bigr\}\, ,
\end{displaymath}
where, as usual, $D^{(\a)}f$ denotes the weak $\a$-partial derivative of $f$.
%
%
%
For notational simplicity, we further write 
$W^m(T):= W^{m,2}(T)$.
 Recall Sobolev's embedding 
theorem, see e.g.~\cite[Theorem 4.12]{AdFo03},
which ensures $W^m(T)\hookrightarrow C(\overline T)$ for all $m>d/2$,
where $C(\overline T)$ denotes the space of continuous functions defined on the closure $\overline T$ of $T$.
For such $m$ we can thus view $W^m(T)$ as an RKHS on $T$.
Following tradition, we will, however, not notationally distinguish between
the cases in which $W^m(T)$ is viewed as a space of equivalence classes or
as a space of functions, since the meaning of the symbol $W^m(T)$  will always be clear from the context.

We further need fractional versions of Sobolev spaces and generalizations of them. To this end
recall from \cite[p.~230]{AdFo03} that 
 the Besov spaces of smoothness $s>0$ are given by 
 \begin{equation}\label{besov-def}
   B^s_{p,q}(T) := \bigl[\Lx p T, W^{m,p}(T)  \bigr]_{s/m,q}\, ,
 \end{equation}
where $m>s$ is an arbitrary natural number and $p,q\in [1,\infty]$. 
Recall that for  $s>d/2$,
we again have a continuous embedding 
 $B^s_{2,2}(T)\hookrightarrow C(\overline T)$, see 
\cite[Theorem 7.37]{AdFo03}. Moreover, we have $B^m_{2,2}(T) = W^m(T)$ for all integers $m\geq 1$, see  
\cite[p.~230]{AdFo03}, and for this reason we often use the notation $W^s(T) :=W^{s,2}(T) := B^s_{2,2}(T)$ for all $s>0$.
Note that with this notation, the equality in \eqref{besov-def} with $p=q=2$ actually holds for all real $m>s$
by the reiteration property of the real interpolation method, 
see
again \cite[p.~230]{AdFo03}.
Finally, for $0<s<1$ and $p\in [1,\infty]$,  we have by \cite[Lemma 36.1 and p.~170]{Tartar07}
\begin{displaymath}
   B^s_{p,p}(T) = \Bigl\{ f\in \Lx p T \, : \,  \snorm f_{T,s,p}<\infty   \Bigr\}\, ,
\end{displaymath}
where 
\begin{displaymath}
   \snorm f_{T,s,p}^p := \int_T\int_T \frac {|f(r)-f(t)|^p}{|r-t|^{d+sp}} \, dr\, dt
\end{displaymath}
with the usual modification for $p=\infty$.
Similarly, if  $s>1$ is not an integer, $B^s_{p,p}(T)$ equals the fractional Sobolev-Slobodeckij spaces, 
i.e.~we have 
\begin{displaymath}
   B^s_{p,p}(T)= \bigl\{f\in W^{\lfloor s\rfloor}_{p,p}(T) : \snorm{D^{(\a)} f}_{T,s-\lfloor s\rfloor,p}<\infty 
\mbox{ for all } \a\in \N_0^d \mbox{ with } |\a| = \lfloor s \rfloor\bigr\}\, .
\end{displaymath}
We refer to \cite[p.~156]{Tartar07} and \cite{DeSh93a} for details.
In particular, 
$B^s_{\infty,\infty}(T)$ is the space of $s$-H\"older continuous functions for all $0<s<1$. 

Let us finally recall some entropy estimates related to fractional Sobolev spaces. 
To this end,  let $T\subset \Rd$ be a bounded subset 
that satisfies 
the strong local Lipschitz condition and $T=\interior \overline T$, where $\interior A$ denotes the interior of $A$.
Then \cite[p.~151]{EdTr96} shows that, for all $s>d/2$, we have 
\begin{equation}\label{sob-entrop}
    \e_i\bigl(\id:  W^s(T)  \to \Lx 2 T\bigr) \asymp i^{-s/d}\, ,
\end{equation}
where we used the notation $W^s(T) = B^s_{2,2}(T)$ introduced above.


\section{Karhunen-Lo{\`e}ve Expansions For Generic Processes}\label{sec:KL-generic}

The goal of this section is to establish a Karhunen-Lo{\`e}ve expansion
that does not require the usual assumptions such as  compact index sets $T$ or continuous kernels $k$.
To this end, we first show that under very generic assumptions the covariance
function of a centered, second-order process satisfies Assumption K, so that the theory developed in Section \ref{sec:prelim}
is applicable. We then 
repeat the classical Karhunen-Lo{\`e}ve 
approach and  combine it with some further aspects from
Section \ref{sec:prelim}.

In the following, let $(\Om, \ca A, P)$ be a probability space and 
 $(T, \ca B, \nu)$ be a $\s$-finite  measure space. Given a stochastic process 
 $(X_t)_{t\in T}$  on $\Om$, we denote the path $t\mapsto X_t(\om)$ of a given 
 $\om \in \Om$ 
by  $X(\om)$.
 Moreover, 
 we call the process $(\ca A\otimes \ca B)$-measurable, if 
 the map $X:\Om\times T\to \R$ defined by 
$(\om, t)\mapsto X_t(\om)$ is measurable. In this case,  each 
path is obviously $\ca B$-measurable.

Let us  assume that $X$ is centered, second-order, that is $X_t\in \sLt P$ and $\E_P X_t = 0$  for all $t\in T$.
Then
the covariance function
$k:T\times T\to \R$ is   given by 
\begin{displaymath}
 k(s,t) := \E_P X_s X_t\, , \qquad \qquad s,t\in T\, .
\end{displaymath}
It is well-known, see e.g.~\cite[p.~57]{BeTA04}, that the covariance function
is symmetric and positive semi-definite, and thus a kernel by the 
Moore-Aronszajn theorem, see e.g.~\cite[Theorem 4.16]{StCh08}.

Let us now additionally  assume that  
 $\nu$ is suitably chosen in the sense of $X\in \sLt {P\otimes \nu}$.
For $P$-almost all $\om \in \Om$, we then have 
$X(\om) \in \sLt \nu$. For such $X$, the
  following lemma collects some additional properties of the covariance function.

\begin{lemma}\label{cov-prop}
Let $(\Om, \ca A, P)$ be a probability space and 
 $(T, \ca B, \nu)$ be a $\s$-finite  measure space. In addition, let $(X_t)_{t\in T} \subset \sLt P$ be a 
 centered and $(\ca A\otimes \ca B)$-measurable stochastic process. 
 Then its covariance function 
$k:T\times T\to \R$ is measurable, and we have   $X\in \sLt {P\otimes \nu}$ if and only if
\begin{displaymath}
   \int_T k(t,t)\, d\nu(t) < \infty\, .
\end{displaymath}
Moreover, the RKHS $H$ of $k$ is compactly embedded into $\Ltn$ and the corresponding integral operator 
$T_k:\Ltn\to \Ltn$ is nuclear.
\end{lemma}

The lemma above in particular shows that for a stochastic process $X\in \sLt {P\otimes \nu}$
and its covariance function $k$ 
Theorem \ref{spectral} applies. Let us therefore assume that
 $(e_i)_{i\in I}\subset H$
and  $(\mu_i)_{i\in I}$ are families  satisfying the assertions of 
Theorem \ref{spectral}.
For $i\in I$ we then define $Z_i:\Om \to \R$ by
\begin{equation}\label{def-Zi}
 Z_i (\om) := \int_T X_t(\om ) e_i(t) \, d\nu(t) 
\end{equation}
for all $\om\in \Om\setminus N$, where $N\subset \Om$ is a measurable subset satisfying
with $P(N) = 0$ and $X(\om)\in \sLtn$ for all $\om \in \Om\setminus N$. For $\om \in N$ we further write  
$Z_i(\om) := 0$. 
Clearly, each $Z_i$ is measurable and 
$Z_i(\om) = \langle [X(\om)]_\sim, [e_i]_\sim\rangle_{\Ltn}$ for $P$-almost all $\om\in \Om$.

Having finished  these preparations we can now formulate our assumptions on the process $X$
that will be used throughout the rest of this work.

\begin{assumptionsx}
Let $(\Om, \ca A, P)$ be a probability space and 
 $(T, \ca B, \nu)$ be a $\s$-finite  measure space such that $\Ltn$ is separable. 
In addition, let $(X_t)_{t\in T} \subset \sLt P$ be a 
 centered and $(\ca A\otimes \ca B)$-measurable stochastic process  such that  $X\in \sLt {P\otimes \nu}$. 
Moreover, let $k$ be its covariance function and $H$ be the RKHS of $k$.
Finally, let $(e_i)_{i\in I}\subset H$
and  $(\mu_i)_{i\in I}$ be as in Theorem \ref{spectral} and $(Z_i)_{i\in I}$ be defined by \eqref{def-Zi}.
\end{assumptionsx}

At first glance, the assumed $(\ca A\otimes \ca B)$-measurability may look restrictive, but 
it is satisfied if, e.g.~$T$ is a Polish space and $X$ has continuous paths.
Moreover, we will later see in Theorem \ref{def-X-by-Z-new} that the $(\ca A\otimes \ca B)$-measurability
is actually necessary for the expansions developed below.

Assumption X explicitly involves the eigenvalues and eigenfunctions of the integral operator $T_k$.
The next simple lemma provides a set of conditions implying Assumption X, which in some cases are easier to check.
Another such alternative set of conditions will be presented later in Theorem \ref{unique-kle}.

\begin{lemma}\label{kernel-instead-of-EV}
Let $(\Om, \ca A, P)$ be a probability space,  $(T,\ca B)$ be a measurable space, and 
 $(X_t)_{t\in T} \subset \sLt P$ be a 
 centered and $(\ca A\otimes \ca B)$-measurable stochastic process.
 Assume that its covariance function $k:T\times T\to \R$ has, for an at most countable index set $I$,  a
 pointwise convergent representation
 \begin{equation}\label{kernel-instead-of-EV-h1}
    k(t,t') = \sum_{i\in I} \mu_i e_i(t) e_i(t')\, ,  \qquad \qquad t,t'\in T,
 \end{equation}
  where the functions  $e_i:T\to \R$ are measurable, and the family 
  $(\mu_i)_{i\in I}$ satisfies both $(\mu_i)_{i\in I}\subset (0,\infty)$ and 
  $\sum_{i\in I} \mu_i < \infty$. If $\nu$ is a $\s$-finite measure on $(T,\ca B)$
  such that $\Ltn$ is separable and
	$([e_i]_\sim)_{i\in I}\subset \Ltn$ is an ONS, then Assumption X is satisfied and
  $(e_i)_{i\in I}\subset H$
	and  $(\mu_i)_{i\in I}$ are the families considered in Assumption X. In particular, these two families 
	satisfy Assumption K
\end{lemma}

The following lemma, which is somewhat folklore,  shows that for processes satisfying Assumption X
an  expansion of the form  \eqref{Loeve-iso-representation} can be obtained if we replace 
 $\xi_i:= \Psi^{-1}(\sqrt{\mu_i} e_i)$ by $\mu_i^{-1/2}Z_i$.
The proof of this lemma  does not deviate much from the one needed for
the classical Karhunen-Lo{\`e}ve expansion, but
since the traditional assumptions for this expansion are more 
restricted and the lemma itself is the very foundation of our following results
we have 
 included it for the sake of completeness.

\begin{lemma}\label{prop-Zi}
Let Assumption X be satisfied. Then, for all $i,j\in I$ and $t\in T$,
 we have $Z_i \in \sLt P$ with 
$\E_P Z_i = 0$ and 
\begin{align}\label{Zi-cov}
 \E_P Z_i Z_j &= \mu_i \d_{i,j} \, ,\\ \label{Zi-Xt-prod}
 \E_P Z_i X_t &= \mu_i e_i(t)\, .
\end{align}
Moreover, for  all finite $J\subset I$ and all $t\in T$ we have 
\begin{equation}\label{LP-diff}
 \bnorm{X_t - \sum_{j\in J} Z_j e_j(t)}_{\sLt P}^2 = k(t,t) - \sum_{j\in  J} \mu_j e_j^2(t)\, ,
\end{equation}
and, for a fixed $t\in T$, the following statements are equivalent:
\begin{enumerate}
   \item With convergence in $\Lx 2P$ we have 
    \begin{equation}\label{KL-L2P-conv}
 [X_t]_\sim = \sum_{i\in I} [Z_i]_\sim e_i(t)\, .
  \end{equation}
  \item We have \begin{equation}\label{kernel-repres-diag}
 k(t,t) = \sum_{i\in  I} \mu_i e_i^2(t)\, .
  \end{equation}
\end{enumerate}
Moreover, if, for some $t\in T$, we have \eqref{KL-L2P-conv}, then 
the convergence in \eqref{KL-L2P-conv} is necessarily unconditional in $\Lx 2 P$  by \eqref{LP-diff}.
Finally, there exists a measurable $N\subset \Om$ with $\nu(N)=0$ such that for all $\om \in \Om\setminus N$ we have 
\begin{equation}\label{path-in-compl-kernel}
 [X(\om)]_\sim \in (\ker T_k)^\perp = \overline{\spann\{[e_i]_\sim:i\in I  \}}^{\Ltn}\, .
\end{equation}
\end{lemma}


Recall that for continuous kernels $k$ over compact metric spaces $T$ and strictly positive measures $\nu$,
Equation \eqref{kernel-repres-diag} is guaranteed 
by the classical theorem of Mercer for \emph{all} $t\in T$. 
Moreover, since the convergence in \eqref{kernel-repres-diag} is 
monotone and $t\mapsto k(t,t)$ is continuous, 
Dini's theorem shows in this case, that the convergence in \eqref{kernel-repres-diag} is uniform in $t$.
By \eqref{LP-diff} we conclude that the $\Lx 2 P$-convergence in \eqref{KL-L2P-conv} is also \emph{uniform} in $t$.
 In the general case, however, \eqref{kernel-repres-diag}
may no longer be true. 
Indeed, the following proposition
 \emph{characterizes} when  
\eqref{kernel-repres-diag} holds. In addition, it shows that for \emph{separable} $H$
Equation \eqref{KL-L2P-conv} holds at least $\nu$-almost everywhere.

\begin{proposition}\label{classical-KL}
Let Assumption X be satisfied.
Then the following statements are equivalent:
\begin{enumerate}
 \item The family $(\sqrt{\mu_i}e_i)_{i\in I}$ is an ONB of $H$.
	\item The operator $I_k:H\to \Ltn$ is injective.
	\item For all $t\in T$ we have \eqref{KL-L2P-conv}.
\end{enumerate}
Moreover, if $H$ is separable, there exists a measurable $N\subset T$ with $\nu(N) = 0$ 
such that \eqref{KL-L2P-conv} holds with unconditional convergence in $\Lx 2P$ for all $t\in T\setminus N$.
\end{proposition}

Note that for $k$-positive measures $\nu$ the injectivity of 
$I_k:H\to \Ltn$ is automatically satisfied by Lemma \ref{Ik-inject}, and thus 
we have \eqref{KL-L2P-conv} for all $t\in T$.
%
%
Moreover note that the injectivity of $I_k$ must not be confound with the injectivity of $T_k$. Indeed, 
the latter is equivalent to $I_k:H\to \Ltn$ having a dense 
image, see \eqref{Sk-ker} and \eqref{Sks-ran}. Moreover, the injectivity of $T_k$ is also equivalent to 
$(|e_i]_\sim)_{i\in I}$ being an ONB of $\Ltn$, see \eqref{Sks-ran-or}.

Due to the particular version of convergence in  \eqref{KL-L2P-conv}, 
Proposition \ref{classical-KL} is useful for approximating the random variable $X_t$ at some given time $t$,
but useless for approximating the paths of the process $X$. This is addressed by the following result, which 
is the generic version of \eqref{KL-expansion-intro} and as such the first new result of this section.


\begin{theorem}\label{pointwise-KL-new}
 Let Assumption X be satisfied.
 Then 
there exists a measurable $N\subset \Om$ with $P(N)=0$ such that 
 for all $\om \in \Om\setminus N$ we have 
	\begin{equation}\label{KL-path-conv-Lt}
	 [X(\om)]_\sim  = \sum_{i\in I} Z_i(\om) [e_i]_\sim\, ,
	\end{equation}
where the convergence is unconditionally in $\Ltn$. Moreover, 
for all   $J\subset I$, we have 
\begin{equation}\label{Ln-diff}
 \int_\Om\bnorm{[X(\om)]_\sim - \sum_{j\in J} Z_j(\om) [e_j]_\sim}_{\Ltn}^2 \, dP(\om) = \sum_{i\in  I\setminus J} \mu_i\, .
\end{equation}
In particular, with unconditional convergence in $\Lx 2 {P\otimes \nu}$,  it holds 
\begin{displaymath}
 [X]_\sim = \sum_{i\in I} [Z_i]_\sim [e_i]_\sim\, .
\end{displaymath}
\end{theorem}

Equation \eqref{KL-path-conv-Lt} shows that almost every path can be approximated using the partial sums
$\sum_{j\in J} Z_j [e_j]_\sim$ while \eqref{Ln-diff} exactly specifies the average speed of 
convergence for such an approximation. In particular,
\eqref{Ln-diff} shows that any meaningful speed of convergence requires
stronger summability assumptions on the sequence $(\mu_i)_{i\in I}$ of eigenvalues.

The next corollary, which again generalizes earlier known results, relates the Lo{\`e}ve isometry to the random variables $Z_i$.

\begin{corollary}\label{Zi-in-CM}
    Let Assumption X be satisfied and $\Psi:\Lx 2 X\to H$ be the Lo{\`e}ve isometry, where
	    \begin{displaymath}
       \Lx 2 X := \overline{\spann\bigl\{[X_t]_\sim: t\in T\bigr\}}^{\Lx 2 P}
    \end{displaymath}
	denotes the Cameron-Martin space. 
Then, for all $i\in I$, we have 
	\begin{displaymath}
	 [Z_i]_\sim =  {\mu_i} \Psi^{-1}(e_i)\, ,
	\end{displaymath}
and  
the family $(\mu_i^{-1/2}[Z_i]_\sim)_{i\in I}$ is an ONS of $\Lx 2 X$.
	Moreover, it is an ONB, if and only if   $(\sqrt{\mu_i}e_i)_{i\in I}$ is an ONB of $H$.
\end{corollary}

Let us now birefly consider the case of Gaussian processes. To this end, let us recall that 
a process $(X_t)_{t\in T}$ is called Gaussian, if, for all $n\geq 1$, $a_1,\dots,a_n\in \R$,
and $t_1,\dots,t_n\in T$, the random variable $\sum_{i=1}^n a_iX_{t_i}$
has a normal distribution. 
The following lemma shows that for Gaussian processes, the $Z_i$'s are independent, normally
distributed random variables.

\begin{lemma}\label{gauss-lemma}
    Let $(X_t)_{t\in T}$ be a Gaussian process for which Assumption X is satisfied. 
    Then the random variables $([Z_i]_\sim)_{i\in I}$ are independent and for all $i\in I$, we have 
    $Z_i \sim \ca N(0,\mu_i)$. 
\end{lemma}

Our next result in this section  in particular shows that all reasonable
sequences of coefficient variables $(Z_i)_{i\in I}$ can occur in the class of processes satisfying Assumption X.
The main difficulty in its proof is the existence of the $(\ca A  \otimes \ca B)$-measurable version.
Note that the existence of this version shows that the $(\ca A  \otimes \ca B)$-measurability  
in Assumption X is necessary for processes having a Karhunen-Lo{\`e}ve expansion.

\begin{theorem}\label{def-X-by-Z-new}
Let $(T,\ca B, \nu)$ be a $\s$-finite   measure space and $k$ be a kernel on $T$ such that 
 Assumption K is satisfied with $\sum_{i\in I}\mu_i<\infty$. Moreover, let $(\Om, \ca A, P)$ be a probability space,
and 
$(Z_i)_{i\in I}\subset \sLx 2 P$ be a sequence of centered random variables such that 
\begin{equation}\label{orth-z-again}
   \E_P Z_i Z_j = \mu_i \d_{i,j}
\end{equation}
for all $i,j\in I$. For $t\in T$, we define 
\begin{equation}\label{def-X-by-Z-h1}
   X_t := \sum_{i\in I} Z_i e_i(t)\, ,
\end{equation}
where we note that the series  converges unconditionally in $\sLx 2 P$.
Then there exists an  $(\ca A  \otimes \ca B)$-measurable version $(Y)_{t\in T} \subset \sLx 2 P$ 
of $(X_t)_{t\in T}$ such that 
 $Y\in \sLx 2{P  \otimes \nu}$. Moreover, the
covariance function of $(Y)_{t\in T}$ is  $k_T^1$, and the $Z_i$'s satisfy \eqref{def-Zi}, i.e.~we
have 
\begin{displaymath}
   Z_i(\om) = \bigl\langle [Y(\om)]_\sim, [e_i]_\sim\bigr\rangle_{\Ltn}
\end{displaymath}
for $P$-almost all $\om\in \Om$ and all $i\in I$.
\end{theorem}

The last result in this section shows that the Karhunen-Lo{\`e}ve expansions we obtained are unique.
 Compared to Theorem \ref{def-X-by-Z-new}, this result does not need   Assumption K, i.e.~we do not need to know 
 the spectral properties of  $(e_i)_{i\in I}$ and $(\mu_i)_{i\in I}$,
while the convergence in \eqref{def-X-by-Z-h1} and the measurability of $X$ is now assumed.

\begin{theorem}\label{unique-kle}
   Let  $(T,\ca B, \nu)$ be a $\s$-finite measure space
for which $\Ltn$ is separable,
$(\Om, \ca A, P)$ be a probability space, and 
   $(X_t)_{t\in T} \subset \sLt P$ be a 
    centered and $(\ca A\otimes \ca B)$-measurable stochastic process. 
    Furthermore, assume that, for $I=\{1,\dots,n\}$ or $I=\N$, we have:
    \begin{enumerate}
       \item A family $(e_i)_{i\in I}$ of functions $T\to \R$ such that $([e_i]_\sim)_{i\in I}$ is an ONS in $\Ltn$.
       \item A family $(\mu_i)_{i\in I} \subset (0,\infty)$ converging monotonously to 0 and $\sum_{i\in I} \mu_i<\infty$.
       \item A family $( Z_i)_{i\in I}\subset \sLx 2 P$   of centered random variables satisfying \eqref{orth-z-again}.
    \end{enumerate}
  If, for all  $t\in T$, we know with convergence in $\Lx 2 P$ that
  \begin{equation}\label{KL-path-conv-Lt-again}
     [X_t]_\sim = \sum_{i\in I} [ Z_i]_\sim e_i(t)\, ,
  \end{equation}
   then we have  $X\in \sLt {P\otimes \nu}$, and 
    $(e_i)_{i\in I}$, $(\mu_i)_{i\in I}$, and $(Z_i)_{i\in I}$ are the families considered in Assumption X. In particular,
    \eqref{KL-path-conv-Lt-again} equals the Karhunen-Lo{\`e}ve expansion and we have \eqref{KL-path-conv-Lt}.
%
\end{theorem}


\section{Sample Paths Contained in Interpolation Spaces}\label{sec:path-in-interpol}

In this section we first characterize when 
the paths of the process are not only contained in $\Lx 2 \nu$ but actually in an 
interpolation spaces between $\Ltn$ and $H$.
In particular, it turns out that  stronger summability assumptions on the sequence $(\mu_i)_{i\in I}$
imply such path behavior, and in this case the average approximation error speed
of the Karhunen-Lo{\`e}ve  expansion  measured in the interpolation space
 can be exactly described 
by the behavior of $(\mu_i)_{i\in I}$.
Moreover, we will see that for Gaussian processes, the summability assumption is actually \emph{equivalent}
to the path behavior. Finally, we apply the developed theory to processes whose RKHS are contained 
in fractional Sobolev spaces, and consider small ball probabilities with respect to the interpolation spaces
considered above.

Let us begin with the following theorem, which characterizes when a single path is contained 
in a suitable interpolation space.

\begin{theorem}\label{interpol-new}
  Let Assumption X be satisfied, $\b\in (0,1)$,
  and $N\subset \Om$ be the measurable $P$-zero set we obtained from 
   Theorem \ref{pointwise-KL-new}.
  Then 
	for all 
	$\om\in \Om\setminus N$ and all finite $J\subset I$ we have 
	      \begin{equation}\label{in-interpol-char-lem-h1-new}
          \bnorm{\sum_{j\in J} Z_j(\om) [e_j]_\sim}_{[H]_\sim^{1-\b}}^2 
	  =\sum_{j\in J} \mu_j^{\b-1} Z_j^2(\om)\, .
      \end{equation}
	Moreover, for each $\om\in \Om\setminus N$, the following statements are equivalent: 
	\begin{enumerate}
	 \item We have $\sum_{i\in I}\mu_i^{\b-1}Z_i^2(\om)<\infty$.
	\item We have $[X(\om)]_\sim \in [H]_\sim^{1-\b}$.
	\item We have $[X(\om)]_\sim \in [\Lx 2 \nu , [H]_\sim]_{1-\b,2}$.
	\end{enumerate}
	Moreover, if one and thus all statements are true for a fixed $\om\in \Om\setminus N$,
	then \eqref{in-interpol-char-lem-h1-new} holds for all $J\subset I$, and
	the convergence in \eqref{KL-path-conv-Lt} is actually unconditional in the interpolation space
	$[\Lx 2 \nu , [H]_\sim]_{1-\b,2}$.
\end{theorem}

By Theorem \ref{interpol-new} we immediately see that almost all paths of the process $X$ are contained 
in the space $[\Lx 2 \nu , [H]_\sim]_{1-\b,2}$, if and only if
\begin{equation}\label{fourier-coeff-cond-interpol-new}
 \sum_{i\in I}\mu_i^{\b-1}Z_i^2(\om)<\infty
\end{equation}
for $P$-almost all $\om \in \Om$. Moreover, in this case 
the convergence in \eqref{KL-path-conv-Lt} is $P$-almost surely unconditional in the space $[\Lx 2 \nu , [H]_\sim]_{1-\b,2}$.
Note that in the case of $[\Lx 2 \nu , [H]_\sim]_{1-\b,2}\hookrightarrow \Lx \infty\nu$ the latter convergence
implies  
$\Lx \infty\nu$-convergence of the  Karhunen-Lo{\`e}ve Expansion in \eqref{KL-path-conv-Lt} for $P$-almost 
all $\om \in \Om$. In Corollary \ref{path-in-rkhs-cor2}, where we will consider this embedding 
situation again,
we will see that significantly more can be said.

To further illustrate Theorem \ref{interpol-new}, let 
 us fix an $\om \in \Om$ for which $[X(\om)]_\sim \in [H]_\sim^{1-\b}$ and \eqref{KL-path-conv-Lt} hold.
 Then, for all $\a\in [\b,1]$, we have both $[X(\om)]_\sim \in [H]_\sim^{1-\a}$ and
 \begin{displaymath}
  [X(\om)]_\sim = \sum_{i\in I} Z_i(\om) [e_i]_\sim = \sum_{i\in I} \mu_i^{(\a-1)/2}Z_i(\om) \, \bigl[\mu_i^{(1-\a)/2}e_i\bigr]_\sim\, .
 \end{displaymath}
Moreover, $([\mu_i^{(1-\a)/2}e_i]_\sim)_{i\in I}$ is an ONB of $[H]_\sim^{1-\a}$, and thus  we see  
that, for each $m\in I$, the sum
\begin{displaymath}
  \sum_{j=1}^m Z_j(\om) [e_j]_\sim
\end{displaymath}
 is the best approximation of $[X(\om)]_\sim$ in $[H]_\sim^{1-\a}$ \emph{for all}  $\a\in [\b,1]$
 \emph{simultaneously}.

Integrating \eqref{fourier-coeff-cond-interpol-new} with respect to $P$ and using \eqref{Zi-cov}
it is not hard to see that \eqref{fourier-coeff-cond-interpol-new} is $P$-almost surely satisfied, 
if $\sum_{i\in I} \mu_i^\b < \infty$. The following theorem characterizes this summability
in terms of the path behavior of the process.

\begin{theorem}\label{path-space}
 Let Assumption X be satisfied. Then, for $0<\b<1$, the following statements are equivalent:
 \begin{enumerate}
    \item We have $\sum_{i\in I} \mu_i^\b < \infty$.
    \item There exists an $N\in \ca A$ with $P(N) = 0$ such that $[X(\om)]_\sim \in [\Lx 2 \nu , [H]_\sim]_{1-\b,2}$ holds
    for all $\om \in \Om\setminus N$.
 Furthermore,   $\Om\setminus N\to [\Lx 2 \nu , [H]_\sim]_{1-\b,2}$
      defined by $\om\mapsto [X(\om)]_\sim$ is Borel measurable and we have 
\begin{displaymath}
\int_\Om\bnorm{[X(\om)]_\sim}_{[\Lx 2 \nu , [H]_\sim]_{1-\b,2}}^2 \, dP(\om) < \infty\, .
\end{displaymath}
 \end{enumerate}
Moreover, 
 there  exist constants $C_1,C_2>0$ such that, 
for all  $J\subset I$, we have 
\begin{equation*}
 C_1 \!\sum_{i\in  I\setminus J} \!\mu_i^\b 
\leq 
\int\limits_\Om\bnorm{\![X(\om)]_\sim - \sum_{j\in J} Z_j(\om) [e_j]_\sim}_{[\Lx 2 \nu , [H]_\sim]_{1-\b,2}}^2 \, dP(\om) 
\leq 
C_2 \!\sum_{i\in  I\setminus J} \!\mu_i^\b\, .
\end{equation*} 
\end{theorem}

In general, almost sure finiteness in \eqref{fourier-coeff-cond-interpol-new} is, of course, 
not equivalent to $\sum_{i\in I} \mu_i^{\b} < \infty$, since by \eqref{Zi-cov}
this summability describes $P$-integrability 
of the random variable in \eqref{fourier-coeff-cond-interpol-new}.
For Gaussian processes, however, we will see below 
that both conditions are in fact equivalent. The following lemma, 
which basically shows the equivalence of both notions under 
a martingale condition on $(Z_i^2)_{i\in I}$, is the key observation in this direction.

\begin{lemma}\label{sum-equiv}
   Let Assumption X be satisfied with $I=\N$. In addition, assume that, for all $i\geq 1$, we have
   $Z_i \in \sLx 4 P$ and 
      \begin{equation}\label{Z-martingale}
      \E_P (Z_{i+1}^2|\ca F_i) = \mu_{i+1}\, ,
   \end{equation}
   where $\ca F_i:= \s(Z_1^2,\dots, Z_i^2)$. Finally, assume that there exist constants $c>0$ and $\a\in (0,1)$
   such that 
   \begin{equation}\label{Z-martingale-var-bound}
      \var Z_i^2  \leq c \mu_i^{2-\a}
   \end{equation}
  for all $i\geq 1$. Then, for all $\b\in (\a,1)$, the following statements are equivalent:
  \begin{enumerate}
     \item We have $\sum_{i\in I} \mu_i^{\b} < \infty$.
     \item There exists an  $N\in \ca A$  with $P(N)=0$ such that 
     \eqref{fourier-coeff-cond-interpol-new} holds for all $\om \in \Om\setminus N$.
  \end{enumerate}
\end{lemma}

Note that \eqref{Z-martingale} is in particular satisfied if the random variables $(Z_i)_{i\in I}$
are independent. Moreover, \eqref{Z-martingale-var-bound} is satisfied, if and only if 
the 4th moments of the normalized variables $\xi_i := \mu_i^{-1/2}Z_i$ do not grow faster than $\mu_i^{-\a}$.

Combining the lemma above with Lemma \ref{gauss-lemma} we now obtain the 
announced equivalence for Gaussian processes.
 It further shows that either almost all or almost no paths  are contained in 
the considered interpolation space.

\begin{corollary}\label{GP-in-interpol}
   Let $(X_t)_{t\in T}$ be a Gaussian process for which Assumption X is satisfied. 
   Then, for $0<\b<1$, the following statements are equivalent:
 \begin{enumerate}
    \item We have $\sum_{i\in I} \mu_i^{\b} < \infty$.
    \item There exists an  $N\in \ca A$  with $P(N)=0$ such that 
    \eqref{fourier-coeff-cond-interpol-new} holds for all $\om \in \Om\setminus N$.
    \item There exists an  $A\in \ca A$  with $P(A)>0$ such that 
     $[X(\om)]_\sim \in [\Lx 2 \nu , [H]_\sim]_{1-\b,2}$ holds for all $\om \in A$.
 \end{enumerate}
  Moreover, all three statements are equivalent to  part ii) of Theorem \ref{path-space}.
\end{corollary}

So far, the developed theory is rather abstract. Our final goal in this 
section is to illustrate how our result can be used to investigate path
properties of certain families of processes. 
These considerations will be based on the following corollary, which, roughly speaking, 
shows that the sample paths of a process
are about $d/2$-less smooth than the functions in its RKHS.

\begin{corollary}\label{sobolev-in-interpol}
   Let $T\subset \Rd$  be a bounded subset 
that satisfies 
the strong local Lipschitz condition and $T=\interior \overline T$. Moreover, let
$\nu$ be the Lebesgue measure on $T$ and 
$(X_t)_{t\in T}$ be 
a stochastic process satisfying Assumption X.
Assume that   $H\hookrightarrow W^m(T)$
for some $m>d/2$. Then, for all $s\in (0,m-d/2)$, we have 
  \begin{equation}\label{sobolev-in-interpol-h1}
     [X(\om)]_\sim \in B^s_{2,2}(T)
  \end{equation}
for $P$-almost all $\om \in \Om$. Moreover, 
 there  exists a constant $C>0$ such that, 
for all  $J\subset I$, we have 
\begin{equation*}
\int\limits_\Om\bnorm{\![X(\om)]_\sim - \sum_{j\in J} Z_j(\om) [e_j]_\sim}_{ B^s_{2,2}(T)}^2 \, dP(\om) 
\leq 
C \sum_{i\in  I\setminus J} \!\mu_i^{1-s/m}\, ,
\end{equation*} 
and if  $H=W^m(T)$ with equivalent norms,  there exist constants $C_1,C_2>0$ such that for all $i\in I$ we have
\begin{equation*}
 C_1   i^{-\frac{2(m-s)}{d}+1}
\leq 
\int\limits_\Om\bnorm{\![X(\om)]_\sim - \sum_{j=1}^i Z_j(\om) [e_j]_\sim}_{ B^s_{2,2}(T)}^2 \, dP(\om) 
\leq 
C_2 i^{-\frac{2(m-s)}{d}+1}\, .
\end{equation*}
Finally, if $(X_t)_{t\in T}$  is a
Gaussian process with $H=W^m(T)$, then the results are sharp in the sense that 
 \eqref{sobolev-in-interpol-h1} does 
not  hold with strictly positive probability for $s:= m-d/2$.
\end{corollary}

Note that for Gaussian processes with $H=W^m(T)$, Corollaries \ref{sobolev-in-interpol}
and \ref{GP-in-interpol}
 show that \eqref{sobolev-in-interpol-h1} holds with some positive probability, 
if and only if, it holds with probability one, and the latter is also equivalent to $m>s+d/2$.

For general processes with $H=W^m(T)$ the smoothness exponent $s$ is 
also sharp in the sense that \eqref{sobolev-in-interpol-h1} does not hold for 
 $s:= m-d/2$ and $P$-almost all $\om\in \Om$, provided that the process satisfies the assumptions of 
Lemma \ref{sum-equiv} for some $\a\in (0, \frac d {2m})$.
The proof of this generalization is an almost literal copy of the proof of 
 Corollary \ref{sobolev-in-interpol}, and thus we decided to omit it.

Corollary \ref{sobolev-in-interpol} provides analytic properties of the sample paths in terms of $B^s_{2,2}(T)$,
whenever $H\hookrightarrow W^m(T)$ is known. 
For \emph{weakly stationary} processes the same result has been recently shown in 
\cite[Theorem 3 and Remark 1]{Scheuerer10a} by completely different techniques. 
Note that for such processes the inclusion $H\hookrightarrow W^m(T)$
can be easily checked using the Fourier transform of the kernel. We refer to \cite[Corollary 10.13]{Wendland05}
for the case $T=\R^d$ from which the general case can be easily deduced.

Let us now present some explicit examples for which 
an inclusion of the form  $H\hookrightarrow W^m(T)$ is known. 
We begin with a class of processes which include L\'evy processes.

\begin{example}\label{levy-example}
 Let $(X_t)_{t\in T}$ be 
a stochastic process satisfying Assumption X for $T=[0,t_0]$ and the Lebesgue measure $\nu$
on $T$. 
Furthermore, assume that 
the kernel is given by
\begin{displaymath}
 k(s,t) = \s^2 \cdot \min\{s,t\}\, ,\qquad \qquad s,t\in  [0,t_0],
\end{displaymath}
where $\s>0$ is some constant. 
 It is well-known, see e.g.~\cite[Example 8.19]{Janson97}, that the RKHS of this kernel 
 is continuously embedded into $W^1(T)$.  
Consequently,  for all $s\in (0,1/2)$, we have 
\begin{displaymath}
  [X(\om)]_\sim \in B^s_{2,2}(T)
\end{displaymath}
for $P$-almost all $\om \in \Om$. 

Note that the considered class of processes include 
L\'evy processes, and for these processes, it has been shown in  \cite{Herren97a}
by different means
that their paths are also contained in $B^s_{p,\infty}(T)$  for all $s\in (0,1/2)$ and $p>2$ with 
$sp<1$. 
Interestingly,  this is
is equivalent to our result above.
Indeed, if we fix a pair of $s$ and $p$ satisfying the assumptions of 
 \cite{Herren97a}, we have 
$s_0:= s- 1/p +   1/2 < 1/2$. For $\e>0$ with $s_0+\e<1/2$ we then obtain $s_0+\e -   1/2 > s-1/p$
and $s_0+\e>s$
and thus 
\begin{displaymath}
   B^{s_0+\e}_{2,2}(T) \hookrightarrow B^s_{p,\infty}(T)
\end{displaymath}
by \cite[p.~82]{RuSi96}.
Consequently, our result implies that of \cite{Herren97a}.
Conversely, if we fix an $s<1/2$, there is an $\e>0$ with $s+2\e<1/2$, and for 
$s_0:= s+\e$ and $p_0:= (s+2\e)^{-1}$, we have $s_0>s$ and $s_0-1/p_0>s-1/2$, so that 
\begin{displaymath}
   B^{s_0}_{p_0,\infty}(T) \hookrightarrow B^s_{2,2}(T)
\end{displaymath}
by \cite[p.~82]{RuSi96}.
Since we also have $s_0<1/2$, $p_0>2$ and $s_0p_0<1$, we then see that the result of
\cite{Herren97a}   implies ours.

Finally, for the 
Brownian motion, it is well-known that there exists a version whose sample paths are 
contained in $B^s_{\infty,\infty}(T)$ for all $s\in (0,1/2)$, and finer results 
can be found in \cite{Roynette93a}.
\end{example}

The following example includes the Ornstein-Uhlenbeck processes. Note that 
although the kernel in this example looks quite different to the one of Example \ref{levy-example}
the results on the smoothness properties of the paths are identical.

\begin{example}\label{ornstein-uhlenbeck-example}
 Let $(X_t)_{t\in T}$ be 
a stochastic process satisfying Assumption X for $T=[0,t_0]$ and the Lebesgue measure $\nu$
on $T$. Furthermore, assume that 
the kernel is given by
\begin{displaymath}
 k(t,t') = ae^{-\s |t-t'|}\, ,\qquad \qquad t,t'\in  [0,t_0],
\end{displaymath}
where $a,\s>0$ are some constants. 
 It is well-known, see e.g.~\cite[p.~316]{BeTA04} and \cite[Example 5C]{Parzen61a}, that the RKHS of this kernel 
 equals $W^1(T)$ up to equivalent norms. Consequently, 
  for all $s\in (0,1/2)$, we have 
\begin{equation}\label{ousp}
  [X(\om)]_\sim \in B^s_{2,2}(T)
\end{equation}
for $P$-almost all $\om \in \Om$. Note that the considered class of processes include a specific form of
the Ornstein-Uhlenbeck process, see \cite[Example 8.4]{Janson97}. 

By subtracting
the one-dimensional $C^\infty$-kernel $(t,t')\mapsto ae^{-\s (t+t')}$ from $k$, we see that 
\eqref{ousp} also holds for processes having the kernel 
\begin{displaymath}
   \tilde k(t,t') = ae^{-\s |t-t'|} - ae^{-\s (t+t')}\, ,\qquad \qquad t,t'\in  [0,t_0].
\end{displaymath}
Recall that the classical Ornstein-Uhlenbeck processes belong to this class of processes.
\end{example}

The following example considers processes on higher dimensional domains with potentially
smoother sample paths. It is in particular interesting for certain 
statistical methods, see \cite{Stein99,RaWi06,PiWuLiMuWo07a,VaZa08a,GnKlSc10a,VaZa11a,ScScSc13a},
since the considered family of covariance functions allows for a high flexibility in these methods.
Moreover, 
note that for $d=1$ and $\a=1/2$ the previous example is recovered.

\begin{example}\label{matern-example}
Let $T\subset \Rd$ be an open and bounded subset 
satisfying the strong local Lipschitz condition
and $\nu$ be the Lebesgue measure on $T$.
Furthermore, let
 $(X_t)_{t\in T}$ be 
a stochastic process satisfying Assumption X.
Assume that 
its covariance is a Mat\'ern  kernel of order $\a>0$, that is 
\begin{displaymath}
 k(s,t) = a \bigr(\s\snorm{s-t}_2\bigr)^\a H_\a \bigl( \s\snorm{s-t}_2 \bigr) \, ,\qquad \qquad s,t\in  T,
\end{displaymath}
where $a, \s>0$ are some constants and $H_\a$ denotes the modified Bessel function of the second type  
of order $\a$. Then up to equivalent norms
the RKHS $H_{\a,\s}(T)$ of this kernel is $B^{\a+d/2}_{2,2}(T)$,
see \cite[Corollary 10.13]{Wendland05} together with \cite[Theorem 5.3]{ScWe01a}, as well as
\cite{BoRoSc13a} for a generalization.
Consequently, 
for all $s\in (0,\a)$, we have 
\begin{displaymath}
  [X(\om)]_\sim \in B^s_{2,2}(T)
\end{displaymath}
for $P$-almost all $\om \in \Om$.

For $d=1$ and $\a=k+r$ with $k\in \N$ and $r\in (1/2,1]$, it was shown in \cite{CrLe67}, cf.~also \cite{HaSt93a},
that there exists a version of
the process with $k$-times continuously differentiable 
paths. 
Our result  improves this. Indeed, for $d$ and $\a$ as above, 
we clearly find an $s\in (0,\a)$ with $s-k>1/2$ and since, for this 
$s$, we have 
$B^s_{2,2}(T)  \hookrightarrow C^{k}(T)$, see e.g.~\cite[Theorem 8.4]{UlrichXXa}, we see that 
$P$-almost all paths 
$X(\om)$ equal $\nu$-almost everywhere a $k$-times continuously differentiable function.
We will show in the next section that there actually exist
a version $(Y_t)_{t\in T}$ of the process with $Y(\om) \in  B^s_{2,2}(T)$ almost surely, 
so that our result does improve 
  the above mentioned  classical result in \cite{CrLe67}.
\end{example}

Theorem \ref{path-space} and its Corollary \ref{sobolev-in-interpol}
showed that the average interpolation space norm of the process's paths 
is finite. In the following we illustrate how our results, in particular 
\eqref{in-interpol-char-lem-h1-new} can be used to estimate  small ball probabilities.
For the sake of simplicity, we will mostly restrict our considerations to Gaussian processes, 
but at the end of the discussion we will briefly indicate possible generalizations.


Let us begin by recalling that the small ball problem considers 
probabilities of the form
\begin{displaymath}
   P\bigl( \{ \om \in \Om :\snorm{X(\om)} \leq \e \}  \bigr)
\end{displaymath}
for $\e\to 0^+$ and various norms including the two standard examples $\inorm \cdot$ and 
$\snorm\cdot_{\Ltn}$. Recall that for centered continuous Gaussian processes
on compact $T$
 with continuous
covariance function
there is a nice link between 
$\inorm \cdot$-small ball probabilities and entropy numbers of the inclusion $H\to C(T)$ via 
Gaussian measures, which was 
discovered by \cite{KuLi93a}, see also \cite[Theorem 3.3]{LiSh01a} in combination with 
\cite[Example 2.4, p.~24, \& p.~33]{Lifshits12}.
Let us further recall that $\Ltn$-small ball probabilities can be estimated 
using the Karhunen-Lo{\`e}ve expansion of the process. 
To this end, we write $\xi_i := \mu_i^{-1/2}Z_i$ for all $i\in I$. Then $(\xi_i)_{i\in I}$
is an ONS in $\sLx 2 P$ by \eqref{Zi-cov}, and 
\eqref{KL-path-conv-Lt} together with Parseval's identity gives 
\begin{equation}\label{basic-kle-small-ball}
   \snorm {[X(\om]_\sim}_\Ltn^2 = \sum_{i\in I} Z_i^2(\om) = \sum_{i\in I} \mu_i\, \xi_i^2(\om)
\end{equation}
for $P$-almost all $\om \in \Om$. Consequently, it suffices to understand the small 
ball behavior of the right-hand side of \eqref{basic-kle-small-ball}.
To proceed, we now assume that $X$ is a Gaussian process satisfying Assumption X.
Then the random variables $(\xi_i)_{i\in I}$ are i.i.d.~with 
$\xi_i \sim \ca N(0,1)$ by Lemma \ref{gauss-lemma}, and 
right-hand side of \eqref{basic-kle-small-ball} is somewhat well-understood, see the references in  
\cite[Section 6]{LiSh01a}. In particular, if $\mu_i \sim i^{-\a}$ for some $\a>1$, i.e.~$I=\N$ and $\lim_{i\to \infty} i^\a \mu_i = 1$,
then \cite[Theorem 4.2]{Aurzada07a} describes  the exact  small ball behavior 
of the right-hand side of \eqref{basic-kle-small-ball}.

Obviously, the same arguments work if we replace 
\eqref{basic-kle-small-ball}
by \eqref{in-interpol-char-lem-h1-new}. The following theorem
presents the corresponding result.

\begin{corollary}\label{cor:small-ball-1}
 Let $(X_t)_{t\in T}$ be a Gaussian process satisfying Assumption $X$. Furthermore, assume that 
 the countably many eigenvalues satisfy $\mu_i \sim i^{-\a}$ for some $\a>1$. Then for all $\b\in (0,1]$ 
 satisfying $\a\b>1$, there exists a constant $C_\b\in (0,\infty)$ such that 
 \begin{displaymath}
  \lim_{\e\to 0^+} \e^{\frac 2 {\a\b-1}} \log P\bigl(\bigl\{\om\in \Om: \mnorm{[X(\om]_\sim}_{[H]_\sim^{1-\b}} \leq \e  \bigr\}  \bigr) = -C_\b
 \end{displaymath}
 Moreover, if we only have $\mu_i \asymp i^{-\a}$, then, for $\e\to 0^+$, we obtain
 \begin{equation}\label{cor:small-ball-1-h1}
 - \log P\Bigl(\bigl\{\om\in \Om: \mnorm{[X(\om]_\sim}_{[\Ltn, [H]_\sim]_{1-\b,2}} \leq \e  \bigr\}  \Bigr) \asymp \e^{-\frac 2 {\a\b-1}}\, .
 \end{equation}
 Finally, in the cases of $\mu_i\preceq i^{-\a}$, respectively $\mu_i\succeq i^{-\a}$, we only have 
 ``$\preceq$'', respectively ``$\succeq$'' in \eqref{cor:small-ball-1-h1}.
\end{corollary}

As indicated above, 
the case $\b=1$ in Corollary \ref{cor:small-ball-1} reproduces 
well-known $L_2$-small ball estimates, which can be established using  \eqref{basic-kle-small-ball}.
Links to such bounds can be found in 
\cite[Section 6]{LiSh01a} and \cite[p.~94]{Lifshits12}.
Corollary \ref{cor:small-ball-1} shows that this technique can be 
extended to interpolation space norms without any technical hurdles.

Note that the constant $C_\b$ above can be explicitly calculated, see \cite[Theorem 1.1]{Aurzada07a}.
Moreover, similar results can be obtained if $\mu_i$ behaves like $i^{-\a} J(i)$, where $J(i)$
is a slowly varying sequence, see again  \cite[Theorem 1.1]{Aurzada07a} in the case of $J(i) = (\log i)^\g$ and
\cite{BoRu08a} for the general case. In this regard, we also note that the Gaussianity in 
Corollary \ref{cor:small-ball-1} is not necessary. Indeed, the results from \cite{Aurzada07a,BoRu08a}
remain valid if  $(\xi_i)_{i\in I}$ are i.i.d.~and $\xi_1$ has a continuous density that does not vanish at $0$.

The following corollary applies Corollary \ref{cor:small-ball-1} to the situation considered in 
Corollary \ref{sobolev-in-interpol}.

\begin{corollary}\label{cor:small-ball-2}
   Let $T\subset \Rd$  be a bounded subset 
that satisfies 
the strong local Lipschitz condition and $T=\interior \overline T$. Moreover, let
$\nu$ be the Lebesgue measure on $T$ and 
$(X_t)_{t\in T}$ be 
a Gaussian process satisfying Assumption X.
Assume that   $H\hookrightarrow W^m(T)$
for some $m>d/2$. Then, for all $s\in (0,m-d/2)$, we have 
  \begin{equation}\label{cor:small-ball-2-h1}
      - \log P\Bigl(\bigl\{\om\in \Om: \mnorm{[X(\om]_\sim}_{B^s_{2,2}(T)} \leq \e  \bigr\}  \Bigr) \preceq \e^{-\frac {2d} {2m-2s-d}}\, .
  \end{equation}
Moreover, if $H=W^m(T)$ with equivalent norms, then we have ``$\asymp$'' in \eqref{cor:small-ball-2-h1}.
\end{corollary}

To illustrate the last corollary, let us consider Examples \ref{levy-example} 
and \ref{ornstein-uhlenbeck-example}. In these examples we have $m=d=1$ and thus 
\eqref{cor:small-ball-2-h1} holds for all $s\in (0,1/2)$ where the exponent on the right hand-side 
of \eqref{cor:small-ball-2-h1} is given by $\frac 2 {1 -2s}$.
Note that for the Wiener process considered in Example \ref{levy-example} 
this has already been proved by \cite{Stolz93a,LiSh99a} in a more general context.
Moreover, in Example \ref{matern-example} we have $m=\a+d/2$ and $H=W^m(T)$, so that 
\eqref{cor:small-ball-2-h1} holds with ``$\asymp$'' and exponent 
$\frac d {\a-s}$ for all $s\in (0,\a)$.

%

\section{Sample Paths Contained in RKHSs}\label{sec:paths-in-rkhs}

So far we have seen that, under some summability assumptions,
the   $\nu$-equivalence classes of the process are contained in a suitable
interpolation space. Now recall from Section  \ref{sec:prelim}
that these interpolation spaces can sometimes be  viewed as RKHSs, too.
The goal of this section is   to present conditions under which
a suitable version of the process has actually its paths in this RKHS.
In particular, we will see that under stronger summability conditions on the 
eigenvalues such a path behavior occurs, in a certain sense, automatically.


Let us begin by fixing the following set of assumptions, which in particular ensure 
that $k_S^{1-\b}$ can be constructed.

\begin{assumptionsks}
  Let Assumption K be satisfied. Moreover, let $0<\b<1$ and $S\subset T$ be a measurable set 
	with $\nu(T\setminus S)=0$
	such that, 
   for all $t\in S$, we have 
  \begin{align}\label{k-eq-kn}
   \sum_{i\in I} \mu_i e_i^2(t) & = k(t,t)\\ \label{k-series-diag-1-b}
 \sum_{i\in I} \mu_i^{1-\b} e_i^2(t) & < \infty\, .
\end{align} 
\end{assumptionsks}

Note that if $H$ is separable, we can always find a set $S$ of full measure $\nu$ for which
\eqref{k-eq-kn} holds, see \cite[Corollary 3.2]{StSc12a}. For such $H$, Assumption KS thus reduces to assuming 
that we can construct $\kxy S{1-\b}$, and the latter is possible, if, e.g.~$\sum_{i\in I}\mu_i^{1-\b}<\infty$,
see \eqref{when-kS}.
Moreover, recall from Lemma \ref{Ik-inject} that \eqref{k-eq-kn} holds for $S=T$ if  
Assumption CK is satisfied. Finally, if,  in addition,  we have 
 $[\Ltn , [H]_\sim]_{1-\b,2}\hookrightarrow \Lx \infty {\nu}$, then Theorem \ref{bounded-powers-of-cont-kernels}
 shows that \eqref{k-series-diag-1-b} also holds for $S=T$.

Our first result   characterizes when a suitable version of our process $(X_t)_{t\in T}$
has its paths in the corresponding RKHS $\Hxy S{1-\b}$.

\begin{theorem}\label{path-in-rkhs-1}
   Let Assumptions X and KS be satisfied. 
Then the following statements are equivalent:
\begin{enumerate}
   \item There exists a measurable $N\subset \Om$  with $P(N)=0$ such that for all $\om \in \Om\setminus N$ we have
    \begin{equation}\label{fourier-coeff-cond}
      \sum_{i\in I}  \mu_i^{\b-1}Z_i^2(\om) < \infty\, .
    \end{equation}
\item  There exists an $(\ca A\otimes \ca B)$-measurable version $(Y_t)_{t\in T}$ of $(X_t)_{t\in T}$  
    such that,  for $P$-almost all $\om \in \Om$, we have 
      \begin{equation}\label{path-in-rkhs-1-h1}
         Y(\om)_{|S} \in \Hxy S{1-\b}\, .
      \end{equation}
\end{enumerate}
Moreover, if one and thus both statements are true, we have for $P$-almost all $\om \in \Om$
\begin{equation}\label{path-in-rkhs-1-h2}
   Y(\om)_{|S}  = \sum_{i\in I} Z_i(\om) (e_i)_{|S}\, ,
\end{equation}
where the convergence is unconditional in $\Hxy S{1-\b}$.
\end{theorem}

If \eqref{fourier-coeff-cond} is $P$-almost surely satisfied then 
Theorem \ref{path-in-rkhs-1} strengthens Theorem \ref{interpol-new} in the sense that
$[X(\om)]_\sim \in [H]_\sim^{1-\b}$ is replaced by   $Y(\om)_{|S} \in \Hxy S{1-\b}$.
Moreover, unlike \eqref{in-interpol-char-lem-h1-new}, which only gives 
$[H]_\sim^{1-\b}$-convergence of the  Karhunen-Lo{\`e}ve Expansion in \eqref{KL-path-conv-Lt},
the expansion \eqref{path-in-rkhs-1-h2} converges in $\Hxy S{1-\b}$, which in particular
implies pointwise convergence at all $t\in S$, since $\Hxy S{1-\b}$ 
is an RKHS.

We already know that the Fourier coefficient condition \eqref{fourier-coeff-cond} 
can be ensured by a summability condition on the eigenvalues. Like in Theorem \ref{path-space},
this summability can be characterized by the path behavior of the version $(Y_t)_{t\in T}$
as the following theorem shows.

\begin{theorem}\label{path-in-rkhs-2}
  Let Assumptions X and KS be satisfied.
Then the following statements are equivalent:
\begin{enumerate}
 \item We have $\sum_{i\in I}  \mu_i^{\b}  < \infty$.
 \item We have $\kxy S1 \ll \kxy S{1-\b}$. 
	\item   There exists an $(\ca A\otimes \ca B)$-measurable version $(Y_t)_{t\in T}$ of $(X_t)_{t\in T}$  
    such that,  for $P$-almost all $\om \in \Om$, we have 
     $Y(\om)_{|S} \in \Hxy S{1-\b}$, and 
	\begin{equation}\label{path-in-rkhs-norm-est}
	 \int_\Om \snorm{Y(\om)_{|S}}_{\Hxy S{1-\b}}^2 \dxy P\om < \infty\, .
	\end{equation}
\end{enumerate}
\end{theorem}

Let us compare the previous two theorems in the case of $S=T$
with the results
of Luki{\'c} and   Beder in  \cite{LuBe01a}.
Their Theorem 5.1 shows that $\kxy T1 \ll \kxy T{1-\b}$ implies
\eqref{path-in-rkhs-1-h1},
and,  their Corollary 3.2 conversely shows that 
\eqref{path-in-rkhs-norm-est}  implies $\kxy T1 \ll \kxy T{1-\b}$.
Clearly, the difference between these two implications is exactly the difference between 
\eqref{path-in-rkhs-norm-est} and \eqref{path-in-rkhs-1-h1}, and this difference is exactly described 
by Theorems \ref{path-in-rkhs-1} and \ref{path-in-rkhs-2}.
While in this sense, the latter two theorems   clarified the situation for the space $\Hxy T{1-\b}$,
 it seems fair to say that  the less exact results 
 in  \cite{LuBe01a} are more general as \emph{arbitrary} RKHS $\bar H$  satisfying $H\hookrightarrow \bar H$
 are considered.

 The following corollary, which considers the case of Gaussian processes,
 basically recovers the findings of \cite[Section 7]{LuBe01a}.
We mainly state it here  for the sake of completeness.

\begin{corollary}\label{GP-in-RKHS}
       Let $(X_t)_{t\in T}$ be a Gaussian process 
	for which Assumptions X and KS are satisfied.
    Then the following statements are equivalent:
    \begin{enumerate}
      \item We have $\sum_{i\in I}  \mu_i^{\b}  < \infty$.
      \item We have $\kxy S1 \ll \kxy S{1-\b}$. 
      \item There exists an $(\ca A\otimes \ca B)$-measurable version $(Y_t)_{t\in T}$ of $(X_t)_{t\in T}$  
    such that,  for $P$-almost all $\om \in \Om$, we have 
      \begin{displaymath}
         Y(\om)_{|S} \in \Hxy S{1-\b}\, .
      \end{displaymath}
            \item There exists an $(\ca A\otimes \ca B)$-measurable version $(Y_t)_{t\in T}$ of $(X_t)_{t\in T}$ 
            and an $A\in \ca A$ with $P(A)>0$ 
    such that,  for all $\om \in A$, we have 
      \begin{displaymath}
         Y(\om)_{|S} \in \Hxy S{1-\b}\, .
      \end{displaymath}
    \end{enumerate}
\end{corollary}

For general processes satisfying the assumptions made in Lemma \ref{sum-equiv} for some $\a\in (0,1)$,
the equivalences \emph{i)} $\Leftrightarrow$ \emph{ii)} $\Leftrightarrow$ \emph{iii)} of Corollary \ref{GP-in-RKHS}
also hold  for all $\b\in (\a,1)$. Indeed, the implication \emph{iii)} $\Rightarrow$ \emph{i)}
can be shown by Lemma \ref{sum-equiv}, and the remaining implications actually do not require 
the Gaussian assumption at all.

If we wish to find an RKHS $\bar H$ that contains the paths of a suitable version of the process by
the results presented so far, we need to know   the eigenvalues and eigenfunctions
as well as the interpolation spaces \emph{exactly}.
However,  
obtaining the exact eigenvalues and -functions of $T_k$ is often a very difficult, if not impossible, task,
and the interpolation spaces may not be readily available, either.
The following two corollaries address this issue by presenting a sufficient condition for 
the existence of such an RKHS $\bar H$.

\begin{corollary}\label{path-in-rkhs-cor1}
   Let Assumption X be satisfied, $H$ be separable, and $\bar H$ be an RKHS on $T$ with kernel $\bar k$ such that 
    $H\hookrightarrow \bar H$.  Let us further assume that $\bar H$ is 
   compactly embedded into $\Ltn$ and  that 
   \begin{equation}\label{entropy-sum}
        \sum_{i=1}^\infty \e_i^{\a}(I_{\bar k})<\infty
   \end{equation}
  for some $\a\in (0,1]$. Then, for all $\b\in [\a/2,1-\a/2]$, there exists   a measurable  $S\subset T$
	with $\nu(T\setminus S)=0$ such that the following statements are true:
	\begin{enumerate}
	   \item Both $\Hxy S{1-\b}$ and $\bbHxy S{1-\b}$ exist, and we have
	$\Hxy S{1-\b} \hookrightarrow \bbHxy S{1-\b}$.
	  \item There exists an $(\ca A\otimes \ca B)$-measurable version $(Y_t)_{t\in T}$ of $(X_t)_{t\in T}$  
    such that  
     $Y(\om)_{|S} \in \Hxy S{1-\b}$ for $P$-almost all $\om \in \Om$, and \eqref{path-in-rkhs-norm-est} holds. 
	\end{enumerate} 
\end{corollary}

Corollary \ref{path-in-rkhs-cor1} shows that 
 in order to construct an RKHS containing paths on a set $S$ of full measure $\nu$ we do not necessarily 
	need to know the eigenvalues and -functions exactly. Instead, it suffices to have 
	an RKHS $\bar H$ with $H\hookrightarrow \bar H$ for which we know both, entropy number estimates 
	of the map $I_{\bar k}$ and the interpolation spaces of $\bar H$ with $\Ltn$.
	Namely, if \eqref{entropy-sum} is satisfied, then the version $(Y_t)_{t\in T}$
	obtained by Corollary \ref{path-in-rkhs-cor1} satisfies 
	$Y(\om)_{|S} \in \bbHxy S{1-\b}$ for $P$-almost all $\om \in \Om$ and combining this with 
	$\Hxy S{1-\b} \hookrightarrow \bbHxy S{1-\b}$ we see that we have $\bbHxy S{1-\b}$-convergence in 
	\eqref{path-in-rkhs-1-h2}. Similarly to Corollary \ref{sobolev-in-interpol} it is further possible
	to upper bound the average speed of $\bbHxy S{1-\b}$-convergence in  \eqref{path-in-rkhs-1-h2}
	with the help of the 
	entropy numbers of $I_{\bar k}$. We omit the details for the sake of brevity.
	Moreover,   $Y(\om)_{|S} \in \Hxy S{1-\b}\subset \bbHxy S{1-\b}$ for $P$-almost all $\om \in \Om$ 
	shows that 
	$P$-almost all paths  also enjoy a representation  of the form
	\begin{displaymath}
	    Y(\om)_{|S}  = \sum_{j\in J} \bar Z_j(\om) \bar e_j\, ,
	\end{displaymath}
	where the convergence is unconditional in   $ \bbHxy S{1-\b}$, $(\bar e_j)_{j\in J}$ is the family 
	obtained by Theorem \ref{spectral} for the operator $T_{\bar k_S^1}$ and $(\bar Z_j)_{j\in J}$
	is a suitable family of random variables such that 
	$\sum_{j\in J}\bar \mu_j^{\b-1} \bar Z_j^2(\om)<\infty$ for $P$-almost all $\om \in \Om$.
	Following the logic above, this representation may be easier at hand than the 
	standard Karhunen-Lo{\`e}ve expansion, but its deeper investigation is beyond the scope of this paper.

The following corollary provides a result in the same spirit for the case $S=T$.
In particular, it provides two sufficient conditions under which 
there exists an RKHS containing almost all paths of a suitable version.
This answers a question raised in \cite{Lukic04a}.

\begin{corollary}\label{path-in-rkhs-cor2}
   Let Assumption X be satisfied, $H$ be separable,   
   and $\bar H$ be an RKHS on $T$ with kernel $\bar k$ such that both
    $H\hookrightarrow \bar H$ and  $\bar H$ is compactly embedded into $\Ltn$.
    Furthermore, assume that  $(T,\ca B,\nu)$ and   $\bar k$ satisfy Assumption CK, and that, 
     for  some $\b\in (0,1/2]$, one of following assumptions are satisfied:
    \begin{enumerate}
        \item The
       eigenfunctions $(\bar e_j)_{j\in J}$ of   $T_{\bar k}$ are uniformly bounded, 
     i.e.~$\sup_{j\in J}\inorm{\bar e_j}<\infty$, and 
        we have
        \begin{displaymath}
          \sum_{i=1}^\infty \e_i^{2\b}(I_{\bar k})<\infty, .
        \end{displaymath}
      \item We have $[\Ltn , [\bar H]_\sim ]_{1-\b,2}\hookrightarrow \Lx \infty {\nu}$.
    \end{enumerate}
     The the following statements hold:
     \begin{enumerate}
      \item The kernels $k_T^{1-\b}$ and $\bar k_T^{1-\b}$  exist,  are bounded, and we have 
      $\Hxy T{1-\b} \hookrightarrow \bbHxy T{1-\b}$.
      \item There exists an $(\ca A\otimes \ca B)$-measurable version $(Y_t)_{t\in T}$ of $(X_t)_{t\in T}$  
	such that   
      $Y(\om)  \in \Hxy T{1-\b}$ for all $\om \in \Om$, and \eqref{path-in-rkhs-norm-est} holds. 
      \item All paths of $Y$ are bounded and $\t(H)$-continuous.
      \item For $P$-almost all $\in \Om$, the expansion \eqref{path-in-rkhs-1-h2} converges uniformly in $t$
      on $S=T$.
      \item If there is a separable and metrizable topology $\t$ on $T$ such that $\t(H)\subset \t$ and
       almost all paths of $X$ are $\t$-continuous, then  $X(\om) = Y(\om)$ for $P$-almost all 
      $\om \in \Om$. In particular, this holds if almost all paths of $X$ are $\t(H)$-continuous
      and $\t(H)$ is Hausdorff.
     \end{enumerate}
\end{corollary}

Note that in the situation of part \emph{iv)} of Corollary \ref{path-in-rkhs-cor2}
the 
Karhunen-Lo{\`e}ve Expansion in \eqref{KL-path-conv-Lt} 
converges in $\ell_\infty(T)$ 
for $P$-almost 
all $\om \in \Om$. Moreover, note the $\t(H)$-continuity of the paths obtained in 
\emph{iii)} and \emph{iv)} is potentially stronger than the $\t$-continuity,
where $\t$ is a ``natural'' topology of $T$.

The last result of this section improves Corollary \ref{sobolev-in-interpol}.
Note that it directly applies to the processes considered in Example \ref{matern-example}.

\begin{corollary}\label{sobolev-in-rkhs}
   Let $T\subset \Rd$  be a bounded subset 
that satisfies 
the strong local Lipschitz condition and $T=\interior \overline T$. Moreover, let
$\nu$ be the Lebesgue measure on $T$ and 
$(X_t)_{t\in T}$ be 
a stochastic process satisfying Assumption X.
Assume that   $H\hookrightarrow W^m(T)$
for some $m>d$. Then the following statements hold:
\begin{enumerate}
   \item For all $s\in (d/2,m-d/2)$, 
 there exists an $(\ca A\otimes \ca B)$-measurable version $(Y_t)_{t\in T}$ of $(X_t)_{t\in T}$   such that, for all $\om \in \Om$, we have  
  \begin{equation}\label{sobolev-in-rkhs-h1}
     Y(\om) \in B^s_{2,2}(T)\, .
  \end{equation}
  Moreover, for $P$-almost all $\om \in \Om$  we have with  unconditional convergence in  $ B^s_{2,2}(T)$:
\begin{equation}\label{sobolev-in-rkhs-h2}
   Y(\om)   = \sum_{i\in I} Z_i(\om) e_i \, .
\end{equation}
  \item If $(X_t)_{t\in T}$  is a Gaussian process with $H=W^m(T)$, then the results are sharp in the sense that 
 \eqref{sobolev-in-rkhs-h1} does 
not  hold with strictly positive probability for $s:= m-d/2$.
\end{enumerate}
\end{corollary}

By \cite[Theorem 7.37]{AdFo03}, we immediately see that the convergence in \eqref{sobolev-in-rkhs-h2}
is uniform in $t$. Moreover, if $s>k+1/2$ for some $k\in \N$, then the convergence is also 
in $C^k(T)$, see e.g.~~\cite[Theorem 8.4]{UlrichXXa}.

Finally, like for Corollary \ref{sobolev-in-interpol}, 
the sharpness result in \emph{ii)} can be extended to a broader class of processes. We refer to our remarks 
following Corollary \ref{sobolev-in-interpol}.

%
%
%

\section{Final Remarks}\label{sec:final-remaks}

Summarizing our findings, we see that the following objects, which 
describe the relationship between $H$ and $\nu$, are crucial for our approach:
\emph{a)} the eigenvalues and -functions of the integral operator $T_k:\Ltn\to \Ltn$, and 
\emph{b)} the interpolation spaces between $H$ and $\Ltn$. 
Here $H$, respectively $k$, is given to us by the considered stochastic process, whereas
$\nu$ can be chosen by us. For example, when $H$ is embedded into a Sobolev space,
we could consider measures of the form $wd\lb$,
where $\lb$ is the Lebesgue measure and $w\geq 0$ is an ``interesting'' weight, see e.g.~\cite{DeMa03a,BaIg11a} for a few specific examples.
This observation indicates that a closer investigation of the relationship between $H$ and $\nu$,
for more general $\nu$ than classical results have focused on, could be fruitful.

Interestingly, such an investigation would  also be interesting for  quite a different reason. 
Indeed, the mathematical theory 
of one of the most successful machine learning algorithms of the last two decades, namely
support vector machines (SVMs), requires knowledge about the relationship between $H$ and $\nu$, too.
For example, if a least squares loss is used in SVMs, then their learning ability can be exactly described
by the eigenvalues of $T_k$ and the interpolation spaces between $H$ and $\Ltn$, see \cite{SmZh03a,SmZh07a,StCh08,StHuSc09b}
for this and related results.
Unlike in this paper, however, the machine learning setting allows us to pick $H$, while $\nu$ is 
given to us by the data generating distribution. This setting thus naturally  demands for considering more general $\nu$.

In conclusion, it seems fair to say that a deeper investigation of the relationship between RKHSs $H$ and (probability) measures $\nu$
would influence both our understanding of certain aspects of stochastic processes and  
the mathematical theory of one of the state-of-the-art learning algorithms.

\section{Proofs}\label{sec:proofs}

\subsection{Some Auxiliary Results and Proofs of Section \ref{sec:prelim}}

Our first result investigates isometries between some spaces of the form $\bHSb$, $\HSb$, and $[H]^\b_\sim$.

\begin{lemma}\label{basic-isometries}
    Let Assumption K be satisfied,  $\b\in (0,1]$, and  $R\subset S\subset T$ be measurable subsets such that $R$ satisfies 
$\nu(T\setminus R)=0$ and \eqref{k-series-diag}. Then the following operators are isometric isomorphisms:
\begin{enumerate}
 \item The multiplication operator $\eins_R:   \bHSb\to\bHxy R\b$ defined by $f\mapsto \eins_Rf$.
	\item The zero-extension operator $\,\hat{\cdot}:\HSb\to \bHSb$.
	\item The restriction operator $\cdot_{|R}: \bHSb\to \Hxy R\b$.
	\item The equivalence-class operator $[\,\cdot\,]_\sim : \bHSb\to [H]^\b_\sim$.
\end{enumerate}
\end{lemma}

\begin{proofof}{Lemma \ref{basic-isometries}}
\ada i Let us pick an $f\in \bHSb$. Then there exists a sequence $(a_i)\in \ell_2(I)$
	such that $f= \sum_{i\in I} a_i\mu_i^{\b/2}\eins_S e_i$, where the convergence is 
	in $\bHSb$ and thus also pointwise. Consequently, we find
	\begin{displaymath}
	 \eins_R f 
	= \eins_R \sum_{i\in I} a_i \mu_i^{\b/2}\eins_S e_i 
	=  \sum_{i\in I} a_i \mu_i^{\b/2}\eins_R \eins_S e_i 
	=  \sum_{i\in I} a_i \mu_i^{\b/2}\eins_R  e_i \, .
	\end{displaymath}
	Now the assertion easily follows from the definitions of the spaces $\bHSb$ and $\bHxy R\b$.

\ada {ii} Can be shown analogously to \emph{i).}

\ada {iii} Again, this can be shown analogously to \emph{i).}

\ada {iv} We obviously have 
$[\hat e_i]_\sim = [e_i]_\sim$ for all $i\in I$.
Moreover, for $(a_i)\in \ell_2(I)$ we have 
\begin{displaymath}
 \Bigl[ \sum_{i\in I} a_i \mu_i^{\b/2} \hat e_i\Bigr]_\sim
=
\sum_{i\in I}  a_i \mu_i^{\b/2} [\hat e_i]_\sim
\end{displaymath}
with convergence in $\Ltn$ by the continuity of $I_{\bkSb}:\bHSb\to \Ltn$. 
Combining both with the definition
 of the spaces $\bHSb$ and $[H]_\sim^\b$ yields the assertion.
\end{proofof}

Our next result investigates the   different notions of continuity for $k$ introduced in Definition 
\ref{k-conts-def}.

\begin{lemma}\label{simple-topol}
   Let $(T,\ca B)$ be a measure space and $k$ be a kernel on $T$ with RKHS $H$ and canonical feature map $\P:T\to H$. Then the 
   following statements are true:
   \begin{enumerate}
      \item The topology $\t_k$ is the smallest topology $\t$ on $T$ for which $k$ is $\t$-continuous. Moreover, we have
     \begin{displaymath}
   \t_k =   \t\bigl(\P:T\to (H,\snorm\cdot_H)\bigr)\, ,
  \end{displaymath}
    where $\t (\P:T\to (H,\snorm\cdot_H) )$ denotes the initial topology of $\P$ with respect to the norm-topology on $H$.
    \item The topology $\t(H)$ is the smallest topology $\t$ on $T$ for which $\P$ is continuous with respect to the 
    weak topology $w$ on $H$, that is 
    \begin{displaymath}
   \t(H) =  \t\bigl(\P:T\to (H, w)\bigr)\, .
  \end{displaymath}
    In particular, we have $\t(H)\subset \t_k$, and in general, the converse inclusion is not even true for $T=[0,1]$.
    \item If $H$ is separable and $k$ is bounded, then there exists a pseudo-metric on $T$ that generates $\t(H)$
	and $\t(H)$ is separable. Moreover, we have $\t(H)\subset \s(H)$.
    \item If  $\t(H)\subset \ca B$, then all $f\in H$ are $\ca B$-measurable.
  \end{enumerate}
\end{lemma}

\begin{proofof}{Lemma \ref{simple-topol}}
   \ada i Both assertions are shown in \cite[Lemma 4.29]{StCh08}.
   
   \ada {ii} Let $\iota:H\to H'$  be the Fr\'echet-Riesz isometric isomorphism. Then we have 
   $f= (\iota f) \circ \P$ for all $f\in H$ by the reproducing property. Let us first prove the inclusion 
   ``$\subset$''. To this end, we fix an $f\in H$ and an open $U\subset \R$. We define $O:= (\iota f)^{-1}(U)$. Then we have 
   $O\in w$ and thus
   \begin{displaymath}
      f^{-1}(U) = \bigl( (\iota f) \circ \P\bigr)^{-1}(U) =\P^{-1} \bigl( (\iota f)^{-1}(U) \bigr) =  \P^{-1}(O) 
      \in \t(\P:T\to (H, w))\, .
   \end{displaymath}   
   The inclusion ``$\subset$'' then follows from the fact that the set of considered pre-images $f^{-1}(U)$
   is a sub-base of $\t(H)$.
   To show   the converse inclusion, we 
    fix an $O\in w$ for which there exist an $f\in H$ and an open $U\subset \R$
   with $O= (\iota f)^{-1}(U)$. Then we find 
   \begin{displaymath}
      \P^{-1}(O) = \P^{-1} \bigl( (\iota f)^{-1}(U) \bigr) = \bigl( (\iota f) \circ \P\bigr)^{-1}(U) = f^{-1}(U)\in \t(H)\, .
   \end{displaymath}
   Since the set of such pre-images $\P^{-1}(O)$
   is a sub-base of  $\t(\P:T\to (H, w))$ we obtained the desired inclusion.
   
   Finally, $\t(H)\subset \t_k$ directly follows from combining part \emph{i)} and \emph{ii)} with the fact that the 
   norm topology on $H$ is finer than the weak topology. 
   To show that the converse inclusion does not hold for $T=[0,1]$, we denote the usual topology on this
   $T$ by $\t$.
   Then \cite{Lehto50a} showed that there exists 
    a bounded 
    separately $\t$-continuous kernel $k$ on $T$ that is not $\t$-continuous.
    This gives  $\t(H)\subset \t$ by 
    \cite[Lemma 4.28]{StCh08} and $\t_k\not\subset \t$, and thus $\t_k\not\subset \t(H)$.
   
  \ada {iii} 
  Since $H'$ is separable, we know that
  for every bounded subset $A'\subset H'$ the relative topology $w^*_{|A'}$  on $A'$, 
	where  $w^*$ denotes the weak* topology
  on $H'$,
  is induced be a metric, see e.g.~\cite[Corollary 2.6.20]{Megginson98}. Moreover, we have $\iota^{-1}(w^*) = w$, where 
  $w$ is the weak topology on $H$. For all bounded $A\subset H$, 
	the relative topology $w_{|A}$ on $A$ is thus induced by a metric. Now $k$ is bounded by assumption, 
	and hence $A:=\P(T)$ is bounded, see e.g.~\cite[p.~124]{StCh08}.
  Consequently, there exists a metric $d$ on $A$ that generates $w_{|A}$.
  Let us  consider the map $\tilde \P:T\to A$, defined by $\tilde \P(t) := \P(t)$ for all $t\in T$.
  By the already proven part \emph{ii)} and   the universal property of 
  the initial topology $\t(\id:A\to (H,w))=w_{|A}$ we then find
      \begin{displaymath}
   \t(H) =  \t\bigl(\P:T\to (H, w)\bigr) = \t\bigl(\tilde \P:T\to (A, w_{|A})\bigr)\, .
  \end{displaymath}
 From this we easily derive that $(t,t')\mapsto d(\P(t), \P(t'))$ is the desired pseudo-metric.
	To see that $\t(H)$ is separable, we recall that closed unit ball $B_{H'}$ 
	of $H'$ is $w^*$-compact by Alaoglu's theorem. Consequently, $(B_{H'}, w^*_{|B_{H'}})$ is a compact
	metric space, and thus separable. Arguing as above, 
	and using that $w_{|B_H} = \iota^{-1}(w^*_{|B_{H'}})$ is metrizable,
	we see that $w_{|A}$ is separable for $A:=\P(T)$,
	and hence so is $\t(H)$.
	
  Finally, since $\t(H)$ is the initial topology of $H$, the collection of 
  sets $f^{-1}(O)$, where $f\in H$ and $O\subset \R$ open,
    form a sub-base of $\t(H)$, and since open $O\subset \R$ are Borel measurable, we also have  $f^{-1}(O)\in \s(H)$
    for all such $f$ and $O$. Consequently, finite intersections taken from this sub-base are contained in $\s(H)$, too,
    and the collection of these intersections form a base of $\t(H)$. Now every $\t(H)$-open set is the 
    union of such intersections. However, we have just seen that $\t(H)$ is separable and generated by a pseudo-metric,
    which by a standard argument shows that $\t(H)$ is second countable. Consequently, $\t(H)$ is 
    Lindel\"of, see \cite[p.~49]{Kelley55}, that is each open cover has a countable sub-cover.
    Consequently, each $\t(H)$-open set is a countable 
    union of the above intersections, and thus contained in $\s(H)$.
 
 \ada {iv} From $\t(H)\subset \ca B$ we conclude that $\s(H) \subset \s(\t(H)) \subset \ca B$, which shows the assertion.
\end{proofof}

\begin{proofof}{Lemma \ref{Ik-inject}}
Let us pick an $f\in H$ with $f\neq 0$.
Then $\{f\neq 0\}$ is $\t(H)$-open and non-empty, and  thus we have $\nu(\{f\neq 0\}) >0$, that is
$I_kf = [f]_\sim \neq 0$. Now, $k=k_T^1$ follows from \cite[Theorem 3.1]{StSc12a}. 
\end{proofof}

The following results investigates the behavior of series of non-negative, continuous functions.

\begin{lemma}\label{lower-bound-on-diag}
Let $(T,\tau)$ be a topological space, $I\subset \N$, and $(g_i)_{i\in I}$ be a family 
of continuous functions $g_i:T\to \R$. Then, for all $t\in T$, the following statements hold:
\begin{enumerate}
   \item If $\sum_{i\in I} g_i^2(t) = \infty$, then, for all $M>0$, there exists an open $O\subset T$
    with $t\in O$ and 
    \begin{displaymath}
       \sum_{i\in I} g_i^2(s) > M\,, \qquad \qquad s\in O\, .
    \end{displaymath}
   \item If $\sum_{i\in I} g_i^2(t) < \infty$, then, for all $\e>0$, there exists an open $O\subset T$
    with $t\in O$ and 
    \begin{displaymath}
       \sum_{i\in I} g_i^2(s) > \sum_{i\in I} g_i^2(t) - \e \,, \qquad \qquad s\in O\, .
    \end{displaymath}
\end{enumerate}
\end{lemma}

\begin{proofof}{Lemma \ref{lower-bound-on-diag}}
   \ada i By assumption, there exists a finite $J\subset I$ such that
   \begin{displaymath}
      \sum_{i\in J} g_i^2(t) > 2M\, .
   \end{displaymath}
   Since the $g_i^2$ are continuous, there then exist, for all $i\in J$, an open $O_i\subset T$ with $t\in O_i$
   and $|g_i^2(s) - g_i^2(t)| < M/|J|$ for all $s\in O_i$. For the open set $O:= \bigcap_{i\in J}O_i$ and $s\in O$
   we then obtain $t\in O$ and 
   \begin{displaymath}
   \biggl|  \sum_{i\in J} g_i^2(s)  -  \sum_{i\in J} g_i^2(t)   \biggr| \leq \sum_{i\in J} |g_i^2(s) - g_i^2(t)| < M\, .
   \end{displaymath}
  This yields 
  \begin{displaymath}
      \sum_{i\in I} g_i^2(s) \geq  \sum_{i\in J} g_i^2(s) >  \sum_{i\in J} g_i^2(t) - M > M\, .
  \end{displaymath}
  \ada {ii} Let us fix an $\e>0$. Then 
  there exists a finite $J\subset I$ such that
  \begin{displaymath}
     \sum_{i\in J} g_i^2(t) > \sum_{i\in I} g_i^2(t) -\e\, .
  \end{displaymath}
  This time we pick open $O_i\subset T$ with $t\in O_i$
   and $|g_i^2(s) - g_i^2(t)| < \e/|J|$ for all $s\in O_i$. Repeating the calculations above, we obtain the assertion for $2\e$.
\end{proofof}

\begin{proofof}{Theorem \ref{bounded-powers-of-cont-kernels}}
	By our assumption and \eqref{interpol-repres} we have $\eHnb \hookrightarrow \Lx \infty {\nu}$, and thus
	\cite[Theorem 5.3]{StSc12a} shows that there exist an $N\in {\ca B}$ and a constant $\k\in [0,\infty)$
	such that $\nu (N) =0$ and
	\begin{equation}\label{bounded-powers-of-cont-kernels-powers-h1}
	 \sum_{i\in I}  \mu_i^{\b} e_i^2(t) \leq \kappa^2 \, , \qquad \qquad t\in T\setminus  N.
	\end{equation}
	Moreover, by the definition of $\t(H)$ we know that all $e_i$ are $\t(H)$-continuous.

  Let us first show that  \eqref{k-series-diag} holds for $S:=T$. To this end, we assume the converse, that 
  is, there exists a $t\in T$ with 
  \begin{displaymath}
     \sum_{i\in I} \mu_i^\b e_i^2(t) = \infty\, .
  \end{displaymath}
	By Lemma \ref{lower-bound-on-diag} there then exists an  $O\in \t(H)$ with $t\in O$ and 
	\begin{equation}\label{bounded-powers-of-cont-kernels-powers-h2}
	 \sum_{i\in I} \mu_i^\b e_i^2(s) > \k^2\, , \qquad \qquad s\in O\, .
	\end{equation}
    Since $\nu$ is assumed to be $k$-positive, we conclude that $\nu(O)>0$, and hence there 
	exists a $t_0\in O\setminus  N$. For this $t_0$ we have both \eqref{bounded-powers-of-cont-kernels-powers-h1}
	and \eqref{bounded-powers-of-cont-kernels-powers-h2}, and thus we have found a contradiction.

	To show that $\kTb$ is bounded, we again assume the converse. Then there exists a $t\in T$ such that 
	\begin{displaymath}
	 \sum_{i\in I} \mu_i^\b e_i^2(t) > \k^2 + 1\, ,
	\end{displaymath}
	so that by
	using $\e:=1$ in part \emph{ii)} of
	Lemma \ref{lower-bound-on-diag} we again find an $O\in \t(H)$ with $t\in O$ and 
	\eqref{bounded-powers-of-cont-kernels-powers-h2}. Repeating  the arguments above we then obtain a contradiction.

	Let us now show that $\t(\HTb)= \t(H)$. To this end, we first fix an $f\in \HTb$. Since $(\mu_i^{\b/2}e_i)_{i\in I}$
	is an ONB of $\HTb$, see \cite[Lemma 2.6 and Proposition 4.2]{StSc12a}, we then have
	\begin{displaymath}
	 f = \sum_{i\in I} \bigl\langle f, \mu_i^{\b/2}e_i\bigr\rangle_{\HTb} \mu_i^{\b/2}e_i \, ,
	\end{displaymath}
	where the convergence is unconditionally in $\HTb$. Since $\kTb$ is bounded, convergence in $\HTb$
	implies uniform convergence, see e.g.~\cite[Lemma 4.23]{StCh08}, and thus the above series 
	also converges unconditionally  with respect to  $\inorm\cdot$. Consequently, $f$ is 
	a $\inorm\cdot$-limit of a sequence of $\t(H)$-continuous functions, and thus itself $\t(H)$-continuous.
	From this we easily conclude that $\t(\HTb)\subset \t(H)$. To show the converse inclusion $\t(H)\subset \t(\HTb)$
	let us recall that the embedding $I_k:H\to \Ltn$ is injective and $H=H_T^1$
	by Lemma \ref{Ik-inject}.
	Now the inclusion $\t(H)\subset \t(\HTb)$ trivially follows from 
	the inclusion $H_T^1\subset \HTb$ established in  \cite[Lemma 4.3]{StSc12a}. 
\end{proofof}

The next result characterizes under which conditions we have $\kxy S1 \ll \kxy S\b$.

\begin{lemma}\label{nuclear-dom}
    Let Assumption K be satisfied.
Then, for all $\b\in (0,1)$ and all measurable $S\subset T$ 
satisfying
$\nu(T\setminus S)=0$ and \eqref{k-series-diag}, the restriction operator 
$
 \cdot_{|S}:\Hxy T1\to \HSb
$
is compact,
and
the following 
statements are equivalent:
\begin{enumerate}
   \item The operator $\cdot_{|S}:\Hxy T1\to \HSb$ is Hilbert-Schmidt.
   \item We have $\sum_{i\in I}\mu_i^{1-\b}<\infty$.
	\item We have $\kxy S1 \ll \kxy S\b$.
\end{enumerate}
\end{lemma}

For the proofs of Lemma \ref{nuclear-dom} and Lemma \ref{entropy-versus-ew}, we need to recall 
some basics on singular numbers. To begin with, 
let us recall that for an arbitrary compact operator $S:H_1\to H_2$ acting between two Hilbert spaces $H_1$ and $H_2$,
the $i$-th singular number, see e.g.~\cite[p.~242]{BiSo87}  is defined by 
\begin{equation}\label{def-snumber}
  s_i(S) := \sqrt{\mu_i(S^*S)}\, ,
\end{equation}
where $\mu_i(S^*S)$ denotes the $i$-th non-zero eigenvalue of the compact, positive and self-adjoint operator $S^*S$.
As usual, these eigenvalues are assumed to be ordered with duplicates according to their
geometric multiplicities. In addition, we extend the sequence of eigenvalues by zero, if we only have finitely
many non-zero eigenvalues. Now, for a compact, self-adjoint and positive $T:H\to H$,
this definition gives
\begin{equation}\label{snumber-forT}
  s_i(T) = \sqrt{\mu_i(T^*T)} = \sqrt{\mu_i(T^2)} = \mu_i(T)\, , \qquad \qquad i\geq 1,
\end{equation}
where the last equality follows from the classical spectral theorem for such $T$, see 
e.g.~\cite[Theorem V.2.10 on page 260]{Kato76} or \cite[Satz VI.3.2]{Werner95}.
For compact $S:H_1\to H_2$ and $T:= S^*S$ we thus find 
\begin{equation}\label{snumber-SvsT}
  s_i^2(S) = \mu_i(S^*S)=\mu_i(T)=s_i(T)
\end{equation}
for all $i\geq 1$.
Consequently, we have $(s_i(S))\in \ell_2$ if and only if $(s_i(T))\in \ell_1$.
Moreover, $T$ is nuclear, if and only if  $(s_i(T))\in \ell_1$, see e.g.~\cite[Satz VI.5.5]{Werner95}
or \cite[p.~245ff]{BiSo87}, while $S$ is Hilbert-Schmidt if and only if $(s_i(S))\in \ell_2$, see 
e.g.~\cite[p.~250]{BiSo87}, \cite[Prop.~2.11.17]{Pietsch87}, or \cite[p.~246]{Werner95}.
\newline

\begin{proofof}{Lemma \ref{nuclear-dom}}
	We first observe that, for $i\in I$, we have 
	\begin{equation}\label{restrict-op}
	 \cdot_{|S}(\mu_i^{1/2}e_i) = \mu_i^{1/2}e_{i|S} = \mu_i^{(1-\b)/2} \mu_i^{\b/2} e_{i|S} \, .
	\end{equation}
	Since $(\mu_i^{1/2}e_i)_{i\in I}$ and $(\mu_i^{\b/2} e_{i|S} )_{i\in I}$ are ONBs of $H_T^1$ and $H_S^\b$, respectively, we 
	obtain the following commutative diagram

	\quadiasn {H_T^1}{H_S^\b}{\ell_2}{\ell_2} {\cdot_{|S}}{\Psi_1}{\Psi_\b}{D}

	\noindent
	where $\Psi_i$ denote the isometric isomorphisms that map each Hilbert space element 
	to its sequence of Fourier coefficients with respect to the ONBs above, and 
	$D$ is the diagonal operator with respect to the sequence $(\mu_i^{(1-\b)/2} )_{i\in I}$. Since the latter sequence
	converges to zero, $D$ is compact, and thus so is the restriction operator.

	\aeqb i {ii} We first observe that \eqref{restrict-op} yields 
   \begin{displaymath}
      \mnorm{\cdot_{|S}(\sqrt{\mu_i}e_i)}_{\HSb}^2 = \mu_i^{1-\b} \, , \qquad \qquad i\in I.
   \end{displaymath}
  Since   $(\sqrt{\mu_i} e_i)_{i\in I}$ is an ONB of $\Hxy T1$, 
  the equivalence \emph{i)} $\Leftrightarrow$ \emph{ii)} immediately follows from the fact, 
	see e.g.~\cite[p.~243f]{Wojtaszczyk91}, that 
	$\cdot_{|S}:\Hxy T1\to \HSb$  is Hilbert-Schmidt, if and only if 
	\begin{displaymath}
	 \sum_{i\in I}  \snorm{\cdot_{|S}(\sqrt{\mu_i} e_i)}_{\bHSb}^2< \infty\, .
	\end{displaymath}

	\aeqb i {iii} The restriction operator admits the following natural factorization

	\tridia {\Hxy T1}{\HSb}{\Hxy S1} {\cdot_{|S}}{\cdot_{|S}}{I_{\kxy S1, \kxy S\b}}

	\noindent
	where there restriction operator $\cdot_{|S}:\Hxy T1\to \Hxy S1$ is an isometric isomorphism.
	Consequently, $\cdot_{|S}:\Hxy T1\to \HSb$ is Hilbert-Schmidt, if and only if
	$I_{\kxy S1, \kxy S\b}$ is Hilbert-Schmidt. In view of the desired equivalence, it suffices to show that 
	$I_{\kxy S1, \kxy S\b}$ is Hilbert-Schmidt, if and only if 
	$I_{\kxy S1, \kxy S\b}\circ S_{\kxy S1, \kxy S\b}$ is nuclear.
	However, since  $S_{\kxy S1, \kxy S\b} = I_{\kxy S1, \kxy S\b}^*$, this equivalence is a  simple consequence
	of the remarks on singular numbers made in front of this proof, if we consider
	the compact operator  $I_{\kxy S1, \kxy S\b}:\Hxy S1 \to \Hxy S\b$ for $S^*$.
\end{proofof}

\begin{proofof}{Lemma \ref{entropy-versus-ew}}
Let us denote the $i$-th approximation number of a bounded linear operator $T:E\to F$
between Banach spaces $E$ and $F$ by $a_i(T)$, that is 
\begin{displaymath}
   a_i(T) := \inf\bigl\{ \snorm{T-A}\, \bigl|\,  A:E\to F \mbox{ bounded linear with } \rank A<i  \bigr\}\, .
\end{displaymath}
Moreover,
we write $s_i(I_k)$ for the $i$-th singular number of $I_k$, see \eqref{def-snumber}.
Since $I_k$ is compact, we actually have $a_i(I_k) = s_i(I_k)$ for all $i\geq 1$,
see \cite[Theorem 7 on p.~240]{Wojtaszczyk91}, and using 
\eqref{snumber-forT} and \eqref{snumber-SvsT} we thus find 
$$
\mu_i = \mu_i(T_k) = s_i(T_k) =  s_i^2(I_k) = a_i^2(I_k)
$$
for all $i\in I$. Moreover, if $|I|<\infty$, then we clearly have $a_i(I_k) = 0$ for all $i>|I|$
by the spectral representation of $T_k$.
From Carl's inequality, see \cite[Theorem 3.1.2]{CaSt90}, we then obtain \eqref{carl-ineq}.
Moreover, \eqref{inverse-car-ineq} follows from the relation 
$$
a_i(R:H_1\to H_2) \leq 2\e_i(R:H_1\to H_2)
$$
that holds for all compact linear operators $R$ between Hilbert spaces $H_1$ and $H_2$, see \cite[p.~120]{CaSt90}.
The second to last equivalence is a direct consequence of \eqref{inverse-car-ineq} and another application
of Carl's inequality, while the last equivalence follows from \eqref{inverse-car-ineq} and 
Carl's inequality with the help of a little trick employed in the proof of \cite[Proposition 2]{Steinwart00a}.
\end{proofof}

\subsection{Proofs Related to Generic KL-Expansions}

\begin{proofof}{Lemma \ref{cov-prop}}
   Since $X$ is ($\ca A\otimes \ca B$)-measurable, the map $(\om,s,t)\mapsto X_s(\om)X_t(\om)$ is 
   $\ca A\otimes \ca B\otimes \ca B$-measurable. From this we easily conclude that $k$ is measurable.
   Moreover, a simple application of Tonelli's theorem  shows 
   \begin{displaymath}
      \int_T k(t,t) \,d\nu(t) = \int_T \E_P X_t^2 \,d\nu(t) = \int_{\Om\times T} X^2\, dP\otimes\nu\, ,
   \end{displaymath}
    from which the conclude the equivalence.
  The remaining assertions  now follow from \cite[Lemma 2.3]{StSc12a}.
\end{proofof}

\begin{proofof}{Lemma \ref{kernel-instead-of-EV}}
   By \cite[Lemma 2.6]{StSc12a} we know that $(\sqrt{\mu_i}e_i)_{i\in I}$ is an ONB of the RKHS $H$ of $k$.
    Furthermore, this lemma shows that $H$ is compactly embedded into $\Ltn$ and 
    \cite[Theorem 2.11]{StSc12a} gives the spectral representation \eqref{Tk-spectral} of $T_k$. 
    Finally, \eqref{Sk-EW}-\eqref{Sks-ran} follow from \cite[Lemma 2.12]{StSc12a} and $X\in \sLx 2{P\otimes \nu}$
    follows from Lemma \ref{cov-prop} in combination with $\sum_{i\in I}\mu_i<\infty$ and part (iii) of \cite[Theorem 2.11]{StSc12a}.
\end{proofof}

\begin{proofof}{Lemma \ref{prop-Zi}}
 For $i\in I$ and $\om\in \Om$, we define 
\begin{displaymath}
   Y_i(\om) := \int_T \bigr| X_t(\om ) e_i(t)\bigr| \, d\nu(t)\, ,
\end{displaymath}
where we note that the measurability of $(\om,t)\mapsto X_t(\om)e_i(t)$ together with Tonelli's theorems 
shows that $Y_i:\Om\to [0,\infty]$ is measurable. Moreover, since we have 
$e_i\in \sLt \nu$  with $\snorm{e_i}_{\sLtn} = 1$ as well as
 $X(\om) \in \sLt \nu$
for $P$-almost all $\om \in \Om$,  
 Cauchy-Schwarz inequality implies 
\begin{align}
 \E_P Y_i^2 \nonumber
& =
\int_\Om \biggl( \int_T \bigr| X_t(\om ) e_i(t)\bigr| \, d\nu(t)  \biggr)^2 dP(\om) \\  \label{Zi-in-Ltn-arg}
& \leq
\int_\Om \biggl( \int_T X_t^2 (\om ) \, d\nu(t)\biggr) \biggl(\int_T e_i^2(t) \, d\nu(t)\biggr) dP(\om) \\ \nonumber
& = 
\int_{\Om\times T} X^2 \, dP\otimes \nu \\ \nonumber
& < \infty\, .
\end{align}
Since $|Z_i| \leq |Y_i|$, we then obtain $Z_i\in \sLx 2 P$. Furthermore, we have 
  $Xe_i \in \sLo{P\otimes \nu}$ since another application of 
 the Cauchy-Schwarz inequality gives
\begin{align}\label{Xei-norm-}
 \int_{\Om\times T} \bigl| Xe_i \bigr|\, dP\otimes \nu 
\leq 
\biggl(\int_{\Om\times T} X^2 \, dP\otimes \nu \biggr)^{1/2} 
\biggl(\int_{\Om\times T}  e_i^2 dP\otimes \nu \biggr)^{1/2} 
= 
\snorm X_{\sLt{P\otimes\nu}} 
 < \infty\, .
\end{align}
Consequently, we can apply  Fubini's theorem, which yields
\begin{align*}
 \E_P Z_i 
=
\int_\Om   \int_T X_t(\om ) e_i(t) \, d\nu(t)\,  dP(\om)   
=
 \int_T  \int_\Om   X_t(\om ) e_i(t) \, dP(\om) \, d\nu(t)   
 = 0\, ,
\end{align*}
where in the last step we used $\E_P X_t = 0$. 
To show \eqref{Zi-cov}, we first observe that 
\begin{align} \nonumber
& \int_{\Om\times T\times T}  \bigl| X_s(\om ) e_i(s)    X_t(\om ) e_j(t)\bigr| \, dP\otimes\nu\otimes\nu(w,s,t) \\\nonumber
& = 
\int_\Om \int_T \int_T  \bigl| X_s(\om ) e_i(s)\bigr| \cdot   \bigl|  X_t(\om ) e_j(t)\bigr|  \, d\nu(s)  \, d\nu(t)  \, dP(\om)\\ \nonumber
 & = 
\int_\Om \biggl( \int_T \bigl|X_s(\om ) e_i(s)  \bigr| \, d\nu(s) \biggr)  \biggl( \int_T \bigl|X_t(\om ) e_i(t)  \bigr| \, d\nu(t) \biggr) dP(\om) \\ \label{prop-Zi-h1}
& =
\E_P Y_i^2 < \infty\, .
%
\end{align}
where in the last inequality we used the arguments from \eqref{Zi-in-Ltn-arg}.
Using Fubini's theorem, we then obtain
\begin{align}\nonumber
 \E_P Z_i Z_j 
& = 
\int_\Om \biggl( \int_T X_s(\om ) e_i(s) \, d\nu(s)  \biggr) \biggl( \int_T X_t(\om ) e_j(t) \, d\nu(t)  \biggr) dP(\om) \\ \nonumber
& = 
\int_\Om   \int_T\int_T  X_s(\om ) e_i(s)    X_t(\om ) e_j(t) \, d\nu(s) \, d\nu(t)\,    dP(\om) \\ \nonumber
& = 
 \int_T\int_T \E_P \bigl(  X_s     X_t\bigr)  e_i(s)  e_j(t)  \, d\nu(s) \, d\nu(t)\\ \label{prop-Zi-h2}
& =
\int_T S_k ([e_i]_\sim)(t)\,  e_j(t) \, d\nu(t)\\ \nonumber
& =
\int_T \mu_i  e_i(t) e_j(t) \, d\nu(t)\\ \nonumber
& = \mu_i \d_{i,j}\, ,
\end{align}
where in the second to last step we used \eqref{Sk-EW}. 

Let us now show \eqref{Zi-Xt-prod}.
To this end, note that 
  the already established $Y_j \in \sLx 2 P$ together with 
$X_t\in \sLx 2 P$ and Tonelli's theorem
implies 
\begin{align*}
   \int_{\Om\times T}  \bigl| X_t (\om)     X_s(\om ) e_j(s) \bigr|  \, dP\otimes\nu(\om,s) 
   =\int_\Om \bigl|X_t(\om) Y_j(\om)\bigr| \, dP(\om) < \infty
\end{align*}
for all $t\in T$.
Consequently, the map $(\om,s)\mapsto X_t (\om)  X_s(\om ) e_j(s)$ is 
$P\otimes \nu$-integrable for 
 each $t\in T$, and by
  Fubini's theorem we thus  obtain
\begin{align*}
 \E_P X_t  Z_j  
&= 
\int_\Om X_t (\om)    \int_T X_s(\om ) e_j(s) \, d\nu(s) \, dP(\om) \\
&= 
\int_T    e_j(s)  \int_\Om     X_t (\om) X_s(\om ) \, dP(\om)\, d\nu(s)  \\
& =
 \int_T   e_j(s)  k(s,t) \, d\nu(s)  \\
& = 
 \mu_j e_j(t)\, ,
\end{align*}
where in the last step we used the definition of $S_k$ and \eqref{Sk-EW}.

Moreover,    \eqref{LP-diff} immediately follows from
\begin{align*} 
 \bnorm{X_t - \sum_{j\in J} Z_j e_j(t)}_{\sLt P}^2
& =
\E_P X_t^2 - 2 \E_P X_t  \sum_{j\in J} Z_j e_j(t) + \E_P \sum_{i,j\in J} Z_i e_i(t) Z_j e_j (t)\\
& =
k(t,t) - 2  \sum_{j\in J}\E_P X_t  Z_j e_j(t) + \sum_{i,j\in J}  e_j  (t) e_i(t) \E_P Z_i  Z_j \\
& =
k(t,t) - 2  \sum_{j\in J} \mu_j e_j^2(t) + \sum_{j\in  J} \mu_j e_j^2(t)\, ,
\end{align*}
where in the  last step we used the already established \eqref{Zi-cov} and \eqref{Zi-Xt-prod}.

\aeqb i {ii} Follows directly from \eqref{LP-diff}.

Finally, to show \eqref{path-in-compl-kernel}, 
we fix a measurable $N\subset \Om$ with $P(N) = 0$ and
$X(\om) \in\sLtn$ for all $\om \in \Om\setminus N$.
Furthermore, we 
 fix an $f\in \sLtn$ with $[f]_\sim\in \ker T_k$.
For  $\om \in N$ we now write $Z(\om) := 0$
and   
\begin{displaymath}
Z(\om) := \int_T X_t(\om ) f(t) \, d\nu(t) 
\end{displaymath}
 otherwise. Then, repeating \eqref{Zi-in-Ltn-arg} and \eqref{Xei-norm-} with $e_i$ replaced by $f$ we 
obtain $Z\in \sLx 2P$ and $Xf\in \sLx 1 {P\otimes \nu}$. Moreover, repeating  
\eqref{prop-Zi-h1} and
\eqref{prop-Zi-h2} in the same way,
we obtain
\begin{displaymath}
 \E_P Z^2 = \int_T S_k ([f]_\sim)(t)  f(t) \, d\nu(t) = 0
\end{displaymath}
since $[f]_\sim\in \ker T_k = \ker S_k$ by \eqref{Sk-ker}.
This gives a measurable $N_f\subset \Om$ with $N\subset N_f$, $P(N_f)=0$, and 
%
%
$\langle [X(\om)]_\sim, [f]_\sim\rangle_{\Ltn} = Z(\om) = 0$ for all $\om \in \Om\setminus N_f$.
Now, since $\Ltn$ is separable by Assumption X, there exists a 
countable family $(f_n)_{n\geq 1}\subset \sLtn$ such that $([f_n]_\sim)_{n\geq 1}\subset \ker T_k$ is dense.
We define $N^* :=  \bigcup_{n\geq 1} N_{f_n}$. Clearly, $N^*$ is measurable with $N\subset N^*$, $P(N^*) = 0$,
and 
$\langle [X(\om)]_\sim, [f_n]_\sim\rangle_{\Ltn} =   0$ for all $\om \in \Om\setminus N^*$ and all $n\geq 1$.
Now let $f\in \sLtn$ with $[f]_\sim\in \ker T_k$. Then there exists a sub-sequence $(f_{n_m})_{m\geq 1}$ with 
$\lim_{m\to \infty} [f_{n_m}]_\sim = [f]_\sim$ in $\Ltn$ and for $\om\in\Om\setminus N^*$ we conclude that 
\begin{displaymath}
   \langle [X(\om)]_\sim, [f]_\sim\rangle_{\Ltn} 
= \lim_{m\to \infty} \langle [X(\om)]_\sim, [f_{n_m}]_\sim\rangle_{\Ltn} = 0\, .
\end{displaymath}
%
%
Thus we have found the first part of \eqref{path-in-compl-kernel}.
The second part of \eqref{path-in-compl-kernel}, namely,
\begin{displaymath}
 (\ker T_k)^\perp = \overline{\spann\{[e_i]_\sim:i\in I  \}}^{\Ltn}\, ,
\end{displaymath}
follows from combining \eqref{Sk-ker} with \eqref{Sks-ran} and \eqref{Sks-ran-or}.
\end{proofof}

\begin{proofof}{Proposition \ref{classical-KL}}
Recall that \cite[Theorem 3.1]{StSc12a} showed that both \emph{i)} and \emph{ii)}
are equivalent to 
\begin{equation}\label{kernel-mercer-eq}
k(t, t') = \sum_{i\in I}  \mu_i e_i(t) e_i(t')\, .
\end{equation}
for all $t,t'\in T$. In view of \eqref{LP-diff} it thus 
suffices to show that $\emph{iii)} \Rightarrow \emph{i)}$.
To show this implication we assume that \eqref{KL-L2P-conv} holds for all $t\in T$, but 
$(\sqrt{\mu_i}e_i)_{i\in I}$ is not an ONB of $H$. Let $(\tilde e_j)_{j\in J}$ be an ONS of $H$
such that the union of $(\sqrt{\mu_i}e_i)_{i\in I}$ and $(\tilde e_j)_{j\in J}$ is an ONB of $H$.
By assumption we know that $J\neq \emptyset$, so we can fix a $j_0\in J$. Since
$\snorm{\tilde e_{j_0}}_H= 1$, there further exists a $t\in T$ with $\tilde e_{j_0}(t) \neq 0$.
Now, it is well-known that the kernel $k$ can be expressed in terms of our ONB, 
see e.g.~\cite[Theorem 4.20]{StCh08}, and hence we obtain 
\begin{align*}
 k(t,t) =  \sum_{i\in I}  \mu_i e_i^2(t)  + \sum_{j\in J} \tilde e_j^2(t) 
\geq \sum_{i\in I}  \mu_i e_i^2(t) + \tilde e_{j_0}^2(t) 
> \sum_{i\in I}  \mu_i e_i^2(t)  
= k(t,t)\, ,
\end{align*}
where the last equality follows from the equivalence of \eqref{KL-L2P-conv} and \eqref{kernel-repres-diag}.
In other words, we have found a contradiction, and hence $\emph{iii)} \Rightarrow \emph{i)}$ is true.

Let us finally consider the case in which $H$ is separable. By 
\cite[Corollary 3.2 and Theorem 3.3]{StSc12a} we then see 
that there exists  a measurable $N\subset T$ with $\nu(N)=0$ such that 
\begin{displaymath}
 k(t, t') = k_T^1(t,t')\, , \qquad \qquad t,t'\in T\, .
\end{displaymath}
Consequently, \eqref{kernel-repres-diag} holds for all $t\in T\setminus N$,
and we obtain the assertion by \eqref{LP-diff}.
\end{proofof}

\begin{proofof}{Theorem \ref{pointwise-KL-new}}
Equation \eqref{path-in-compl-kernel} shows that there exists a measurable $N_1\subset \Om$ with $P(N_1) = 0$
such that for all $\om \in \Om\setminus N_1$ 
the path 
$[X(\om)]_\sim$ is contained in the space spanned by the ONS $([e_i]_\sim)_{i\in I}$.
Moreover, by the definition of $Z_i$ there exists another measurable $N_2\subset \Om$ with $P(N_2)= 0$
and 
\begin{equation}\label{pointwise-KL-h1}
 Z_i(\om) = \langle [X(\om)]_\sim, [e_i]_\sim\rangle_{\Ltn}
\end{equation}
for $\om \in \Om\setminus N_2$.
Let us define $N:= N_1 \cup N_2$. For $\om \in \Om\setminus N$ we then obtain \eqref{KL-path-conv-Lt}.

To show \eqref{Ln-diff}, we again pick an $\om \in \Om\setminus N$. Using Parseval's identity and
\eqref{pointwise-KL-h1}, we obtain 
\begin{displaymath}
 \bnorm{[X(\om)]_\sim - \sum_{j\in J} Z_j(\om) [e_j]_\sim}_{\Ltn}^2 = \sum_{i\in I\setminus J} Z_i^2(\om)
\end{displaymath}
Furthermore, Lemma \ref{prop-Zi} 
implies
\begin{equation}\label{pointwise-KL-h2}
 \E_P \sum_{i\in I\setminus J} Z_i^2 = \sum_{i\in I\setminus J} \E_P  Z_i^2 = \sum_{i\in I\setminus J} \mu_i\, .
\end{equation}
Combining both equations then yields \eqref{Ln-diff} and the last assertion is a trivial consequence of \eqref{Ln-diff}.
%
\end{proofof}

\begin{proofof}{Corollary \ref{Zi-in-CM}}
	Our first goal is to show that $[Z_i]_\sim \in \Lx 2 X$ for all $i\in I$. To this end, 
	recall from e.g.~\cite[p.~65]{BeTA04} and \cite[Chapter 8.4]{Janson97}
	that the  Lo{\`e}ve isometric isomorphism $\Psi:\Lx 2 X\to H$ is the unique continuous extension 
	of the well-defined linear map $\Psi_0 :\spann\{[X_t]_\sim: t\in T\} \to  \spann \{ k(t,\cdot): t\in T\}$
	described by 
	  \begin{displaymath}
     \Psi_0\Bigl(\sum_{i=1}^n a_i[X_{t_i}]_\sim  \Bigr) := \sum_{i=1}^n a_ik(t_i, \cdot)\, .
  \end{displaymath}
  Now let $(\tilde e_j)_{j\in J}$ be an ONS in $H$ such that 
  $(\sqrt{\mu_i}e_i)_{i\in I} \cup (\tilde e_j)_{j\in J}$ is an ONB of $H$. For an arbitrary $t\in T$ 
  and all $i\in I$ and $j\in J$, we
  then find $\langle k(t,\cdot) , \sqrt{\mu_i} e_i\rangle_H = \sqrt{\mu_i} e_i(t)$
  and  $\langle k(t,\cdot) , \tilde e_j\rangle_H = \tilde e_j(t)$ and thus we obtain
  \begin{displaymath}
     k(t,\cdot ) = \sum_{i\in I} \mu_i e_i(t) e_i + \sum_{j\in J} \tilde e_j(t) \tilde e_j\, ,
  \end{displaymath}
  where the series converge unconditionally in $H$. Applying $\Psi^{-1}$ on both sides yields
  \begin{displaymath}
     [X_t]_\sim
     = \Psi^{-1} \bigl( k(t,\cdot) \bigr)
     = \sum_{i\in I} \mu_i e_i(t) \Psi^{-1}(e_i) + \sum_{j\in J} \tilde e_j(t) \Psi^{-1}(\tilde e_j)\, ,
  \end{displaymath}
  where the series converge unconditionally in $\Lx 2 P$. Let us fix $\xi_i, \tilde \xi_j\in \sLx 2 P$ with 
  $[\xi_i]_\sim =  \mu_i \Psi^{-1}(e_i)$
  and $[\tilde \xi_j]_\sim =   \Psi^{-1}(\tilde e_j)$.
  Then our constructions ensures
  \begin{equation}
     [X_t]_\sim =  \sum_{i\in I} [\xi_i]_\sim e_i(t)  + \sum_{j\in J} [\tilde\xi_j]_\sim \tilde e_j(t) \, ,
  \end{equation}
  where, for all $t\in T$, the series converge unconditionally in $\Lx 2 P$.
	For some fixed  finite sets $I_0\subset I$ and $J_0\subset J$, we further have 
	\begin{align*}
	 \int_\Om \bnorm{[X(\om)]_\sim - \sum_{i\in I_0} \xi_i(\om)  [e_i]_\sim }_\Ltn^2  \dxy P\om
	& =\int_\Om \int_T  \Bigl|  X_t (\om) - \sum_{i\in I_0} \xi_i(\om) e_i(t) \Bigr|^2  \dxy \nu t\dxy P\om \\
	& = \int_T \bnorm{[X_t]_\sim - \sum_{i\in I_0} \mu_i e_i(t) \Psi^{-1}(e_i)}_{\Lx 2 P}^2   \dxy \nu t\\
	& = \int_T \bnorm{k(t,\cdot) - \sum_{i\in I_0} \mu_i e_i(t) e_i }_{H}^2   \dxy \nu t \\
	& = \int_T \Bigl( \sum_{i\in I\setminus I_0} \mu_i e_i^2(t) + \sum_{j\in J} \tilde e_j^2(t) \Bigr) \dxy \nu t\\
	& = \sum_{i\in I\setminus I_0}\mu_i \, \mnorm{[e_i]_\sim}_\Ltn^2 + \sum_{j\in J} \,\mnorm{[\tilde e_j]_\sim}_\Ltn^2 \\
	& = \sum_{i\in I\setminus I_0} \mu_i\, ,
\end{align*}
		where in the last step we used  Theorem \ref{spectral}, which  implies
  \begin{displaymath}
     \tilde e_j \in \overline{\spann\{\sqrt{\mu_i}e_i:i\in I  \}}^\perp
     =
     (\overline{\ran S_k})^\perp
     = \ker S_k^*
     = \ker I_k\, .
  \end{displaymath}
	Consequently, there exists a measurable $N\subset \Om$ with $P(N)=0$ such that for all $\om \in \Om\setminus N$
	we have 
	\begin{displaymath}
	 [X(\om)]_\sim = \sum_{i\in I} \xi_i(\om)  [e_i]_\sim \, ,
	\end{displaymath}
	where the series converges in $\Ltn$. By Theorem \ref{pointwise-KL-new}
	we may assume without loss of generality that \eqref{KL-path-conv-Lt}
	also holds for $\om\in \Om\setminus N$.
Since $([e_i]_\sim)_{i\in I}$ is an ONS, we then see that 
	\begin{displaymath}
	 \xi_i(\om) = \langle [X(\om)]_\sim , [e_i]_\sim\rangle_{\Lx 2 P} = Z_i(\om)
	\end{displaymath}
	for such $\om$, and thus we finally obtain 
	$[Z_i]_\sim = [\xi_i]_\sim \in \Lx 2 X$.

	Now, \eqref{Zi-cov} shows that $(\mu_i^{-1/2}[Z_i]_\sim)_{i\in I}$ is an ONS of $\Lx 2 X$, and
	\eqref{LP-diff} together with Proposition \ref{classical-KL} shows that it is an ONB, if and only if 
	$(\sqrt{\mu_i}e_i)_{i\in I}$ is an ONB of $H$.
\end{proofof}

\begin{proofof}{Lemma \ref{gauss-lemma}}
  By Lemma \ref{prop-Zi} we know that the random variables $(Z_i)_{i\in I}$ are mutually uncorrelated and centered
  with $\var Z_i = \mu_i$ for all $i\in I$. 
	Moreover, by Corollary \ref{Zi-in-CM} we know $\sum_{i\in I_0}^n a_iZ_i\in \Lx 2 X$ for all
	finite $I_0\subset I$ and $a_i \in \R$. Since $\Lx 2 X$ consists of normally
	distributed random variables, which can be easily checked by L\'evy's continuity theorem, 
	we conclude that $(Z_i)_{i\in I}$ are jointly normal.
	Consequently, they are independent, and $Z_i \sim \ca N(0,\mu_i)$ becomes obvious.
\end{proofof}

\begin{proofof}{Theorem \ref{def-X-by-Z-new}}
   Let us first show that the series defining each $X_t$ do converge. To this end, we fix a finite $J\subset I$.
   Then an easy calculation shows
   \begin{align}\nonumber
      \int_\Om  \biggl( \sum_{j\in J} Z_j(\om) e_j(t) \biggr)^2 \dP \om
      &=
      \int_\Om \sum_{i,j\in J} Z_i(\om) Z_j(\om) e_i(t) e_j(t) \dP\om \\ \nonumber
      &= \sum_{i,j\in J} e_i(t) e_j(t) \E_P Z_i Z_j\\ \label{def-X-by-Z-h2}
      &=\sum_{j\in J} \mu_j e_j^2(t)\, .
   \end{align}
  By \eqref{diag-sum} we thus see that 
  the sequence of partial sums on the right-hand side of \eqref{def-X-by-Z-h1} is a Cauchy sequence in
  $\sLx 2 P$. Consequently, it converges, and by repeating the argument above we see that the series also 
  converges unconditionally.

  Let us now construct the $(\ca A \otimes \ca B)$-measurable version 
  $(Y)_{t\in T} \subset \sLx 2 P$ 
  of $(X_t)_{t\in T}$. Clearly, if $I$ is finite, there is nothing to prove, and hence we may assume without loss
  of generality that $I=\N$. Our first step in the construction of $Y$ is to show that the map 
  \begin{align*}
      \tilde X:T& \to \Lt P\\
      t&\mapsto [X_t]_\sim
  \end{align*}
  is $\Lt P$-measurable. To this end, we write $\xi_i := \mu_i^{-1/2} [Z_i]_\sim$. Clearly, $(\xi_i)_{i\in I}$ is an ONS 
  in $\Lt P$ and by  \eqref{diag-sum} and \eqref{def-X-by-Z-h1} we conclude that 
  \begin{displaymath}
     [X_t]_\sim \in \overline{\spann\{\xi_i: i\in I   \}}^{\Lt P}\, , \qquad \qquad t\in T.
  \end{displaymath}
  Consequently, the image of $\tilde X$ is contained in a separable subspace of $\Lt P$. Moreover, 
  the already established $\sLx 2 P$-convergence in \eqref{def-X-by-Z-h1} guarantees that, 
  for $f\in \Lt P$ and $t\in T$, we have 
  \begin{displaymath}
     \langle h, \tilde X\rangle_{\Lt P}
     =
     \Bigl \langle h , \sum_{i\in I} [Z_i]_\sim e_i(t)  \Bigr\rangle_{\Lt P}
     =
     \sum_{i\in I} e_i(t) \bigl\langle h , [Z_i]_\sim\bigr\rangle_{\Lt P}
     =
     \sum_{i\in I} \langle h, \xi_i\rangle_{\Lt P} \sqrt{\mu_i} e_i(t) \, .
  \end{displaymath}
  Since $( \langle h, \xi_i\rangle_{\Lt P} )_{i\in I} \in \ell_2(I)$ and \eqref{diag-sum},  the latter series converges for all $t\in T$, and 
  therefore, the map  $t\mapsto \langle h, \tilde X\rangle_{\Lt P}$ is measurable. By Petti's measurability theorem, see e.g.~\cite[p.~9]{Dinculeanu00}
  or \cite[p.~42]{DiUh77}, we conclude that $\tilde X$ is indeed $\Lt P$-measurable. By \cite[Proposition 11 on p.~6]{Dinculeanu00}
  there then exists a sequence $(\tilde X_n)_{n\geq 1}$ of $\Lt P$-measurable functions $T\to \Lt P$, which are of the form
  \begin{equation}\label{apprx-seq-form}
     \tilde X_n = \sum_{m=1}^\infty \eins_{A_{n,m}} [h_{n,m}]_\sim
  \end{equation}
  for suitable $h_{n,m} \in \sLx 2 P$, and, for $n\geq 1$, mutually disjoint $A_{n,1}, A_{n,2}, \dots\in  \ca B$, such that 
  \begin{equation}\label{uni-conv}
     \sup_{t\in T} \snorm{\tilde X_n(t) - \tilde X(t)}_{\Lt P} \to 0.
  \end{equation}
  Let us write $\hat X_n(\om , t) :=  \sum_{m=1}^\infty \eins_{A_{n,m}}(t) h_{n,m}(\om )$, 
  so that we have $[\hat X(\cdot, t)]_\sim = \tilde X_n(t)$ for all $t\in T$.
  Clearly, each $\hat X_n$ is ($\ca A\otimes \ca B$)-measurable. Our next goal is 
  to show that there exists a subsequence $(\hat X_{n_l})_{l\geq 1}$ such that, for all $t\in T$, 
  there exists an $\hat N_t\in \ca A$ with $P(\hat N_t)=0$ and  
  \begin{equation}\label{almost-conv}
      \hat X_{n_l} (\om , t)\to X_t(\om) \, , \qquad \qquad \om \in \Om\setminus \hat N_t\, .
  \end{equation}
    To this end, we first observe that by \eqref{uni-conv},  for all $l\geq 1$, there exists an $n_l\geq 1$ such that for all 
    $n\geq n_l$ we have 
    \begin{displaymath}
       \sup_{t\in T} \snorm{\tilde X_n(t) - \tilde X(t)}_{\Lt P}^2 \leq 2^{-l}. 
    \end{displaymath}
    Let us fix a $t\in T$. By Markov's inequality we then obtain 
    \begin{displaymath}
       P\bigl( \{ \om \in \Om:    | \hat X_{n_l}(\om,t) - X_t( \om)| \geq l^{-1} \}  \bigr) 
       \leq l^2 \snorm{\tilde X_{n_l}(t) - \tilde X(t)}_{\Lt P}^2 
       \leq l^2 2^{-l}
    \end{displaymath}
    for all $l\geq 1$. A standard application of the Borel-Cantelli lemma then gives the desired $P$-zero sets $\hat N_t\in \ca A$ for which 
    \eqref{almost-conv} holds. 
    
    Let us now write 
    \begin{displaymath}
       D:= \bigl\{ (\om,t): \exists \lim_{l\to \infty} \hat X_{n_l}(\om,t)   \bigr\}\, .
    \end{displaymath}
    Clearly, $D$ is ($\ca A\otimes \ca B$)-measurable, 
    and therefore, there exists an ($\ca A\otimes \ca B$)-measurable function $\hat X:\Om\times T\to \R$
    such that $\hat X(\om , t) = \lim_{l\to \infty} \hat X_{n_l}(\om,t)$ for all $(\om,t)\in D$.
   Moreover, 
  we have 
    $N_t :=\{\om : (\om,t)\not\in D\}\in \ca A$ for all $t\in T$, and 
     our construction ensures $N_t\subset \hat N_t$ for all $t\in T$.  Therefore \eqref{almost-conv}
    yields 
        \begin{displaymath}
       P\bigl(\{\om: \hat X(\om , t) = X_t(\om)\}\bigr) = 1\, , \qquad \qquad t\in T, 
    \end{displaymath}
    i.e.~$(Y_t)_{t\in T}$ defined by $Y_t(\om) := \hat X(\om,t)$
      is indeed an ($\ca A\otimes \ca B$)-measurable version of $(X_t)_{t\in T}$.
%

 Let us finally verify the remaining properties of $(Y_t)_{t\in T}$. To this end, we first 
     observe that the already established $X_t\in \sLx 2 P$ implies $Y_t\in \sLx 2 P$
     and the $\sLx 2 P$-convergence in \eqref{def-X-by-Z-h1} yields
  \begin{align*}
   \E_P Y_s Y_t
   =
   \E_P X_s X_t
      =
      \langle X_s, X_t\rangle_{\sLx 2 P}
      = \sum_{i,j\in I} e_i(s) e_j(t) \E_P Z_i Z_j
      = k_T^1(s,t)
   \end{align*}
   for all   $s,t\in T$. 
   Moreover, for finite $J\subset I$ we obtain by 
   Fubini's theorem, that 
   \begin{align*}
     \int_\Om \int_T  \biggl( Y_t(\om) -  \sum_{j\in J} Z_j(\om) e_j(t) \biggr)^2 \dxy {  \nu} t \, \dxy{  P} \om 
     & =  \int_T\int_\Om  \biggl( X_t(\om) -  \sum_{j\in J} Z_j(\om) e_j(t) \biggr)^2 \dP \om \dxy \nu t \\
     &= \int_T \Bigl\langle \sum_{j\in I\setminus J} Z_j  e_j(t), \sum_{j\in I\setminus J} Z_j  e_j(t)\Bigr\rangle_{\sLx 2 P} \dxy \nu t \\
     & = \int_T \sum_{j\in I\setminus J}  \mu_j e_j^2(t)  \dxy \nu t \\
     & = \sum_{j\in I\setminus J} \mu_j\, ,
   \end{align*}
   and hence we conclude both $Y\in \sLx 2 {P\otimes \nu}$  and 
   \begin{displaymath}
      \bigl[Y(\om)\bigr]_\sim = \sum_{i\in I} Z_i(\om) [e_i]_\sim
   \end{displaymath}
   with convergence in $\Ltn$
    for $P$-almost all $\om \in \Om$.
   For these $\om$, we then find  \eqref{def-Zi} since $([e_i]_\sim)_{i\in I}$
   is an ONS in $\Ltn$.

\end{proofof}

\begin{proofof}{Theorem \ref{unique-kle}}
   For $s,t\in T$ the assumed $\Lx 2 P$-convergence in \eqref{KL-path-conv-Lt-again} together with \eqref{orth-z-again} implies
   \begin{displaymath}
      k(s,t) = \E_PX_sX_t = \bigl\langle [X_s]_\sim, [X_t]_\sim\bigr\rangle_{\Lx 2 P}  
      = \sum_{i,j\in I} e_i(t) e_j(s) \bigl\langle [Z_i]_\sim, [Z_j]_\sim\bigr\rangle_{\Lx 2 P}  
      = \sum_{i\in I} \mu_ie_i(s)e_i(t)\, .
   \end{displaymath}
  By Lemma \ref{kernel-instead-of-EV} we conclude that Assumption X is satisfied and that 
  $(e_i)_{i\in I}\subset H$
	and  $(\mu_i)_{i\in I}$ are the families considered in Assumption X. Consequently, they satisfy Assumption K, and 
	repeating the last part of the proof of 
	Theorem \ref{def-X-by-Z-new} with $Y = X$ thus shows that the  $Z_i$'s satisfy \eqref{def-Zi}.
\end{proofof}

\subsection{Proofs Related to Almost Sure Paths in Interpolation Spaces}

\begin{proofof}{Theorem \ref{interpol-new}}
  Let us  begin by some preliminary remarks. To this end, we  define, for all
$i\in I$, random variables $\xi_i:\Om\to \R$ by
      \begin{equation}\label{def-xi}
         \xi_i(\om ) := \mu_i^{(\b-1)/2}Z_i(\om)\, , \qquad\qquad \om\in \Om.
      \end{equation}
      This definition immediately yields $Z_i(\om)[e_i]_\sim = \xi_i(\om )\mu_i^{(1-\b)/2} [e_i]_\sim$ 
      for all $\om  \in \Om$.
   
   Let us begin by proving \eqref{in-interpol-char-lem-h1-new}. To this end, we simply note that 
   the definition
    of the norm of $[H]_\sim^{1-\b}$ gives
    \begin{align*}
       \bnorm{\sum_{j\in J} Z_j(\om) [e_j]_\sim}_{[H]_\sim^{1-\b}}^2 
      = \bnorm{\sum_{j\in J} \xi_i(\om )\mu_i^{(1-\b)/2} [e_j]_\sim}_{[H]_\sim^{1-\b}}^2 
      =\sum_{j\in J}  \xi_j^2(\om)
      = \sum_{j\in J} \mu_i^{\b-1} Z_j^2(\om)\, ,
    \end{align*}    
    which shows the assertion.  
    
    \aeqb i {ii}
    This immediately follows from \eqref{in-interpol-char-lem-h1-new}, the definition of $[H]_\sim^{1-\b}$, and the equality 
    $[X(\om)]_\sim = \sum_{i\in I} Z_i(\om) [e_i]_\sim$.
    
    \aeqb {ii} {iii} This is a trivial consequence of \eqref{interpol-repres}.
    
    Let us now fix an $\om \in \Om\setminus N$ for which we have $\sum_{i\in I}\mu_i^{\b-1}Z_i^2(\om)<\infty$.
    For an arbitrary $J\subset I$, we then  have $\sum_{j\in J} \mu_j^{\b-1} Z_j^2(\om)<\infty$, and hence we find
    $\sum_{j\in J} Z_j(\om) [e_j]_\sim \in [H]_\sim^{1-\b}$ by using the fact that 
    $(\mu_j^{(\b-1)/2} Z_j(\om))_{j\in J}$ is the sequence of Fourier coefficients of $\sum_{j\in J} Z_j(\om) [e_j]_\sim$
    in $[H]_\sim^{1-\b}$.
     The definition of the norm of $[H]_\sim^{1-\b}$
    then yields \eqref{in-interpol-char-lem-h1-new}. Finally, the unconditional convergence is a direct
    consequence of \eqref{in-interpol-char-lem-h1-new} and the fact that $[H]_\sim^{1-\b}$ and 
    $[\Lx 2 \nu , [H]_\sim]_{1-\b,2}$ have equivalent norms.
\end{proofof}

\begin{proofof}{Theorem \ref{path-space}}
\atob i {ii}  By our assumptions, Lemma \ref{prop-Zi}, and Beppo Levi's theorem we obtain 
\begin{equation}\label{basic-Zi-argument}
   \E_P \sum_{i\in I} \mu_i^{\b-1}Z_i^2 
   =  \sum_{i\in I} \mu_i^{\b-1}\E_PZ_i^2 
   = \sum_{i\in I} \mu_i^\b < \infty\, .
\end{equation}
Consequently, there exists a measurable $\tilde N\subset \Om$ with $P(\tilde N)=0$ such that for all 
 $\om\in \Om\setminus \tilde N$ we have 
$\sum_{i\in I} \mu_i^{\b-1}Z_i^2 (\om) < \infty$.
By Theorem \ref{interpol-new}, we then obtain 
\begin{displaymath}
   [X(\om)]_\sim \in \eHnx{1-\b} = [\Lx 2 \nu , [H]_\sim]_{1-\b,2}
\end{displaymath}
for all $w\in \Om\setminus (N\cup \tilde N)$, which shows the first assertion.
Moreover, choosing $J:=I$ in \eqref{in-interpol-char-lem-h1-new}, we find 
\begin{equation}\label{path-space-h4}
\int_\Om\bnorm{[X(\om)]_\sim}_{[H]_\sim^{1-\b}}^2 \, dP(\om) 
= \int_\Om \sum_{i\in I}  \mu_i^{\b-1}Z_i^2(\om)\, dP(\om) 
= \sum_{i\in I} \mu_i^\b
< \infty\, ,
\end{equation}
where we note that measurability is not an issue as the right-hand side of 
\eqref{in-interpol-char-lem-h1-new} is measurable.
Since the norms of $[\Lx 2 \nu , [H]_\sim]_{1-\b,2}$ and $[H]_\sim^{1-\b}$
are equivalent as discussed around \eqref{interpol-repres}, it thus remains to show that  
 the map $\Om\setminus N\to [H]_\sim^{1-\b}$
      defined by $\om\mapsto [X(\om)]_\sim$ is Borel measurable. To this end, 
      we consider the map $\xi:\Om\setminus  (N\cup \tilde N)\to \ell_2(I)$ defined by
      \begin{displaymath}
         \xi(\om) :=\bigl(\mu_i^{\b-1}Z_i^2 (\om)\bigr)_{i\in I}
      \end{displaymath}
      for all $\om\in \Om\setminus  (N\cup \tilde N)$. Note that our previous considerations
      showed that  $\xi$ indeed maps into $\ell_2(I)$. Consequently,
      $\langle a, \xi\rangle_{\ell_2(I)}:\Om\setminus (N\cup \tilde N)\to \R$
 is well-defined for all $a\in \ell_2(I)$. In addition, this map is clearly measurable, and since $\ell_2(I)$ is 
 separable, the combination of Petti's measurability theorem, cf.~\cite[p.~9]{Dinculeanu00}, with 
 \cite[Theorem 8 on p.~8]{Dinculeanu00} shows that $\xi$ is Borel measurable. Using the isometric relation
 \eqref{Hnk-norm} we conclude that the map $\Om\setminus (N\cup \tilde N)\to [H]_\sim^{1-\b}$ defined by
 \begin{displaymath}
    \om \mapsto \sum_{i\in I} \xi_i(\om) \mu_i^{(1-\b)/2} [e_i]_\sim = [X(\om)]_\sim
 \end{displaymath}
is Borel measurable. 

\atob {ii} i Let $N\subset \Om$ be a $P$-zero set with $[X(\om)]_\sim \in  [\Lx 2 \nu , [H]_\sim]_{1-\b,2}$ for 
all  $\om \in \Om\setminus N$.  
By Theorem \ref{pointwise-KL-new}
we may again assume without loss of generality  that \eqref{KL-path-conv-Lt} is also satisfied 
for all  $\om \in \Om\setminus N$. Using Beppo Levi's theorem and the discussion around \eqref{interpol-repres},
as well as  Lemma \ref{prop-Zi} and \eqref{in-interpol-char-lem-h1-new}, we then obtain 
\begin{displaymath}
 \sum_{i\in I} \mu_i^\b
 = \E_P \sum_{i\in I} \mu_i^{\b-1}Z_i^2 
=\int_\Om\bnorm{[X(\om)]_\sim}_{[H]_\sim^{1-\b}}^2 \, dP(\om) 
< \infty\, .
\end{displaymath}

Let us finally assume that \emph{i)} and \emph{ii)} are true. By Theorem \ref{pointwise-KL-new} there then
exists a measurable $N\subset \Om$ with $P(N)=0$ such that 
$\sum_{i\in I} Z_i(\om) [e_i]_\sim = [X(\om)]_\sim $ in $\Ltn$,  and
$[X(\om)]_\sim \in  [\Lx 2 \nu , [H]_\sim]_{1-\b,2}$ for all $\om\in \Om\setminus N$.
For these $\om$, Theorem \ref{interpol-new}
 immediately yields  
 \begin{equation}\label{path-space-hxxx}
    \sum_{i\in I}\mu^{\b-1}Z_i^2(\om)<\infty\, .
 \end{equation}
Now, to show  
the stronger  $[\Lx 2 \nu , [H]_\sim]_{1-\b,2}$-convergence in \eqref{KL-path-conv-Lt} 
we observe that for all $J \subset I$ and for all $\om \in \Om\setminus N$
we have \eqref{in-interpol-char-lem-h1-new} by \eqref{path-space-hxxx}.
By \eqref{path-space-hxxx} and  \eqref{in-interpol-char-lem-h1-new} we then conclude that the sequence of partial 
sums of $\sum_{i\in I} Z_i(\om) [e_i]_\sim$ is a Cauchy sequence in $[H]_\sim^{1-\b}$ and thus convergent
in $[H]_\sim^{1-\b}$. Moreover, since $[H]_\sim^{1-\b}\hookrightarrow\Ltn$ and
$\sum_{i\in I} Z_i(\om) [e_i]_\sim = [X(\om)]_\sim $ in $\Ltn$, 
its limit is $[X(\om)]_\sim$, which shows 
the $[H]_\sim^{1-\b}$-convergence in \eqref{KL-path-conv-Lt}. Finally, because of \eqref{path-space-hxxx},
the formula \eqref{KL-path-conv-Lt} 
equals the ONB representation of $[X(\om)]_\sim$ with respect to the 
ONB 
$(\mu_i^{(1-\b)/2} [e_i]_\sim)_{i\in I}$ of $[H]_\sim^{1-\b}$, and hence the convergence is also unconditionally.
Now using that $[H]_\sim^{1-\b}$ and $[\Lx 2 \nu , [H]_\sim]_{1-\b,2}$ have equivalent norms, we 
see that the convergence in \eqref{KL-path-conv-Lt}  is indeed 
unconditionally in 
$\bigl[\Lx 2 \nu , [H]_\sim\bigr]_{1-\b,2}$.

To show the last assertion, we combine  \eqref{in-interpol-char-lem-h1-new} with
the just established $[H]_\sim^{1-\b}$-convergence in \eqref{KL-path-conv-Lt} and
 a calculation that is  analogous to 
\eqref{path-space-h4}  to obtain
\begin{displaymath}
 \int_\Om\bnorm{[X(\om)]_\sim - \sum_{j\in J} Z_j(\om) [e_j]_\sim}_{[H]_\sim^{1-\b}}^2  dP(\om) 
=
\sum_{i\in I\setminus J} \mu_j^\b \, .
\end{displaymath}
Again, using that $[H]_\sim^{1-\b}$ and $[\Lx 2 \nu , [H]_\sim]_{1-\b,2}$ have equivalent norms, we then obtain 
 the assertion.
\end{proofof}

\begin{lemma}\label{sums-for-martingales}
   Let $(\xi)_{i\geq 1}$ be a sequence of $\R$-valued random variables on some probability space $(\Om,\ca A,P)$
   and $(\mu_i)_{i\geq 1}\subset (0,\infty)$ be a monotonically decreasing sequence.
   We define $\ca F_i:= \s(\xi_1^2,\dots,\xi_i^2)$ and assume that $\E_P\xi_1^2=1$ and both $\xi_i\in \sLx 4P$ and 
   \begin{equation}\label{martingale-h1}
      \E_P (\xi_{i+1}^2|\ca F_i) = 1
   \end{equation}
    for all $i\geq 1$. Furthermore, assume that, for some $\b\in (0,1)$, we have
     \begin{equation}\label{var-bounds}
       \sum_{i=1}^\infty \mu_i^{2\b}\var \x_i^2 <\infty\, .
     \end{equation}
   Then, the following statements are equivalent:
  \begin{enumerate}
     \item We have $\sum_{i=1}^\infty \mu_i^\b< \infty$.
     \item There exists an $N\in \ca A$ with $P(N)=0$ such that for all $\om \in \Om\setminus N$
     we have 
     \begin{equation}\label{almost-sure-conv}
	    \sum_{i=1}^\infty \mu_i^\b \xi_i^2(\om) <\infty\, .
     \end{equation}
  \end{enumerate}
\end{lemma}

\begin{proofof}{Lemma \ref{sums-for-martingales}}
  Before we begin with the actual proof we 
      note that, for all $i\geq 1$, we have $\E_P \xi_{i+1}^2 = \E_P\E_P (\xi_{i+1}^2|\ca F_i)=1$ by 
      \eqref{martingale-h1}.
      Moreover, for $i>j+1$ an elementary calculation shows
      \begin{equation}\label{sums-for-martingales-h1}
	  \E_P(\xi_i^2|\ca F_j) = \E_P\bigl( \E_P(\xi_i^2|\ca F_{i-1})   |\ca F_j \bigr) = 1\, ,
      \end{equation}
    and by 
      \eqref{martingale-h1} we thus have $\E_P(\xi_i^2|\ca F_j)=1$ for all $i>j$.

     \atob i {ii} This simply follows from 
     \begin{displaymath}
        \E_P \sum_{i=1}^\infty \mu_i^\b \xi_i^2 
        =
          \sum_{i=1}^\infty \mu_i^\b \E_P\xi_i^2 
         =  \sum_{i=1}^\infty \mu_i^\b < \infty\, .
     \end{displaymath}
    \atob {ii} i For $i,n\geq 1$, we write $X_i := \mu_i^\b(\xi_i^2-1)$ and $Y_n:=\sum_{i=1}^n X_i$.
    Then, our first simple observation is that, for $i>j$, we have 
   \begin{equation}\label{sums-for-martingales-h2}
      \E_P (X_i|\ca F_j) = \mu_i^\b  \E_P(\xi_i^2-1|\ca F_j) = 0
   \end{equation}
    by our preliminary considerations. Moreover,  for all $n\geq 1$, the random variable 
    $Y_n$ is $\ca F_n$-measurable and satisfies
    $Y_n\in \sLx 2 P$. In addition, we have 
    \begin{displaymath}
       \E_P(Y_{n+1}|\ca F_n)
       =
       \E_P(X_{n+1}|\ca F_n) + Y_n
       = 
       Y_n
    \end{displaymath}
      by \eqref{sums-for-martingales-h2}, and thus $(Y_n)_{n\geq 1}$ is a martingale with 
      respect to the filtration $(\ca F_n)_{n\geq 1}$. Our next goal is to show that it is uniformly bounded in 
      $\sLt P$. To this end, we first observe that for $i>j$ we have 
      \begin{displaymath}
         \E_P (X_iX_j) = \E_P \E_P(X_iX_j|\ca F_j)
         =
         \E_P\bigl( X_j  \E_P(X_i|\ca F_j) \bigr)  
         =0
      \end{displaymath}
      since $X_j$ is $\ca F_j$-measurable and \eqref{sums-for-martingales-h2}.
      Consequently, we obtain
      \begin{align*}
         \E_P Y_n^2 
         = 
         \sum_{i=1}^n \E_P X_i^2 + 2 \sum_{i=1}^n\sum_{j=1}^{i-1} \E_P(X_iX_j)
          = \sum_{i=1}^n \mu_i^{2\b} \E_P(\xi_i^2-1)^2 
          \leq \sum_{i=1}^\infty \mu_i^{2\b} \var \xi_i^2 \, ,
      \end{align*}
      which by \eqref{var-bounds} shows that  $(Y_n)_{n\geq 1}$ is indeed
      uniformly bounded in $\sLx 2P$. By martingale convergence, see e.g.~\cite[Theorem 11.10]{Klenke08}, 
      there thus exists a random variable
      $Y_\infty\in \sLx 2P$ such that $Y_n\to Y_\infty$ in $\sLx 2P$ and $P$-almost surely.
      In particular, there exists an $\om \in \Om$ with $Y_\infty(\om)\in \R$
      such that we have both \eqref{almost-sure-conv} and  $Y_n(\om)\to Y_\infty(\om)$,
      where the latter simply means that $\sum_{i=1}^\infty X_i(\om)$ converges.
      For this $\om$, we thus obtain
      \begin{align*}
         \sum_{i=1}^\infty \mu_i^\b  
         = \sum_{i=1}^\infty \mu_i^\b  \bigl(\xi_i^2(\om) - \xi_i^2(\om) +1\bigr)
         = \sum_{i=1}^\infty \mu_i^\b  \xi_i^2(\om)   -  \sum_{i=1}^\infty \mu_i^\b  \bigl(\xi_i^2(\om) - 1\bigr)
         = \sum_{i=1}^\infty \mu_i^\b  \xi_i^2(\om) - Y_\infty(\om)\, ,
      \end{align*}
  and since the last difference is a real number we have proven the assertion.
\end{proofof}

\begin{proofof}{Lemma \ref{sum-equiv}}
   \atob i {ii} Follows from a literal repetition of \eqref{basic-Zi-argument}.
   
   \atob {ii} i Our first goal is to show that the random variables $\xi_i := \mu_i^{-1/2} Z_i$
   satisfy the  assumptions of Lemma \ref{sums-for-martingales}. Indeed, we clearly, have 
   $\xi_i\in \sLx 4P$ and the definition of the $\s$-algebras $\ca F_i$ is 
   consistent with Lemma \ref{sums-for-martingales}.
   Moreover, \eqref{Z-martingale} implies \eqref{martingale-h1}, 
   and, for all $\b\in (0,1)$,  condition \eqref{fourier-coeff-cond-interpol-new} implies \eqref{almost-sure-conv}.
   Furthermore, 
    our definitions yields
   \begin{equation}\label{bound-on-xi}
      \var \xi_i^2 = \mu_i^{-2} \var Z_i^2 \leq c \mu_i^{-\a} 
   \end{equation}
  for all $i\geq 1$, and  consequently, we find 
  \begin{displaymath}
      \sum_{i=1}^\infty \mu_i^{2\b}\var \x_i^2
      \leq
      c \sum_{i=1}^\infty \mu_i^{2\b-\a} 
      < \infty
  \end{displaymath}
  whenever $2\b\geq \a+1$, i.e.~\eqref{var-bounds} is satisfied for such $\b$.
  Using Lemma \ref{sums-for-martingales}, we then see that the implication  
  \emph{ii)} $\Rightarrow$ \emph{i)} is true for all 
  $\b\in [\b_1, 1)$, where $\b_1:= (\a+1)/2$. To treat the case $\b\in (\a,\b_1)$, we
   define a sequence $(\b_n)_{n\geq 1}$ by $\b_{n+1} := (\a+\b_n)/2$ for all $n\geq 1$.
  By induction and the definition of $\b_1$, we then see that 
  \begin{displaymath}
     \b_n = 2^{-n} + \a \sum_{i=1}^n 2^{-i}
  \end{displaymath}
  for all $n\geq 1$. Consequently, we have both $\b_n\in (\a,1)$ for all $n\geq 1$
  and $\b_n\searrow \a$.
  
  Our next goal is to show that the implication  
  \emph{ii)} $\Rightarrow$ \emph{i)} is true for all $\b_n$.
  To this end, we first observe that we have already seen that the implication is true for $\b_1$.
  To proceed by induction, we now assume that the implication is true for $\b_n$, so that our goal is to 
  show that it is also true for $\b_{n+1}$. To this end, let us assume that 
  there exists a measurable $N\subset \Om$  with $P(N)=0$ such that 
     \eqref{fourier-coeff-cond-interpol-new}, and thus \eqref{almost-sure-conv}, 
     holds for $\b_{n+1}$ and all $\om \in \Om\setminus N$.
     Here we note that in the absence of such an $N$ there is nothing to prove.
  Now, since $\mu_i\to 0$ and $\b_n>\b_{n+1}$, 
  it is easy to see that  \eqref{fourier-coeff-cond-interpol-new} also holds
  for $\b_{n}$ and all $\om \in \Om\setminus N$, and hence our induction hypothesis yields 
  $\sum_{i=1}^\infty\mu_i^{\b_n}<\infty$. This in turn shows 
  \begin{equation}\label{var-estimate}
     \sum_{i=1}^\infty \mu_i^{2\b_{n+1}}\var \x_i^2
     =
      \sum_{i=1}^\infty \mu_i^{\a+\b_{n}}\var \x_i^2
      \leq 
      c  \sum_{i=1}^\infty \mu_i^{\a+\b_{n}} \mu_i^{-\a} 
      < \infty
  \end{equation}
  by \eqref{bound-on-xi}.
  Consequently, applying Lemma \ref{sums-for-martingales} gives  $\sum_{i=1}^\infty\mu_i^{\b_{n+1}}<\infty$,
  which finishes the induction.
  
  Finally, let us fix a  $\b\in (\a,\b_1)$ for which there exists a measurable $N\subset \Om$  with $P(N)=0$ such that 
     \eqref{fourier-coeff-cond-interpol-new} holds for $\b$ and all $\om \in \Om\setminus N$.
   By the construction of $(\b_n)$, there then exists an $n\geq 1$ such that  $\b\in [\b_{n+1},\b_n)$.   
    Using the same arguments as above, we then see that \eqref{fourier-coeff-cond-interpol-new} also holds
  for $\b_{n}$ and all $\om \in \Om\setminus N$, and hence we find $\sum_{i=1}^\infty\mu_i^{\b_n}<\infty$
  by our preliminary result. Repeating \eqref{var-estimate}, we find 
  \begin{displaymath}
     \sum_{i=1}^\infty \mu_i^{2\b}\var \x_i^2 
     \leq \sum_{i=1}^\infty \mu_i^{2\b_{n+1}}\var \x_i^2 
     \leq  \sum_{i=1}^\infty \mu_i^{\a+\b_{n}} \mu_i^{-\a} 
     <\infty\, ,
  \end{displaymath}
  and consequently  Lemma \ref{sums-for-martingales} gives  $\sum_{i=1}^\infty\mu_i^{\b}<\infty$.
\end{proofof}

\begin{proofof}{Corollary \ref{GP-in-interpol}}
	Clearly, if $I$ is finite, there is nothing to prove, and hence we solely focus on the case $I=\N$. 
	
	\aeqb i {ii}
	By Lemma \ref{gauss-lemma} we know that the $(Z_i)_{i\in I}$ are independent, and thus we find
	$ \E_P (Z_{i+1}^2|\ca F_i) = \E_P Z_{i+1}^2 = \mu_{i+1}$ by Lemma \ref{prop-Zi}.
	Consequently, \eqref{Z-martingale} is satisfied. Moreover,
	since we have $Z_i \sim \ca N(0,\mu_i)$ for all $i\in I$ by  Lemma \ref{gauss-lemma} 
	there exists a constant $c>0$ such that 
	\begin{displaymath}
	 \mu^{-2} \var Z_i^2= \var (\mu_i^{-1/2}Z_i)^2 \leq c
	\end{displaymath}
	for all $i\in I$. This shows that \eqref{Z-martingale-var-bound} holds for all $\a\in (0,1)$.
	Applying Lemma \ref{sum-equiv} then yields the assertion.
	
	\atob {ii} {iii} trivial.
	
	\atob {iii} {ii} Assume that there exists an  $A\in \ca A$  with $P(A)>0$ such that 
     $[X(\om)]_\sim \in [\Lx 2 \nu , [H]_\sim]_{1-\b,2}$ holds for all $\om \in A$. 
     Without loss of generality we may additionally assume that $A\subset \Om\setminus N$, where 
     $N\subset \Om$ is the measurable $P$-zero set obtained from 
   Theorem \ref{pointwise-KL-new}.
     By 
     Theorem \ref{interpol-new} we then know that 
     $\sum_{i\in I}\mu^{\b-1}Z_i^2(\om)<\infty$ for all $\om \in A$, and hence 
     \begin{displaymath}
      P\Bigl(\Bigl\{ \sum_{i\in I}\mu^{\b-1}Z_i^2<\infty  \Bigr\}\Bigr) > 0\, .
     \end{displaymath}
     However, the $(Z_i)_{i\in I}$ are independent by Lemma \ref{gauss-lemma} and hence we conclude 
     by Kolmogorov's 
     zero-one law  that $ \sum_{i\in I}\mu^{\b-1}Z_i^2(\om)<\infty$ 
     actually holds for $P$-almost all $\om\in \Om$.
\end{proofof}

\begin{proofof}{Corollary \ref{sobolev-in-interpol}}
   Let us write $I$ for the embedding $H\hookrightarrow W^m(T)$. 
   Using \eqref{sob-entrop} and the multiplicativity of the dyadic entropy numbers, see \cite[p.~21]{CaSt90},
   we then find
   \begin{displaymath}
      \e_i\bigl(I_k:H\to \Ltn\bigr) \leq \snorm I\cdot \e_i\bigl(\id :W^m(T)\to \Ltn\bigr) \leq c\,  i^{-m/d}\, ,
   \end{displaymath}
   where $c>0$ is a suitable constant.
   Lemma \ref{entropy-versus-ew} then gives $\mu_i \leq 4 c\,  i^{-2m/d}$ for all $i\geq 1$, and 
   hence we have $\sum_{i\in I}\mu_i^\b<\infty$ for all $\b>\frac d {2m}$. 
   Let us fix an $0<s<m-d/2$. For $\b:= 1-s/m$, we then have  $\b\in (\frac d {2m},1)$, and by 
    Theorem \ref{path-space}
   we conclude  that 
   \begin{displaymath}
       [X(\om)]_\sim 
       \in [\Lx 2 T , [H]_\sim]_{1-\b,2} 
       \subset \bigl[\Lx 2 T, W^m(T)  \bigr]_{1-\b,2}
       = B^{(1-\b)m}_{2,2}(T)
       = B^s_{2,2}(T)
   \end{displaymath}
   for $P$-almost all $\om \in \Om$.
  Moreover, the first
  norm estimate, including the implicitly assumed measurability of the integrand, also follows 
  from Theorem \ref{path-space}. The second norm estimate follows by combining 
  Theorem \ref{path-space} with \eqref{sob-entrop} and Lemma \ref{entropy-versus-ew}, which is 
  possible by the assumed $H=W^m(T)$.

   Finally, let us assume that $(X_t)_{t\in T}$  is a Gaussian process with $H=W^m(T)$ but
   \eqref{sobolev-in-interpol-h1} \emph{does}  hold   for $s:= m-d/2$ with strictly positive probability $P$.
   Then we have 
   \begin{displaymath}
       [X(\om)]_\sim  
       \in B^s_{2,2}(T) 
       =  \bigl[\Lx 2 T, W^m(T)  \bigr]_{s/m,2}
       = \bigl[\Lx 2 T, W^m(T)  \bigr]_{1- \b,2}\, ,
   \end{displaymath}
  where $\b:= \frac d {2m}$.
   By Corollary \ref{GP-in-interpol} we then see that $\sum_{i\in I}\mu_i^\b<\infty$, 
   and thus 
   \begin{displaymath}
      \sum_{i\in I}\e_i^{d/m}\bigr(\id : W^m(T)\to \Lx 2 T\bigr)
      =  \sum_{i\in I}\e_i^{2\b}\bigr(I_k : H\to \Lx 2 T\bigr)
      <\infty
   \end{displaymath}
  by Lemma \ref{entropy-versus-ew}. However, this contradicts  \eqref{sob-entrop}.
\end{proofof}

\begin{proofof}{Corollary \ref{cor:small-ball-1}}
 Let us write  $\xi_i := \mu_i^{-1/2}Z_i$ for all $i\in I$. Then we have already seen in 
 front of Corollary \ref{cor:small-ball-1} that 
 $(\xi_i)_{i\in I}$ are i.i.d.~with 
$\xi_i \sim \ca N(0,1)$. Moreover, 
\eqref{in-interpol-char-lem-h1-new} gives 
\begin{displaymath}
 \snorm {[X(\om]_\sim}_{[H]_\sim^{1-\b}}^2 = \sum_{i\in I} \mu_i^{\b-1}Z_i^2(\om) 
 = \sum_{i\in I} \mu_i^\b\, \xi_i^2(\om)
 = \mnorm{(\mu_i^{\b/2}\xi_i(\om))_{i\in I}}_{\ell_2}^2
\end{displaymath}
for $P$-almost all $\om \in \Om$.
Setting $\mu:= \a\b/2$, $p:=2$, and $\s_i:= \mu_i^{\b/2}$, we then obtain the first assertion by
\cite[Theorem 1.1]{Aurzada07a}.
To show the second assertion, we first observe that there are constants $c_1$ and $c_2$ such that
\begin{displaymath}
c_2  \mnorm{(i^{-\a\b/2}\xi_i(\om))_{i\in I}}_{\ell_2}\leq
\mnorm{[X(\om]_\sim}_{[\Ltn, [H]_\sim]_{1-\b,2}} 
\leq c_2  \mnorm{(i^{-\a\b/2}\xi_i(\om))_{i\in I}}_{\ell_2}
\end{displaymath}
for $P$-almost all $\om \in \Om$. Now the second assertion again follows by \cite[Theorem 1.1]{Aurzada07a}.
The last two assertions can be shown analogously.
\end{proofof}

\begin{proofof}{Corollary \ref{cor:small-ball-2}}
   We have already seen in the proof of Corollary \ref{sobolev-in-interpol} that $\mu_i \preceq i^{-2m/d}$
   and
   $[\Lx 2 T, [H]_\sim]_{1-\b, 2} \subset [\Lx 2 T, W^m(T)]_{1-\b, 2} = B_{2,2}^s(T)$
    for $\b:= 1-s/m$. Since this 
   inclusion is continuous by the assumed $H\hookrightarrow W^m(T)$, we then find the assertion by applying Corollary
   \ref{cor:small-ball-1} for $\a:= 2m/d$.
\end{proofof}

\subsection{Proofs Related to Almost Sure Paths in RKHSs}

\begin{lemma}\label{paths-of-versions}
Let $(\Om, \ca A, P)$ be a probability space,
 $(T, \ca B, \nu)$ be a   measure space, and  $(X_t)_{t\in T} \subset \sLt P$ 
  be a $(\ca A\otimes \ca B)$-measurable stochastic process with $X\in \sLt {P\otimes \nu}$. Then, for every $(\ca A\otimes \ca B)$-measurable version 
 $(Y_t)_{t\in T}$ of  $(X_t)_{t\in T}$,
 we have both $(Y_t)_{t\in T} \subset \sLt P$ and  $Y\in \sLt {P\otimes \nu}$, and, 
 for $P$-almost all $w\in \Om$, we further have 
 \begin{displaymath}
    [Y(\om)]_\sim =  [X(\om)]_\sim\, .
 \end{displaymath}
\end{lemma}

\begin{proofof}{Lemma \ref{paths-of-versions}}
   Since  $(Y_t)_{t\in T} \subset \sLt P$ is a version of  $(X_t)_{t\in T} \subset \sLt P$, we have 
   \begin{displaymath}
      P(Y_t = X_t) =1\, , \qquad \qquad t\in T,
   \end{displaymath}
  and thus we find both $(Y_t)_{t\in T} \subset \sLt P$ and $\snorm{Y_t-X_t}_{\sLt P} = 0$ for all $t\in T$.
  Using the measurability of $Y:\Om\times T\to \R$ and Tonelli's theorem, we thus find  
  \begin{align*}
     \int_P\mnorm{ [Y(\om)]_\sim -  [X(\om)]_\sim}_{\Lt \nu}^2\dP\om
     =
    \int_P \int_T \bigl|Y_t(\om) - X_t(\om) \bigr|^2 \dn t \dP\om  
     = 0\, .
  \end{align*}
  This shows  $[Y(\om)]_\sim =  [X(\om)]_\sim$ for  $P$-almost all $w\in \Om$, and since another application
  of  Tonelli's theorem yields 
  \begin{displaymath}
  \int_{\Om\times T}  \bigl|Y_t(\om) - X_t(\om) \bigr|^2  \dxy{P\otimes \nu}{\om,t}
  =
     \int_T\int_P  \bigl|Y_t(\om) - X_t(\om) \bigr|^2  \dP\om\dn t     
     = 0\, ,
  \end{displaymath}
  we also obtain $Y\in \sLt {P\otimes \nu}$.
\end{proofof}

\begin{proofof}{Theorem \ref{path-in-rkhs-1}}
  \atob i {ii} As in the proof of Theorem \ref{interpol-new},  we define, for all
$i\in I$, random variables $\xi_i:\Om\to \R$ by
      \begin{displaymath}
         \xi_i(\om ) := \mu_i^{(\b-1)/2}Z_i(\om)\, , \qquad\qquad \om\in \Om.
      \end{displaymath}
%
%
    For $t\in S$, we further define $Y_t$ by
    \begin{equation}\label{def-Yt}
       Y_t(\om) := \sum_{i\in I} \xi_i(\om) \mu_i^{(1-\b)/2}e_i(t)\, , \qquad \qquad \om \in \Om\setminus N
    \end{equation}
  and $Y_t(\om):= 0$ otherwise. Here we note that 
  the series \eqref{def-Yt} converges for all $s\in T$ and $\om \in \Om\setminus N$, since
  \eqref{fourier-coeff-cond} ensures $(\xi_i(\om))_{i\in I}\in \ell_2(I)$ for all
  $\om \in \Om\setminus N$, while \eqref{k-series-diag-1-b} ensures $(\mu_i^{(1-\b)/2}e_i(t))_{i\in I}\in \ell_2(I)$ for all $t\in S$.
  Finally, for $t\in T\setminus S$ we simply write $Y_t:= X_t$.
  Obviously, this construction guarantees the  $(\ca A\otimes \ca B)$-measurability of  $Y:\Om\times T\to \R$.
  
    Let us first show that $(Y_t)_{t\in T}$ is a version of $(X_t)_{t\in T}$. Clearly, it suffices to show that 
    \begin{displaymath}
       P\bigl(X_t = Y_t \bigr) = 1
    \end{displaymath}
  for all $t\in S$. However, this immediately follows from 
  \begin{displaymath}
     \mnorm{X_t-Y_t}_{\sLt P}^2 = \bnorm{X_t - \sum_{i\in I} Z_i e_i(t)}_{\sLt P}^2 = k(t,t) - \sum_{i\in  I} \mu_i e_i^2(t) = 0\, ,
  \end{displaymath}
where we used both \eqref{LP-diff} and \eqref{k-eq-kn}.
  
  Let us now show that all paths of $Y$ restricted to $S$ are contained in $\Hxy S{1-\b}$. 
  Clearly,
   for $\om \in  N$ our definition yields $ Y(\om)_{|S}=0$, and hence there is nothing to prove for such $\om$.
   Moreover, in the case $\om\in\Om\setminus N$, we first observe that 
   the family of functions $((\mu_i^{(1-\b)/2}\hat e_{i})_{|S})_{i\in I}$
   forms an ONB of $\Hxy S{1-\b}$ since 
   the restriction operator 
   \begin{displaymath}
      \cdot_{|S}:\bHxy S{1-\b}\to \Hxy S{1-\b}
   \end{displaymath}
 is a isometric isomorphism by Lemma \ref{basic-isometries}. Using 
   $(\mu_i^{(1-\b)/2}\hat e_{i})_{|S} = (\mu_i^{(1-\b)/2} e_{i})_{|S}$ and 
   $(\xi_i(\om))_{i\in I} \in \ell_2(I)$,
   we then find $Y(\om)_{|S} \in \Hxy S{1-\b}$ by the definition \eqref{def-Yt}
   of the 
   random variables $Y_t$ for $t\in S$.
   
   \atob {ii} i By Lemma \ref{paths-of-versions} we find a measurable $N_1\subset \Om$ with $P(N_1)=0$ such that 
   $Y(\om)_{|S} \in \Hxy S{1-\b}$ and 
   \begin{displaymath}
       [Y(\om)]_\sim =  [X(\om)]_\sim
   \end{displaymath}
  for all $\om\in \Om\setminus N_1$. Let us fix an $\om \in \Om\setminus N_1$. Since $Y(\om)_{|S} \in \Hxy S{1-\b}$ there 
  then exists a sequence $(a_i)_{i\in I}\subset \ell_2(I)$ such that 
  \begin{equation}\label{path-in-rkhs-1-h3}
     Y(\om)_{|S} = \sum_{i\in I} a_i \mu_i^{(1-\b)/2} (e_{i})_{|S}\, ,
  \end{equation}
  where the convergence is in $\Hxy S{1-\b}$. Let us write $\hat Y (\om):= \eins_S Y(\om)$. Then we find 
  $\hat Y(\om) \in \bHxy S{1-\b}$ and 
  \begin{displaymath}
     \hat Y(\om) = \sum_{i\in I} a_i \mu_i^{(1-\b)/2} \hat e_{i}\, ,
  \end{displaymath}
  where the convergence is in $\bHxy S{1-\b}$. Since $\bHxy S{1-\b}$ is compactly embedded into $\Ltn$, 
   the operator $[\mycdot]_\sim: \bHxy S{1-\b}\to \Ltn$ is 
  continuous, which in turn yields
    \begin{displaymath}
   [X(\om)]_\sim = [Y(\om)]_\sim  =  [\hat Y(\om)]_\sim 
    = \sum_{i\in I} a_i \mu_i^{(1-\b)/2}  [\hat e_{i}]_\sim
    = \sum_{i\in I} a_i \mu_i^{(1-\b)/2}  [ e_{i}]_\sim\, ,
  \end{displaymath}
  where the convergence is in $\Ltn$.
  On the other hand, Theorem \ref{pointwise-KL-new} showed that 
there exists a measurable $N_2\subset \Om$ with $P(N_2)=0$ such that 
 for all $\om \in \Om\setminus N_2$ we have 
 \begin{displaymath}
    [X(\om)]_\sim  = \sum_{i\in I} Z_i(\om) [e_i]_\sim\, ,
 \end{displaymath}
	where again the convergence is in $\Ltn$.
	Using that $([e_i]_\sim)$ is an ONS in $\Ltn$, we thus find $ Z_i(\om) = a_i \mu_i^{(1-\b)/2}$
	if $\om  \not\in N_1\cup N_2$. 
	Now \eqref{path-in-rkhs-1-h2} follows from \eqref{path-in-rkhs-1-h3}, and since 
	 $(a_i)_{i\in I}\in \ell_2(I)$ we also  obtain \emph{i)}
	for $N:=N_1\cup N_2$.
\end{proofof}

\begin{proofof}{Theorem \ref{path-in-rkhs-2}}
 \aeqb i {ii} This has already been shown in Lemma \ref{nuclear-dom}.

Before we prove the remaining implications, let us assume that we have 
 an $(\ca A\otimes \ca B)$-measurable version $(Y_t)_{t\in T}$ of $(X_t)_{t\in T}$  
    such that  $Y(\om)_{|S} \in \Hxy S{1-\b}$  for $P$-almost all $\om \in \Om$.
    By Lemma \ref{paths-of-versions}
     we then conclude that 
	\begin{displaymath}
	[\hat Y(\om)_{|S}]_\sim =  [Y(\om)]_\sim = [X(\om)]_\sim
	\end{displaymath}
	for $P$-almost all $\om \in \Om$, where $\hat Y(\om)$ denotes the zero-extension of  $Y(\om)_{|S}$ to $T$.
	In addition, we have $\snorm{Y(\om)_{|S}}_{\Hxy S{1-\b}} = \snorm{[\hat Y(\om)_{|S}]_\sim}_{[H_S^{1-\b}]_\sim}$
	by Lemma \ref{basic-isometries}. Together, this yields
	\begin{equation}\label{XY-norm}
	 \int_\Om \mnorm{Y(\om)_{|S}}_{\Hxy S{1-\b}}^2 \dP\om 
	= \int_\Om \mnorm{[X(\om)]_\sim}_{[H]_\sim^{1-\b}}^2 \dP\om
	= \sum_{i\in I}\mu_i^{\b}
	\end{equation}
	where the last identity follows by 
	a repetition of \eqref{path-space-h4}.
	Moreover, note that all three quantities may simultaneously be infinite.

	\atob i {iii} We have 
	\begin{displaymath}
	 \int_\Om \sum_{i\in I}  \mu_i^{\b-1}Z_i^2(\om) \dP\om 
	= \sum_{i\in I} \mu_i^{\b-1}\int_\Om Z_i^2(\om)\dP\om 
	= \sum_{i\in I}  \mu_i^{\b} <\infty\, ,
	\end{displaymath}
	and hence we find a measurable $N\subset \Om$  with $P(N)=0$ such that for all $\om \in \Om\setminus N$ we have
    \eqref{fourier-coeff-cond}. Now the assertion follows from Theorem \ref{path-in-rkhs-1} and \eqref{XY-norm}.
	
	\atob {iii} i Follows directly from \eqref{XY-norm}.
\end{proofof}

\begin{proofof}{Corollary \ref{GP-in-RKHS}}
   \aeqb {i} {ii} This has already been shown in Lemma \ref{nuclear-dom}, see also Theorem \ref{path-in-rkhs-2}.
   
   \atob i {iii} Repeating \eqref{basic-Zi-argument}, we see yet another time  that \eqref{fourier-coeff-cond}
   holds for $P$-almost all $\om\in \Om$. Applying Theorem \ref{path-in-rkhs-1} then yields the assertion.
   
   
   \atob {iii} {iv} trivial
   
   \atob {iv} i For $\om\in A$ we have $[X(\om)]_\sim = [\hat Y(\om)_{|S}]_\sim \in [\Hxy S{1-\b}]_\sim = [H]_\sim^{1-\b}$
  and hence \emph{i)} follows by   Corollary \ref{GP-in-interpol}.
\end{proofof}

\begin{proofof}{Corollary \ref{path-in-rkhs-cor1}}
  Before we begin with the actual proof, let us first 
     note that the factorization
    
    \tridia H \Ltn {\bar H} {I_k}\id{I_{\bar k}}
    
    \noindent
    together with the multiplicativity of the dyadic entropy numbers, see \cite[p.~21]{CaSt90},
    yields 
    \begin{displaymath}
       \e_i(I_k) \leq \snorm{\id: H\to \bar H} \, \e_i(I_{\bar k})
    \end{displaymath}
 for all $i\geq 1$, and therefore
    we find $\sum_{i=1}^\infty \e_i^{\a}(I_{k})<\infty$. Applying  Lemma 
   \ref{entropy-versus-ew} then 
   shows both  $\sum_{j\in J}\bar \mu_j^{\a/2}<\infty$ and $\sum_{i\in I}\mu_i^{\a/2}<\infty$,
   where $(\bar \mu_j)_{j\in J}$ is the sequence of non-zero eigenvalues of $T_{\bar k}$ obtained by 
   Theorem \ref{spectral}.
   
  Moreover, for $\b\in [\a/2,1-\a/2]$, we have $\a/2\leq 1-\b$, and thus  
  we find both 
   $\sum_{j\in J}\bar \mu_j^{1-\b}<\infty$ and $\sum_{i\in I}\mu_i^{1-\b}<\infty$. Analogously, $\b\geq \a/2$
   implies $\sum_{j\in J}\bar \mu_j^{\b}<\infty$ and $\sum_{i\in I}\mu_i^{\b}<\infty$.

   \ada i
   Let us   pick a $\b\in [\a/2,1-\a/2]$. 
   Then, our preliminary considerations showed both 
   $\sum_{j\in J}\bar \mu_j^{1-\b}<\infty$ and $\sum_{i\in I}\mu_i^{1-\b}<\infty$. By \eqref{when-kS}
   we then see that we find a measurable  $S_0\subset T$
	with $\nu(T\setminus S_0)=0$ such that both $\Hxy{S_0}{1-\b}$ and $\bbHxy{S_0}{1-\b}$ exist.

   Our next goal is to find a subset $S$ of $S_0$ with $\nu(T\setminus S)=0$ and $\Hxy S{1-\b}\subset \bbHxy S{1-\b}$.
   To this end, note that   \eqref{interpol-repres} 
   together with $[H]_\sim \subset [\bar H]_\sim \subset \Ltn$ and the definition of interpolation norms 
   shows
   \begin{displaymath}
       \eHnx {1-\b}
      =\bigl[\Lx 2 \nu , [H]_\sim\bigr]_{1-\b,2}
      \hookrightarrow
       \bigl[\Lx 2 \nu , [\bar H]_\sim\bigr]_{1-\b,2}
       = \ebbHnx {1-\b} \, ,
   \end{displaymath}
	and hence the  inclusion operator $I:\eHnx {1-\b} \to \ebbHnx {1-\b}$ is continuous. Now consider 
	the situation

	\seqdia {\Hxy {S_0}{1-\b}}{\eHnx {1-\b}}{\ebbHnx {1-\b}}{\bbHxy {S_0}{1-\b}} {[\hat \mycdot]_\sim}I{[\hat \mycdot]_\sim}
	
	\noindent
	where the operators ${[\hat \mycdot]_\sim}$ are isometric isomorphisms by 	
	Lemma \ref{basic-isometries}.
    Consequently, for all $f\in \Hxy{S_0}{1-\b}$ there exists a unique $g_f\in \bbHxy{S_0}{1-\b}$ such that 
    $[\hat f]_\sim = [\hat g_f]_\sim$, and the map $f\mapsto g_f$ is linear and continuous.
	In other words, for all $f\in \Hxy{S_0}{1-\b}$, there exists a 
     a measurable $N_f\subset S_0$ with $\nu(N_f)=0$ and 
    $f(t) = g_f(t)$ for all $t\in S_0\setminus N_f$.  

  Let us 
 find such a $\nu$-zero set $N$ that is an independent of $f$. To this end, we fix a countable dense
    $D\subset \Hxy{S_0}{1-\b}$ and define $N:= \bigcup_{f\in D} N_f$, where we note that such a $D$ exists 
    since $\Hxy{S_0}{1-\b}$ is separable by construction. Now the definition of $N$  immediately yields 
    $N\subset S_0$ and $\nu(N) = 0$, as well as
    \begin{equation}\label{equal-on-S0}
       f(t) = g_f(t)\, , \qquad \qquad t\in S_0\setminus N
    \end{equation}
    for all $f\in D$. 
	To show the latter for all $f\in \Hxy{S_0}{1-\b}$, we fix such an $f$ and a sequence 
	$(f_n)\subset D$ with $f_n\to f$ in $\Hxy{S_0}{1-\b}$. Then we have $g_{f_n}\to g_f$ in 
	$\bbHxy{S_0}{1-\b}$ by the above mentioned continuity of $f\mapsto g_f$, 
	and since both spaces are reproducing kernel Hilbert spaces, we obtain 
	$f_n(t) \to f(t)$ and $g_{f_n}(t)\to g_f(t)$ for all $t\in S_0$.
	Using $f_n(t) = g_{f_n}(t)$ for all $t\in S_0\setminus N$ and $n\geq 1$, we thus find
	\eqref{equal-on-S0}. Defining $S:= S_0\setminus N$ then gives $\Hxy S{1-\b}\subset \bbHxy S{1-\b}$
	and the continuity of this embedding follows from the continuity of $I$ and Lemma \ref{basic-isometries}.
	
  \ada {ii} Our goal is to apply Theorem \ref{path-in-rkhs-2}. To this end, we first observe that 
    \eqref{k-eq-kn} holds for a set $\tilde S\subset T$ with $\nu(T\setminus \tilde S)=0$ by the assumed separability of $H$ and \cite[Corollary 3.2]{StSc12a}.
%
%
    Consequently, we may assume without loss of generality that \eqref{k-eq-kn} holds for the set $S$ found in 
    part \emph{i)}.
    Moreover,  we have already seen in part \emph{i)} that we have 
    $\sum_{i\in I}\mu_i^{1-\b}<\infty$, which in turn implies \eqref{k-series-diag-1-b}  by \eqref{when-kS}.
    Finally, our preliminary considerations showed that $\b\geq \a/2$ implies $\sum_{i\in I}\mu_i^{\b}<\infty$,
    and thus Theorem \ref{path-in-rkhs-2} is applicable.
\end{proofof}

\begin{proofof}{Corollary \ref{path-in-rkhs-cor2}}
 We first show that assumption \emph{i)} implies assumption \emph{ii)}, so that in the remainder of this proof
 is suffices to work with the latter. To this end, note that 
 \begin{align*}
   \sum_{j\in J} \bar \mu_j^{1-\b} \bar e_j^2(t)
    \leq \sup_{j\in J}\inorm{\bar e_j}  \sum_{j\in J} \bar \mu_j^{1-\b}
   \leq \sup_{j\in J}\inorm{\bar e_j}  \sum_{j\in J} \bar \mu_j^{\b}
   \leq 4\sup_{j\in J}\inorm{\bar e_j}  \sum_{i=1}^\infty \e_i^{2\b}(I_{\bar k})< \infty\, ,
 \end{align*}
  where we used $1-\b\geq \b$ and Lemma \ref{entropy-versus-ew}. Consequently, $\bar k_T^{1-\b}$ exists and is bounded,
  and from the latter we immediately obtain 
  $[\Ltn , [\bar H]_\sim ]_{1-\b,2}  = [\bar H_{T}^{1-\b}]\hookrightarrow \Lx \infty {\nu}$.

 \ada {i} We first note that $H\subset \bar H$ implies $\t(H)\subset \t(\bar H)$, and hence Assumption CK
   is satisfied for $k$, too. Moreover, the
    continuity of 
    the  inclusion operator $I:\eHnx {1-\b} \to \ebbHnx {1-\b}$
    considered in the proof of part \emph{i)} of Corollary \ref{path-in-rkhs-cor1}
    implies $[\Ltn , [H]_\sim]_{1-\b,2}\hookrightarrow \Lx \infty {\nu}$.
    By Theorem \ref{bounded-powers-of-cont-kernels}, we then see that both $\Hxy T{1-\b}$ and $\bbHxy T{1-\b}$ 
    do exist.
    Moreover, the kernels $k_T^{1-\b}$ and $\bar k_T^{1-\b}$ are bounded by Theorem \ref{bounded-powers-of-cont-kernels}.
    
    To show that $\Hxy T{1-\b}\subset \bbHxy T{1-\b}$, we   consider the map 
    $f\mapsto g_f$ from the proof of part \emph{i)} of Corollary \ref{path-in-rkhs-cor1}.
    Then we have seen above that \eqref{equal-on-S0} holds for $S_0=T$ and all $f\in \Hxy T{1-\b}$.
    Let us assume that there exists an $f\in \Hxy T{1-\b}$ and a $t\in T$ such that $f(t)\neq g_f(t)$.
    Then we have $\{|f-g_f|>0  \}\neq \emptyset$ and $\{|f-g_f|>0  \} \in \t(\bar H)$, which together imply
    $\nu(\{|f-g_f|>0  \})>0$, since $\nu$ is assumed to be $\bar k$-positive. In other words, 
    \eqref{equal-on-S0} does not hold for   $f$, which contradicts our earlier findings. This shows $f=g_f$
    for all $f\in \Hxy T{1-\b}$ and thus  $\Hxy T{1-\b}\subset \bbHxy T{1-\b}$. The continuity of the
    corresponding embedding again follows from the continuity of $I$. 
    
    \ada {ii}   Considering the proof of part
    \emph{ii)} of Corollary \ref{path-in-rkhs-cor1},
    we easily see that it suffices to check that  \eqref{k-eq-kn} holds for $S:=T$.
    The latter, however, follows from   Lemma \ref{Ik-inject}.
    
    \ada {iii} All $f\in \Hxy T{1-\b}$  are bounded since the kernel $k_T^{1-\b}$ is bounded.
    Moreover, all $f\in \Hxy T{1-\b}$ are $\t(\Hxy T{1-\b})$-continuous by the very definition of this topology,
    and since 
    Theorem \ref{bounded-powers-of-cont-kernels} showed $\t(\Hxy T{1-\b}) = \t(H)$, they are also $\t(H)$-continuous.
    Now the additional assertions on the paths of $Y$ follow from $Y(\om)\in \Hxy T{1-\b}$ for all $\om \in \Om$.
    
    \ada {iv} Since $\bar k_T^{1-\b}$ is bounded, we have $\bbHxy T{1-\b}\hookrightarrow \ell_\infty(T)$, see 
    e.g.~\cite[Lemma 4.23]{StCh08}. Now 
    the $\ell_\infty(T)$-convergence of \eqref{path-in-rkhs-1-h2}
    follows from the $\bbHxy T{1-\b}$-convergence established in Theorem \ref{path-in-rkhs-1}.
    
    \ada {v} Let us   fix a countable, $\t$-dense subset $D\subset T$. Since $Y$ is a version of $X$,
    we then have $P(\{Y_t \neq X_t\}) = 0$ for all $t\in D$, and hence there exists a $P$-zero set $N\in \ca A$
    such that $X_t(\om) = Y_t(\om)$ for all $t\in D$ and $\om \in \Om\setminus N$. Without loss of 
    generality we may also assume that $X(\om)$ is $\t$-continuous for all $\om \in \Om\setminus N$.
    and since $\t(H)\subset \t$, we further see by part \emph{iii)} that all paths of $Y$ are $\t$-continuous,
    too. Now the assertion follows by  a simple limit argument.
    
    By Lemma \ref{Ik-inject} the operator $I_{\bar k}$ is injective, and thus \cite[Theorem 3.1]{StSc12a}
    shows that $(\bar e_j)_{j\in J}$ is an ONB of $\bar H$. Consequently, $\bar H$ is separable and 
    Lemma \ref{simple-topol} shows that $\t(\bar H)$ is separable and generated by a pseudo-metric. 
    If $\t(H)$ is Hausdorff, this pseudo-metric becomes a metric and the assertion follows from the first part.
\end{proofof}

\begin{proofof}{Corollary \ref{sobolev-in-rkhs}}
  \ada i Let us consider Corollary \ref{path-in-rkhs-cor1} for $\bar H=W^m(T)$.
  Then \eqref{sob-entrop} 
  shows that 
  \begin{displaymath}
      \sum_{i=1}^\infty \e_i^{\a}(I_{\bar k})<\infty
  \end{displaymath}
  holds for all $\a>d/m$. Let us pick an $s\in (d/2,m-d/2)$ and define $\b:= 1 - s/m$.
  This gives $\frac d{2m} <\b<1-\frac d{2m}$, and hence $\b$ satisfies the assumptions of 
  Corollary \ref{path-in-rkhs-cor1} for a suitable 
  $\a\in (d/m,1]$ with $\b\in [\a/2, 1-\a/2]$.
%
%
  Moreover, we have 
  \begin{displaymath}
     [\Lx 2 T , [H]_\sim]_{1-\b,2}  
      \hookrightarrow \bigl[\Lx 2 T, W^m(T)  \bigr]_{1-\b,2}
       = B^{(1-\b)m}_{2,2}(T)
       = B^s_{2,2}(T)
       \hookrightarrow \Lx\infty \nu
  \end{displaymath}
  by Sobolev's embedding theorem for Besov spaces, see e.g.~\cite[Theorem 7.34]{AdFo03},
  and hence we can apply part \emph{iii)} of  Corollary \ref{path-in-rkhs-cor1} and Theorem \ref{path-in-rkhs-1}.
   
   \ada {ii} This follows from Corollary \ref{sobolev-in-interpol} since 
   \eqref{sobolev-in-rkhs-h1} implies \eqref{sobolev-in-interpol-h1}.
\end{proofof}

{\small
\bibliographystyle{plain}
\bibliography{../../../../literatur-DB/steinwart-mine,../../../../literatur-DB/steinwart-books,../../../../literatur-DB/steinwart-proc,../../../../literatur-DB/steinwart-article}
}

\end{document}